\documentclass[11pt,twoside,english]{elsarticle}
\usepackage{amsfonts}
\usepackage[T1]{fontenc}
\usepackage[latin9]{inputenc}
\usepackage{geometry}
\usepackage{amsmath}
\usepackage{amsthm}
\usepackage{amssymb}
\usepackage{graphicx}
\usepackage{subfig}
\usepackage{babel}
\usepackage{float}
\makeatletter
\providecommand{\tabularnewline}{\\}

\newtheorem{thm}{Theorem}
\newtheorem{lem}{Lemma}
\theoremstyle{definition}
\newtheorem{rem}{Remark}

\ifx\proof\undefined
\newenvironment{proof}[1][\protect\proofname]{\par
\normalfont\topsep6\p@\@plus6\p@\relax
\trivlist
\itemindent\parindent
\item[\hskip\labelsep\scshape #1]\ignorespaces
}{
\endtrivlist\@endpefalse
}
\providecommand{\proofname}{Proof}
\fi
\@ifundefined{showcaptionsetup}{}{ \PassOptionsToPackage{caption=false}{subfig}}
\makeatother

\input{tcilatex}
\begin{document}

\begin{frontmatter}{}

\title{Convexification and experimental data for a 3D inverse scattering problem with the moving point source}

\tnotetext[t1]{The work of Khoa, Bidney, Klibanov, L. H. Nguyen and Astratov was supported by US Army Research Laboratory and US Army
Research Office grant W911NF-19-1-0044. The work of Khoa was also
partly supported by the Research Foundation-Flanders (FWO) under the
project named ``Approximations for forward and inverse reaction-diffusion
problem related to cancer models''.}

\author[rvt]{Vo Anh Khoa}
\ead{vakhoa.hcmus@gmail.com}

\author[rvt1]{Grant W. Bidney}
\ead{gbidney@uncc.edu}

\author[rvt]{Michael V. Klibanov\corref{cor1}}
\ead{mklibanv@uncc.edu}
\cortext[cor1]{Corresponding author.}

\author[rvt]{Loc H. Nguyen}
\ead{loc.nguyen@uncc.edu}

\author[rvt2]{Lam H. Nguyen}
\ead{lam.h.nguyen2.civ@mail.mil}

\author[rvt2]{Anders J. Sullivan}
\ead{anders.j.sullivan.civ@mail.mil}

\author[rvt1]{Vasily N. Astratov}
\ead{astratov@uncc.edu}

\address[rvt]{Department of Mathematics and Statistics, University of North Carolina
at Charlotte, Charlotte, North Carolina 28223, USA.}

\address[rvt1]{Department of Physics and Optical Science, University of North Carolina at Charlotte, Charlotte, NC 28223, USA.}

\address[rvt2]{U.S. Army Research Laboratory, Adelphi, Maryland 20783-1197, USA.}

\begin{abstract}
Inverse scattering problems of the reconstructions of physical properties of 
a medium from boundary measurements
are substantially challenging ones. This work aims to 
verify the performance on experimental data of a newly developed convexification 
method for a 3D coefficient
inverse problem for the case of objects buried in a sandbox a
fixed frequency and the  point source moving along an interval of a straight line. Using
a special Fourier basis, the method of this work strongly relies on a new derivation
of a boundary value problem for a system of coupled quasilinear elliptic equations.
This problem, in turn is solved via the minimization of a Tikhonov-like functional  weighted by a Carleman Weight Function. The global convergence of the numerical procedure is established analytically.
The numerical verification
is performed using experimental data, which are raw backscatter data of the
electric field. These data were collected using a microwave scattering facility
at The University of North Carolina at Charlotte. 
\end{abstract}

\end{frontmatter}{}

\section{Introduction}

Inverse and ill-posed problems are ubiquitous in many branches of physical,
biological, ecological and social sciences. The goal of this paper is to
figure out how to reconstruct physical properties of objects buried closely
under the ground on a small depth not exceeding 10 centimeters. Since these
objects are buried, then they present a significant challenge in detection
and identification of suspicious explosive-like targets (antipersonnel land
mines and improvised explosive devices) in military applications.
Mathematically, the challenge in solving such problems is the limitation of
observable quantities of inputs to the physical-based mathematical systems.
In this 3D scenario, only data on a single surface are physically measured
in a fixed and specific discrete setting. Certainly, it is expensive and
time-consuming to try finer refinements of spatial measurements of
observable quantities.

We call a numerical method for a nonlinear ill-posed problem \emph{globally
convergent} if there is a theorem which guarantees that this method delivers
at least one point in a sufficiently small neighborhood of the exact
solution without any advanced knowledge of this neigborhood. In other words,
this theorem guarantees that although a good first guess for the solution is
not required, still a good approximation for the solution will be obtained
by that numerical method. The convexification method we work here is a
globally convergent one.

\subsection{Scopes and novelty}

The goal here is to solve a 3D Coefficient Inverse Problem (CIP) for the
case of the experimentally collected backscattering data. In doing so, we
want to estimate the dielectric constants and shapes of buried targets. We
consider here the case when the point source moves along an interval of a
straight line and the frequency is fixed. In this work, we apply the globally convergent convexification method, which has been developed by the third
author and his coauthors for a significant number of years. On the other hand, a large amount of previous
papers concerning the convexification with numerical investigations consider the scenario when either the point source or the direction of the
incident plane wave is fixed and the frequency is varied; cf. \cite%
{Klibanov2017a,convIPnew,Klibanov2019,Klibanov2019b,Klibanov2019a,Klibhyp,Truong}. While locations
and dielectric constants of targets are imaged accurately in such cases,
shapes of targets are not accurately imaged; see, e.g., images of \cite%
{Klibanov2019b} for the case of experimentally collected data.

Thus, the idea here is to move the source in a hope that this would provide
better images of shapes of targets of interest while still preserving
accurate values of computed dielectric constants. This novel insight was first
implemented in \cite{Khoa2019}, where it was shown that it works
well for computationally simulated backscattering data: both dielectric
constants and shapes of targets were accurately imaged via a new version of
the convexification method. In this paper, we continue to verify the
computational performance of the idea of \cite{Khoa2019} using experimental
backscattering data. These data are raw backscatter data of the electric field that we collected using a microwave scattering facility at our University. We then show here that the version of \cite{Khoa2019} the convexification method
accurately computes values of dielectric constants as well as shapes of
targets buried in a sandbox (including those with even rather complicated shapes). Notice that since
the frequency is fixed both in \cite{Khoa2019} and here, then the data we
use are non overdetermined ones.

We realize that estimates of dielectric constants, shapes and locations of
buried targets mimicking explosives are insufficient to distinguish between
explosives and non-explosives. Nevertheless, these estimates might serve in
the future as a piece of the information, which is an additional one to the
features currently used for classifications of explosive-like targets. Hence, this new information might help to decrease the false alarm rate.

The convexification method \textquotedblleft convexifies\textquotedblright\
a weighted Tikhonov-like functional via using a suitable Carleman Weight
Function (CWF). While starting from the first inception in 1981 (cf. \cite%
{BukhKlib}) with the only goal at that time of proofs of global uniqueness
theorems for CIPs (also, see, e.g., \cite%
{Beilina2012,Klibanov2013,Yamamoto1999}), the notion of applications of
Carleman estimates to CIPs got a new aspect nowadays in terms of numerical
methods for CIPs. This is because, driven by the CWF, the resulting convexified
cost functionals avoid being trapped in local minima and ravines.

The first publications on the convexification were in 1995 and 1997; cf. \cite%
{Klib95,Klib97}. Also, the reader can be referred to some other initial follow-up results in \cite%
{BK,KT,KK}. However, these past papers were purely theoretical ones. The reason was that some
important mathematical facts leading to numerical implementations were
unknown at that time, although some numerical results were published in \cite{KT}. In
2017, the work \cite{Bakushinsky2017} has \ clarified those facts. This led
to a significant number of recent publications on the convexification, where both analytical
and numerical results were present; cf. \cite%
{Khoa2019,Klibanov2017a,convIPnew,Klibanov2019,Klibanov2019b,Klibanov2019a,Klibhyp,Truong}%
. It is worth mentioning that papers \cite{Klibanov2017a,convIPnew,Klibanov2019b} address the
performance of the convexification on experimental data, which is the same research line of this article. In particular, the
paper \cite{Klibanov2019b} considers experimental data for buried
objects for the case of a single location of the point source and multiple
frequencies.

In our four latest publications \cite%
{Klibanov2019,Klibanov2019b,Klibanov2019a,Klibhyp} the convexification
method is based upon the derivation of boundary value problem with
overdetermined boundary conditions for a system of coupled elliptic PDEs. In
fact, some of those boundary conditions are Cauchy data. Next, an
approximate solution of this system is found via the minimization of a
weighted globally strictly convex Tikhonov-like functional with a CWF in it.
In the case of CIPs for the Helmholtz equation considered in \cite%
{Klibanov2019,Klibanov2019b}, multiple frequencies were used and only a
single location of the point source. 

The boundary value problem mentioned in the previous paragraph is about
finding spatially dependent coefficients of a certain truncated Fourier
series with respect to a new orthonormal basis in the $L^{2}$ space, which
was first introduced in \cite{Klibanov2017}. In \cite{Khoa2019}, the theory
of that version of the convexification was developed for the continuous
case. Unlike \cite{Khoa2019}, we consider here a \textquotedblleft partial
finite difference'' case, which is another novelty of this article. This means that partial derivatives with respect to
two out of three variables are written in finite differences. However, we
do not allow the grid step size to tend to zero. This agrees with the fact that it is certainly expensive in real-world applications to decrease the grid step size indefinitely. We point out that an important advantage of
using partial finite differences is that we do not use the regularization
penalty term in the weighted Tikhonov-like functional. The global
convergence analysis in the framework of partial finite differences is
performed here. We believe that this is a significant analytical novelty of our work.

\subsection{Related works and outline of the paper}

The existing literature on this topic is huge to
be singled out. Conventional numerical approaches to CIPs rely on
the minimization of some least squares Tikhonov functionals; see, e.g., \cite{Chavent,Gonch1,Gonch2}. It is well known, however, that these functionals are non convex
and typically have multiple local minima and ravines. Thus, convergence of a
minimization procedure to the exact solution can be sometimes guaranteed in
such a case only if its starting point is located in a sufficiently small
neighborhood of that solution, i.e. this would be a \emph{locally}
convergent numerical method. Unlike this, the concept of the convexification
leads to globally convergent numerical methods.

We now refer only to the closest publications, since this
paper is not a survey. A version of the convexification, which is different from ours (see above-cited references), has been developed by Baudouin,
de Buhan, Ervedoza and Osses (cf. \cite{Baud1,Baud2}) for two CIPs for the
hyperbolic equations and then for quasilinear parabolic equations (cf. 
Boulakia, de Buhan and Schwindt \cite{BBS} and Le and Nguyen \cite{LN}). In
this version, a CWF is used to construct a sort of a contractual mapping operator. A nice feature of these techniques is that the corresponding numerical method converge globally, which again accentuates the importance of the inclusion of CWFs in numerical schemes. The publications \cite%
{Baud1,Baud2,BBS} work within the framework of the Bukhgeim--Klibanov method 
(cf. \cite{Beilina2012,BukhKlib,KT,Klibanov2013}), which assumes that an initial condition is not vanishing in
the entire domain of interest. On the other hand, all our above cited publications about the
convexification for the CIPs for the wave propagation processes essentially
work with the case when the forward problem for the Helmholtz equation can
be reduced via the Fourier transform to the one for a hyperbolic equation,
in which one of initial conditions is identical to zero and the second one is
the $\delta -$function (also see \cite{Klibhyp} for a similar situation in
the time domain case).

CIPs for the fixed frequency case have been extensively studied by Novikov
and his coauthors since 1988; cf. \cite{N1}. In particular, several reconstruction
methods were proposed by this group and numerical results were also accompanied in 
\cite{Ag,Al,N2,N3}. It is worth noting that these CIPs and reconstruction techniques are different
from the ones we are using. A nice feature of these reconstruction methods
is that they are in the category of globally convergent numerical methods,
since they do not require a good first guess; see the second paragraph of
this section for the definition of the global convergence. An interesting feature of \cite{Ag,N3} is that these publications consider
the case of non overdetermined data for the reconstruction of the potential
of the Schrödinger equation at high values of the wavenumber are considered
in \cite{Ag,N3}. Furthermore, the data in \cite{Ag} are phaseless.
Corresponding numerical results of this series of works are published in 
\cite{Ag,Al}.

There is another feature of the techniques of the Novikov's group, which is
philosophically close to one of features of our technique. Their
reconstruction procedures rely on certain steps on truncations of Fourier
series. On the other hand, we also truncate the Fourier series in this work; see Remark \ref{rem:1}.
And also, neither we nor the group of Novikov cannot prove convergence as
the truncation number $N\rightarrow \infty .$

Recently, Bakushinsky and Leonov have proposed an algorithm for solving a CIP
with multi frequency data; cf. \cite{Bakushinsky2019}. Their method is based on
the solution of an integral equation of the first kind generated by the
fundamental solution of the Helmholtz equation in a homogeneous medium. In
the case of the time dependent data and low contrast targets, we refer to
Goncharsky and Romanov \cite{Gonch1} for computationally simulated data and
to Goncharsky, Romanov and Seryozhnikov \cite{Gonch2} for a quite successful
application of the method of \cite{Gonch1} to experimental data.

This paper is organized as follows. Section \ref{sec:Statement-of-the} is
devoted to the mathematical statement of the forward and inverse problems. In
section \ref{sec:A-globally-convergent}, we present the derivation of the
quasi-linear PDE system and then design a weighted Tikhonov-like functional
that we want to minimize in the numerical process. Thereby, an important part
in this section will be principal proofs of convergence of the minimizer towards the
true solution in a finite difference framework. In section \ref%
{sec:Experimental-results}, we delineate our numerical results with experimental data. Here, the
so-called data propagation procedure is revisited.

\section{Statement of the problem\label{sec:Statement-of-the}}

Let $\mathbf{x}=\left( x,y,z\right) \in \mathbb{R}^{3}$. Consider a
rectangular prism $\Omega =\left( -R,R\right) \times \left( -R,R\right)
\times \left( -b,b\right) $ in $\mathbb{R}^{3}$ for $R,b>0$ as our
computational domain of interest. Let $c\left( \mathbf{x}\right) $ be a
spatially distributed dielectric constant of the medium. We assume that the
function $c\left( \mathbf{x}\right) $ is smooth and also that 
\begin{equation}
\begin{cases}
c\left( \mathbf{x}\right) \geq 1 & \text{in }\mathbb{R}^{3}, \\ 
c\left( \mathbf{x}\right) =1 & \text{in }\mathbb{R}^{3}\backslash \Omega .%
\end{cases}
\label{eq:2}
\end{equation}%
The second line of (\ref{eq:2}) means that we are assuming to have vacuum
outside the domain of interest $\Omega $. Let the number $d>b$ and let $%
a_{1}<a_{2}$. We consider a line of sources 
\begin{equation*}
L_{\text{src}}:=\left\{ \left( \alpha ,0,-d\right) :a_{1}\leq \alpha \leq
a_{2}\right\} ,
\end{equation*}%
which is parallel to the $x$-axis and is located outside of the closed
domain $\overline{\Omega }$. The distance from $L_{\text{src}}$ to the $xy$%
-plane is $d$, and the length of the line of sources is $\left(
a_{2}-a_{1}\right) $. Using this setting, for each $\alpha \in \left[
a_{1},a_{2}\right] $ we arrange the point source $\mathbf{x}_{\alpha
}:=\left( \alpha ,0,-d\right) $ located on the straight line $L_{\text{src}}$%
. We also define the near-field measurement site as the lower side of the
prism $\Omega ,$ 
\begin{equation*}
\Gamma :=\left\{ \mathbf{x}:\left\vert x\right\vert ,\left\vert y\right\vert
<R,z=-b\right\} .
\end{equation*}

Throughout the paper, we use either $\alpha $ or $\mathbf{x}_{\alpha }$ to
indicate the dependence of a function/parameter/number on those point
sources. We denote by $u$, $u_{i}$ and $u_{s}$ the total wave, incident wave
and scattered wave, respectively. Also, we note that $u=u_{i}+u_{s}$.

\subsection*{Forward problem}

Given the wavenumber $k>0$ and the function $c\left( \mathbf{x}\right) ,$ the
forward problem is to seek the function $\left. u\left( \mathbf{x},\alpha
\right) \right\vert _{\Gamma }$ such that $u=u\left( \mathbf{x},\alpha
\right) $ satisfying the Helmholtz equation with the radiation boundary
condition, which reads as 
\begin{align}
& \Delta u+k^{2}c\left( \mathbf{x}\right) u=-\delta \left( \mathbf{x}-%
\mathbf{x}_{\alpha }\right) \quad \text{in }\mathbb{R}^{3},  \label{eq:helm}
\\
& \lim_{r\rightarrow \infty }r\left( \partial _{r}u-\text{i}ku\right)
=0\quad \text{for }r=\left\vert x-x_{\alpha }\right\vert ,\text{i}=\sqrt{-1}.
\label{eq:somm}
\end{align}%
Here, the incident wave is 
\begin{equation}
u_{i}\left( \mathbf{x},\alpha \right) =\frac{\text{exp}\left( \text{i}%
k\left\vert \mathbf{x}-\mathbf{x}_{\alpha }\right\vert \right) }{4\pi
\left\vert \mathbf{x}-\mathbf{x}_{\alpha }\right\vert }.  \label{eq:3}
\end{equation}%
Moreover, we can deduce that the scattered wave satisfies the integral
equation: 
\begin{align}
u_{s}\left( \mathbf{x},\alpha \right) & =k^{2}\int_{\mathbb{R}^{3}}\frac{%
\text{exp}\left( \text{i}k\left\vert \mathbf{x}-\mathbf{x}^{\prime
}\right\vert \right) }{4\pi \left\vert \mathbf{x}-\mathbf{x}^{\prime
}\right\vert }\left( c\left( \mathbf{x}^{\prime }\right) -1\right) u\left( 
\mathbf{x}^{\prime },\alpha \right) d\mathbf{x}^{\prime }  \label{eq:4} \\
& =k^{2}\int_{\Omega }\frac{\text{exp}\left( \text{i}k\left\vert \mathbf{x}-%
\mathbf{x}^{\prime }\right\vert \right) }{4\pi \left\vert \mathbf{x}-\mathbf{%
x}^{\prime }\right\vert }\left( c\left( \mathbf{x}^{\prime }\right)
-1\right) u\left( \mathbf{x}^{\prime },\alpha \right) d\mathbf{x}^{\prime
},\quad \text{$\mathbf{x}$}\in \mathbb{R}^{3},  \notag
\end{align}%
since $c-1$ is compactly supported in $\Omega $; see again (\ref{eq:2}). Combining %
\eqref{eq:3} and \eqref{eq:4}, we find that the total wave $u$ satisfies the
Lippmann--Schwinger equation: 
\begin{equation*}
u\left( \mathbf{x},\alpha \right) =u_{i}\left( \mathbf{x},\alpha \right)
+k^{2}\int_{\Omega }\frac{\text{exp}\left( \text{i}k\left\vert \mathbf{x}-%
\mathbf{x}^{\prime }\right\vert \right) }{4\pi \left\vert \mathbf{x}-\mathbf{%
x}^{\prime }\right\vert }\left( c\left( \mathbf{x}^{\prime }\right)
-1\right) u\left( \mathbf{x}^{\prime },\alpha \right) d\mathbf{x}^{\prime
},\quad \text{$\mathbf{x}$}\in \mathbb{R}^{3}.
\end{equation*}

\begin{rem}\label{rem:2}
	Even though the propagation of the electromagnetic wave
	field is governed by the system of the Maxwell equations, we model our
	process by the single Helmholtz equation. The reason is twofold. First, it
	was demonstrated numerically in Appendix for \cite{KNN} that if the incident
	electric wave field $E=\left( E_{x},E_{y},E_{z}\right) $ has only one non
	zero component, then the propagation of this component in a heterogeneous
	medium is governed equally well by the Helmholtz equation and the Maxwell
	equations. This is true for at least a rather simply structured medium.
	Second, both our previous \cite%
	{Beilina2012,Klibanov2019b,Nguyen2018,Thanh2015} and current results for
	experimental data demonstrate a good reconstruction accuracy for the case of
	modeling by the Helmholtz equation. 
\end{rem} 

\begin{rem}
	Physically, $u\left( \mathbf{x},\alpha \right) $ is such a
	component of the electric field $E=\left( E_{x},E_{y},E_{z}\right) $ which
	was non zero when being incident, whereas other two components of the
	incident electric field $E$ equal zero. As it was mentioned in Remark \ref{rem:2},
	the propagation of this component in a heterogeneous medium is governed
	equally well by the Helmholtz equation and the full system of Maxwell's
	equations, as it was shown numerically in Appendix for \cite{KNN}. In the
	particular case of our experimental data, $u=E_{y}$ is indicated.
\end{rem} 

\subsection*{Coefficient inverse problem}

Given $k>0$, the \emph{inverse problem} is to reconstruct the smooth
dielectric constant $c\left( \mathbf{x}\right) $, $\mathbf{x}\in \Omega $
satisfying conditions (\ref{eq:2}) from the boundary measurement as near-field data 
\begin{equation}
F\left( \mathbf{x},\mathbf{x}_{\alpha }\right) =u\left( \mathbf{x},\alpha
\right) \quad \text{for }\mathbf{x}\in \Gamma ,\mathbf{x}_{\alpha }\in L_{%
\text{src}},  \label{eq:bdr}
\end{equation}%
where $u\left( \mathbf{x},\alpha \right) $ is the total wave associated to
the incident wave $u_{i}$ of \eqref{eq:3}. A schematic diagram of locations
of sources and detectors is presented on Figure \ref{fig:Mesh-refinement} in subsection \ref{subsec:Experimental-setup}.

Currently uniqueness theorem for this CIP can be proven only in the case
when the right hand side of equation (\ref{eq:helm}) is not vanishing in $%
\overline{\Omega }$. Such a theorem can be proven by the method of \cite%
{BukhKlib} which was mentioned in the Introduction. In addition, uniqueness can
be proven if working within an approximate mathematical model. In this
regard, uniqueness follows from Theorem 3.2 of \cite{Khoa2019}, and, for the
approximate mathematical model of this paper, it follows from Theorem 4.

\section{A globally convergent numerical method\label%
{sec:A-globally-convergent}}

\subsection{Derivation of a system of coupled quasilinear elliptic PDEs}\label{subsec:3.1}

Since our line of sources $L_{\text{src}}$ is located outside of $\overline{%
\Omega }$, we deduce this system from the Helmholtz homogeneous version of
equation \eqref{eq:helm} and for each $\alpha \in \left[ a_{1},a_{2}\right] 
$ 
\begin{equation}
\Delta u+k^{2}c\left( \mathbf{x}\right) u=0\quad \text{in }\Omega .
\label{eq:helm1}
\end{equation}

Now, we set that 
\begin{equation*}
\log u_{i}\left(\mathbf{x},\alpha\right)=\text{i}k\left|\mathbf{x}-\mathbf{x}%
_{\alpha}\right|-\log\left(4\pi\left|\mathbf{x}-\mathbf{x}%
_{\alpha}\right|\right),
\end{equation*}
which then leads to 
\begin{equation}
\nabla\left(\log u_{i}\left(\mathbf{x},\alpha\right)\right)=\frac{\text{i}%
k\left(\mathbf{x}-\mathbf{x}_{\alpha}\right)}{\left|\mathbf{x}-\mathbf{x}%
_{\alpha}\right|}-\frac{\mathbf{x}-\mathbf{x}_{\alpha}}{\left|\mathbf{x}-%
\mathbf{x}_{\alpha}\right|^{2}}.  \label{eq:gradui}
\end{equation}

Even though we work with a fixed value of $k$, we now temporarily indicate the
dependence of the function $u$ on $k$, i.e. $u=u\left( \mathbf{x},\alpha
,k\right) .$ Using the asymptotic behavior of $u\left( \mathbf{x},\alpha
,k\right) $ at $k\rightarrow \infty $ (cf. \cite{KR}), one can prove, similarly
with \cite{KR}, that there exists a sufficiently large number $\overline{k}>0
$ such that $u\left( \mathbf{x},\alpha ,k\right) \neq 0$ for $k\geq 
\overline{k},$\textbf{\ }$\mathbf{x}\in \overline{\Omega },$ $\alpha \in %
\left[ a_{1},a_{2}\right] .$ Furthermore, it was shown in \cite%
{Klibanov2019,Klibanov2019b} that, using that asymptotic behavior, one can
uniquely define the function $\log u\left( \mathbf{x},\alpha ,k\right) $ for 
$k\geq \overline{k},$\textbf{\ }$\mathbf{x}\in \overline{\Omega },$ $\alpha
\in \left[ a_{1},a_{2}\right] .$ Thus, we assume below that our specific
value of the wave number $k=k_{0}$ we work with is sufficiently large, $%
k_{0}\geq \overline{k},$ and also that the function $\log u\left( \mathbf{x}%
,\alpha ,k_{0}\right) $ is uniquely defined as in \cite%
{Klibanov2019,Klibanov2019b}. In particular, we work with the dimensionless
values $k_{0}\geq 6.62$; see subsection \ref{sec:4.4}. We point out that since we work
below only with the derivatives (with respect to $\mathbf{x}$ and $\alpha$) of $\log u\left( \mathbf{x},\alpha
,k_{0}\right) $, i.e.
\begin{align*}
&  \nabla \left( \log u\left( \mathbf{x},\alpha ,k_{0}\right)
\right) := \nabla _{\mathbf{x}}\left( \log u\left( \mathbf{x},\alpha ,k_{0}\right)
\right)  = \nabla _{\mathbf{x}}u\left( \mathbf{x},\alpha ,k_{0}\right)
/u\left( \mathbf{x},\alpha ,k_{0}\right) , \\
& \Delta \left( \log
u\left( \mathbf{x},\alpha ,k_{0}\right) \right) := \Delta _{\mathbf{x}}\left( \log
u\left( \mathbf{x},\alpha ,k_{0}\right) \right) , \\
&  \partial _{\alpha }\left( \nabla \left( \log u\left( \mathbf{x%
},\alpha ,k_{0}\right) \right) \right) :=  \partial _{\alpha }\left( \nabla _{\mathbf{x}}\left( \log u\left( \mathbf{x%
},\alpha ,k_{0}\right) \right) \right) , \\
& 
\partial _{\alpha }\left( \Delta\left( \log u\left( \mathbf{x},\alpha ,k_{0}\right) \right)
\right)  := \partial _{\alpha }\left( \Delta _{%
\mathbf{x}}\left( \log u\left( \mathbf{x},\alpha ,k_{0}\right) \right)
\right) ,
\end{align*}%
then what we actually need in our work is not that definition of $\log
u\left( \mathbf{x},\alpha ,k_{0}\right) $ but rather that $u\left( \mathbf{x}%
,\alpha ,k_{0}\right) \neq 0.$ However, we have not noticed in our
computations those values of $\left\vert u\left( \mathbf{x},\alpha
,k_{0}\right) \right\vert $ which would be close to zero.

Denote $v_{0}\left( \mathbf{x},\alpha \right) =u\left( \mathbf{x},\alpha
\right) /u_{i}\left( \mathbf{x},\alpha \right) $. We define the function $%
v\left( \mathbf{x},\alpha \right) $ as 
\begin{equation*}
v\left( \mathbf{x},\alpha \right) :=\log \left( v_{0}\left( \mathbf{x}%
,\alpha \right) \right) =\log \left( u\left( \mathbf{x},\alpha \right)
\right) -\log \left( u_{i}\left( \mathbf{x},\alpha \right) \right) \quad 
\text{for }\mathbf{x}\in \Omega ,\alpha \in \left[ a_{1},a_{2}\right] .
\end{equation*}%
Hence, one computes that
\begin{equation}
\nabla v\left( \mathbf{x},\alpha \right) =\frac{\nabla v_{0}\left( \mathbf{x}%
,\alpha \right) }{v_{0}\left( \mathbf{x},\alpha \right) },\quad \Delta
v\left( \mathbf{x},\alpha \right) =\frac{\Delta v_{0}\left( \mathbf{x}%
,\alpha \right) }{v_{0}\left( \mathbf{x},\alpha \right) }-\left( \frac{%
\nabla v_{0}\left( \mathbf{x},\alpha \right) }{v_{0}\left( \mathbf{x},\alpha
\right) }\right) ^{2}.  \label{eq:gradel}
\end{equation}%
Thus, using \eqref{eq:helm1}--\eqref{eq:gradel} we derive the equation for $v$%
, 
\begin{equation}
\Delta v+\left( \nabla v\right) ^{2}+2\nabla v\cdot \nabla \left( \log
\left( u_{i}\left( \mathbf{x},\alpha \right) \right) \right) =-k^{2}\left(
c\left( \mathbf{x}\right) -1\right) \quad \text{for }\mathbf{x}\in \Omega .
\label{eq:v}
\end{equation}

\begin{rem}\label{rem:X}
	Observe that if the function $v\left( \mathbf{x}%
	,\alpha \right) $ is known, then we can immediately find the target
	coefficient $c\left( \mathbf{x}\right) $ via equation (\ref{eq:v}).
\end{rem}
 
We now differentiate \eqref{eq:v} with respect to $\alpha $ and use %
\eqref{eq:gradui} to obtain the following third-order PDE: 
\begin{align*}
 \Delta \partial _{\alpha }v+2\nabla v\cdot \nabla \partial _{\alpha }v
 +2\nabla \partial _{\alpha }v\cdot \left[ \frac{\text{i}k\left( \mathbf{x}-%
\mathbf{x}_{\alpha }\right) }{\left\vert \mathbf{x}-\mathbf{x}_{\alpha
}\right\vert }-\frac{\mathbf{x}-\mathbf{x}_{\alpha }}{\left\vert \mathbf{x}-%
\mathbf{x}_{\alpha }\right\vert ^{2}}\right] +2\partial _{\alpha }\left[ 
\frac{\text{i}k\left( \mathbf{x}-\mathbf{x}_{\alpha }\right) }{\left\vert 
\mathbf{x}-\mathbf{x}_{\alpha }\right\vert }-\frac{\mathbf{x}-\mathbf{x}%
_{\alpha }}{\left\vert \mathbf{x}-\mathbf{x}_{\alpha }\right\vert ^{2}}%
\right] \cdot \nabla v=0.
\end{align*}%
To simplify the presentation, this equation can be rewritten as: 
\begin{equation}
\Delta \partial _{\alpha }v+2\nabla v\cdot \nabla \partial _{\alpha
}v+2\nabla \partial _{\alpha }v\cdot \tilde{\mathbf{x}}_{\alpha }+2\hat{%
\mathbf{x}}_{\alpha }\cdot \nabla v=0\quad \text{for }\mathbf{x\in \Omega },
\label{eq:vv}
\end{equation}%
where, for $\mathbf{x}-\mathbf{x}_{\alpha }=\left( x-\alpha ,y,z+d\right)
\in \mathbb{R}^{3},$ 
\begin{align*}
\tilde{\mathbf{x}}_{\alpha }& =\frac{\text{i}k\left( \mathbf{x}-\mathbf{x}%
_{\alpha }\right) }{\left\vert \mathbf{x}-\mathbf{x}_{\alpha }\right\vert }-%
\frac{\mathbf{x}-\mathbf{x}_{\alpha }}{\left\vert \mathbf{x}-\mathbf{x}%
_{\alpha }\right\vert ^{2}}, \\
\hat{\mathbf{x}}_{\alpha }& =\frac{\text{i}k}{\left\vert \mathbf{x}-\mathbf{x%
}_{\alpha }\right\vert ^{3}}\left( -y^{2}-\left( z+d\right) ^{2},\left(
x-\alpha \right) y,\left( x-\alpha \right) z\right) \\
& -\frac{1}{\left\vert \mathbf{x}-\mathbf{x}_{\alpha }\right\vert ^{4}}%
\left( \left( x-\alpha \right) ^{2}-y^{2}-\left( z+d\right) ^{2},2\left(
x-\alpha \right) y,2\left( x-\alpha \right) z\right) .
\end{align*}

Even though the unknown coefficient $c\left( \mathbf{x}\right) $ is not
present in equation (\ref{eq:vv}) for the function $v\left( \mathbf{x}%
,\alpha \right) $, it is still not easy to solve the nonlinear third-order
PDE \eqref{eq:vv}. To circumvent this, we rely on a special orthonormal
basis with respect to $\alpha $ to turn \eqref{eq:vv} into a system of
coupled quasilinear elliptic PDEs.

For $\alpha \in \left( a_{1},a_{2}\right) $, let $\left\{ \Psi _{n}\left(
\alpha \right) \right\} _{n=0}^{\infty }$ be the special orthonormal basis
in $L^{2}\left( a_{1},a_{2}\right) $, which was first proposed in \cite%
{Klibanov2017}. Herewith, the construction of this basis is shortly
revisited. For each $n\in \mathbb{N}$, let $\varphi _{n}\left( \alpha
\right) =\alpha ^{n}e^{\alpha }$ for $\alpha \in \left[ a_{1},a_{2}\right] $%
. The set $\left\{ \varphi _{n}\left( \alpha \right) \right\} _{n=0}^{\infty
}$ is linearly independent and complete in $L^{2}\left( a_{1},a_{2}\right) $%
. Using the Gram--Schmidt orthonormalization procedure, we can obtain the
orthonormal basis $\left\{ \Psi _{n}\left( \alpha \right) \right\}
_{n=0}^{\infty }$ in $L^{2}\left( a_{1},a_{2}\right) $, which possesses the
following special properties:

\begin{itemize}
\item $\Psi_{n}\in C^{\infty}\left[a_{1},a_{2}\right]$ for all $n\in\mathbb{N%
}$;

\item Let $s_{mn}=\left\langle \Psi_{n}^{\prime },\Psi_{m}\right\rangle $,
where $\left\langle \cdot,\cdot\right\rangle $ denotes the scalar product in 
$L^{2}\left(a_{1},a_{2}\right)$. Then the square matrix $S_{N}=\left(s_{mn}%
\right)_{m,n=0}^{N-1}$ for $N\in\mathbb{N}^{*}$ is invertible in the sense
that 
\begin{equation*}
s_{mn}=%
\begin{cases}
1 & \text{if }n=m, \\ 
0 & \text{if }n<m.%
\end{cases}%
\end{equation*}
\end{itemize}

We note that even though the Gram--Schmidt procedure is unstable if using
the infinite number of functions $\varphi _{n}\left( \alpha \right) ,$ we
have observed numerically that it has good stability properties for rather
small numbers $N$ which we use. Neither classical orthogonal polynomials nor
the classical basis of trigonometric functions do not hold the second
property. The matrix $S_{N}$ is an upper diagonal matrix with $\det \left(
S_{N}\right) =1$. On the other hand, the special second property allows us
to reduce the third-order PDE \eqref{eq:vv} to a system of coupled elliptic
PDEs.

Given $N\in \mathbb{N}^{\ast }$, we consider the following truncated Fourier
series for $v$: 
\begin{equation}
v\left( \mathbf{x},\alpha \right) =\sum_{n=0}^{N-1}\left\langle v\left( 
\mathbf{x},\cdot \right) ,\Psi _{n}\left( \cdot \right) \right\rangle \Psi
_{n}\left( \alpha \right) \quad \text{for }\mathbf{x}\in \Omega ,\alpha \in %
\left[ a_{1},a_{2}\right] .  \label{eq:fourier}
\end{equation}%
Substitution \eqref{eq:fourier} into \eqref{eq:vv} gives 
\begin{align}
& \sum_{n=0}^{N-1}\Psi _{n}^{\prime }\left( \alpha \right) \Delta
v_{n}\left( \mathbf{x}\right) +2\sum_{n=0}^{N-1}\sum_{l=0}^{N-1}\Psi
_{n}\left( \alpha \right) \Psi _{l}^{\prime }\left( \alpha \right) \nabla
v_{n}\left( \mathbf{x}\right) \cdot \nabla v_{l}\left( \mathbf{x}\right)
\label{eq:afterfourier} \\
& +2 \sum_{n=0}^{N-1}\Psi _{n}^{\prime }\left( \alpha \right)\nabla
v_{n}\left( \mathbf{x}\right) \cdot \tilde{\mathbf{x}}_{\alpha }+2 \sum_{n=0}^{N-1}\Psi
_{n}\left( \alpha \right)\hat{\mathbf{x}}_{\alpha }\cdot
\nabla v_{n}\left( \mathbf{x}\right) =0.  \notag
\end{align}%
Multiplying both sides of \eqref{eq:afterfourier} by $\Psi _{m}\left( \alpha
\right) $ for $0\leq m\leq N-1$ and then integrating the resulting equation
with respect to $\alpha $, we obtain the following system of coupled
elliptic equations: 
\begin{align}
& \Delta V\left( \mathbf{x}\right) +K\left( \nabla V\left( \mathbf{x}\right)
\right) =0\quad \text{for }\mathbf{x}\in \Omega ,  \label{eq:pdeV} \\
& \nabla V\left( \mathbf{x}\right) \cdot \text{n}=0\quad \text{for }\mathbf{x%
}\in \partial \Omega \backslash \Gamma ,  \label{eq:bdrV1} \\
& V\left( \mathbf{x}\right) =\psi _{0}\left( \mathbf{x}\right) ,V_{z}\left( 
\mathbf{x}\right) =\psi _{1}\left( \mathbf{x}\right) \quad \text{for }%
\mathbf{x}\in \Gamma .  \label{eq:bdrV2}
\end{align}%
Here, $V\left( \mathbf{x}\right) \in \mathbb{R}^{N}$ is the unknown vector
function given by 
\begin{equation}\label{3}
V\left( \mathbf{x}\right) =%
\begin{pmatrix}
v_{0}\left( \mathbf{x}\right) & v_{1}\left( \mathbf{x}\right) & \cdots & 
v_{N-1}\left( \mathbf{x}\right)%
\end{pmatrix}%
^{T}.
\end{equation}%
Thus, we have obtained a boundary value problem (\ref{eq:pdeV})--(\ref%
{eq:bdrV2}) for a system of coupled quasilinear elliptic PDEs (\ref{eq:pdeV}%
). Boundary conditions (\ref{eq:bdrV2}) are overdetermined ones. Note that
conditions (\ref{eq:bdrV2}) are the Cauchy data. A Lipschitz stability
estimate for problem (\ref{eq:pdeV})--(\ref{eq:bdrV2}) is obtained in \cite%
{Khoa2019}.

In PDE \eqref{eq:pdeV}, we denote by $K\left( \nabla V\left( \mathbf{x}%
\right) \right) =S_{N}^{-1}f\left( \nabla V\left( \mathbf{x}\right) \right) $%
, where $f=\left( \left( f_{m}\right) _{m=0}^{N-1}\right) ^{T}\in \mathbb{R}%
^{N}$ is quadratic with respect to the first derivative of components of $%
V\left( \mathbf{x}\right) $, 
\begin{align}\label{eq:17}
f_{m}\left( \nabla V\left( \mathbf{x}\right) \right) &
=2\sum_{n,l=0}^{N-1}\nabla v_{n}\left( \mathbf{x}\right) \cdot \nabla
v_{l}\left( \mathbf{x}\right) \int_{a_{1}}^{a_{2}}\Psi _{m}\left( \alpha
\right) \Psi _{n}\left( \alpha \right) \Psi _{l}^{\prime }\left( \alpha
\right) d\alpha   \\
& +2\sum_{n=0}^{N-1}\int_{a_{1}}^{a_{2}}\Psi _{m}\left( \alpha \right) \Psi
_{n}^{\prime }\left( \alpha \right) \nabla v_{n}\left( \mathbf{x}\right)
\cdot \tilde{\mathbf{x}}_{\alpha }d\alpha  \notag +2\sum_{n=0}^{N-1}\int_{a_{1}}^{a_{2}}\Psi _{m}\left( \alpha \right) \Psi
_{n}\left( \alpha \right) \hat{\mathbf{x}}_{\alpha }\cdot \nabla v_{n}\left( 
\mathbf{x}\right) d\alpha .  
\end{align}%
Now, we explain how to obtain the Cauchy data \eqref{eq:bdrV2}. First, an
obvious application of the expansion \eqref{eq:fourier} to the near-field
data \eqref{eq:bdr} gives us the function $\psi _{0}$ at $\Gamma $. However,
our experimental data are in fact far-field type; see section \ref%
{sec:Experimental-results}. This means that we collect the experimental data
far from the domain of interest, i.e. on the plane $\left\{ z=-D\right\} $
for $D>b$. However, the far-field data do not look nice as we experienced in
many previous works; cf. e.g. \cite{Nguyen2018}. Therefore, we rely on the
so-called \textquotedblleft data propagation technique'' to get the data
closer to the target's side and to reduce the side of the domain of interest 
$\Omega $ we are considering. In our work, this technique is revisited in
subsection \ref{subsec:Data-propagation-revisited}. As a by-product of this
technique, we can obtain an approximation of the $z$--derivative of the
function $u$ at the near-field measurement site $\Gamma $. This leads to the
presence of the function $\psi _{1}$ in \eqref{eq:bdrV2}.

We now explain the Neumann zero boundary condition \eqref{eq:bdrV1} on $%
\partial \Omega \backslash \Gamma $. Assuming that, due to the radiation
condition (\ref{eq:somm}), $\left( \partial _{n}u-\text{i}ku\right) \left( \mathbf{x}%
,\alpha \right) \mid _{\partial \Omega \backslash \Gamma }\approx 0$ and
also that $\left( \partial _{n}u_{i}-\text{i}ku_{i}\right) \mid _{\partial \Omega
\backslash \Gamma }\approx 0,$ we easily obtain $\partial _{n}v\left( 
\mathbf{x},\alpha \right) \mid _{\partial \Omega \backslash \Gamma }\approx
0 $ for $v\left( \mathbf{x},\alpha \right) =\log \left( u\left( \mathbf{x}%
,\alpha \right) /u_{i}\left( \mathbf{x},\alpha \right) \right) .$ Hence,
condition \eqref{eq:bdrV1} follows from the latter and is an approximate one.

\begin{rem}\label{rem:1}
	Using the truncated Fourier series (\ref{eq:fourier}), we
	actually replace the original CIP with its approximation. Thus, we work with
	an approximate mathematical model. Substantial difficulties in solving our
	CIP are linked with both its nonlinearity and ill-posedness. Due to these
	difficulties, we cannot prove convergence at $N\rightarrow \infty $, where $N
	$ is the number of Fourier coefficients in that truncated Fourier series.
	Besides, we recall that another feature of our approximate mathematical
	model is that we do not allow the grid step size of our partial finite
	differences to tend to zero. We believe that good reconstruction results,
	which we demonstrate below, justify our approximate mathematical model. It
	is worthy to mention here hat it is quite \textbf{typical} in the field of
	CIPs to consider approximate mathematical models caused by truncations of
	certain Fourier series without proofs of convergence at $N\rightarrow \infty
	.$ In such cases good numerical results serve as justifications of these
	models. In this regard we refer to works \cite{Ag,Al,Kab1,Kab2,N2,N3}. More
	detailed discussions of the issue of approximate mathematical models can be
	found in \cite{Khoa2019,Klibhyp}.
\end{rem}

\subsection{A weighted cost functional}

We set 
\begin{equation}
L\left( V\right) \left( \mathbf{x}\right) =\Delta V\left( \mathbf{x}\right)
+K\left( \nabla V\left( \mathbf{x}\right) \right) .  \label{eq:LL}
\end{equation}%
We use partial finite difference setting via considering finite differences
with respect to $x,y$ and keeping the standard derivatives with respect to $%
z $. In this setting, we look for an approximate solution of the system %
\eqref{eq:pdeV}--\eqref{eq:bdrV2} using the minimization of a Tikhonov-type
functional weighted by a suitable Carleman Weight Function (CWF).

Let the numbers $\theta >b$ and $\lambda \geq 1$. We define our CWF as 
\begin{equation}
\mu _{\lambda }\left( z\right) =\exp \left( 2\lambda \left( z-\theta \right)
^{2}\right) \quad \text{for }z\in \left[ -b,b\right] .  \label{eq:CWF}
\end{equation}%
The choice $\theta >b$ is based on the fact that the gradient of the CWF
should not vanish in the closed domain $\overline{\Omega }$. Observe that
the CWF is decreasing for $z\in \left( -b,b\right) $ and 
\begin{align}
	& \max_{z\in \left[ -b,b\right] }\mu _{\lambda }\left( z\right) =\mu
	_{\lambda }\left( -b\right) =e^{2\lambda \left( b+\theta \right) ^{2}},
	\label{max} \\
	& \min_{z\in \left[ -b,b\right] }\mu _{\lambda }\left( z\right) =\mu
	_{\lambda }\left( b\right) =e^{2\lambda \left( b-\theta \right) ^{2}}.
	\label{max1}
\end{align}%
This means that the CWF \eqref{eq:CWF} attains its maximal value on the
measurement site $\Gamma $, and it attains its minimal value on the opposite
side. This notion of using the CWF is essential because it \textquotedblleft
maximizes\textquotedblright\ the influence of the actual measured boundary
data at $z=-b$. Furthermore, the CWF plays the vital role in convexifying
the cost functional globally in both this and other types of the
convexification method. In fact what the CWF does is that it controls the
nonlinear term. In our particular case this term is $K\left( \nabla V\left( 
\mathbf{x}\right) \right) $ in \eqref{eq:pdeV}.

\subsubsection{Preliminaries of the partial finite differences setting \label%
	{subsec:Preliminaries-on-settings}}

We use the same grid step size $h$ in $x$ and $y$ directions. Introduce two
partitions of the closed interval $\left[ -R,R\right] ,$ 
\begin{align*}
	& -R=x_{0}<x_{1}<\ldots <x_{Z_{h}-1}<x_{Z_{h}}=R,\quad x_{p}-x_{p-1}=h, \\
	& -R=y_{0}<y_{1}<\ldots <y_{Z_{h}-1}<y_{Z_{h}}=R,\quad y_{q}-y_{q-1}=h.
\end{align*}%
We write the differential operator in \eqref{eq:LL} in the following partial
finite difference form. For any $N$-dimensional vector function $u\left( 
\mathbf{x}\right) $, we denote by $u_{p,q}^{h}\left( z\right) =u\left(
x_{p},y_{q},z\right) $ the corresponding semi-discrete function defined at
grid points $\left\{ \left( x_{p},y_{q}\right) \right\} _{p,q=0}^{Z_{h}}$.
Thus, the interior grid points are $\left\{ \left( x_{p},y_{q}\right)
\right\} _{p,q=1}^{Z_{h}-1}$. Denote 
\begin{equation*}
	\Omega _{h}=\left\{ \left( x_{p},y_{q},z\right) :\left\{ \left(
	x_{p},y_{q}\right) \right\} _{p,q=0}^{Z_{h}-1}\subset \left[ -R,R\right]
	\times \left[ -R,R\right] ,z\in \left( -b,b\right) \right\} ,
\end{equation*}%
\begin{equation*}
	\Gamma _{h}=\left\{ \left( x_{p},y_{q},-b\right) :\left\{ \left(
	x_{p},y_{q}\right) \right\} _{p,q=0}^{Z_{h}-1}\subset \left[ -R,R\right]
	\times \left[ -R,R\right] \right\} .
\end{equation*}%
Henceforth, the corresponding Laplace operator in the partial finite
differences is given by 
\begin{equation*}
	\Delta ^{h}u^{h}=u_{zz}^{h}+u_{xx}^{h}+u_{yy}^{h},
\end{equation*}%
where, for interior points of $\Omega _{h}$, we use 
\begin{equation*}
	u_{xx}^{h}=h^{-2}\left( u_{p+1,q}^{h}\left( z\right) -2u_{p,q}^{h}\left(
	z\right) +u_{p-1,q}^{h}\left( z\right) \right) ,\text{ }p,q\in \left[
	1,Z_{h}-1\right]
\end{equation*}%
and similarly for the finite difference operator $u_{yy}^{h}$. As to the
gradient operator, we write for interior points $\nabla ^{h}u_{p,q}\left(
z\right) =\left( \partial _{x}^{h}u_{p,q}\left( z\right) ,\partial
_{y}^{h}u_{p,q}\left( z\right) ,\partial _{z}u_{p,q}^{h}\left( z\right)
\right) $, where 
\begin{equation*}
	\partial _{x}^{h}u_{p,q}^{h}\left( z\right) =\left( 2h\right) ^{-1}\left(
	u_{p+1,q}^{h}\left( z\right) -u_{p-1,q}^{h}\left( z\right) \right) .
\end{equation*}

To simplify the presentation, we consider any $N$--D complex valued function 
$W=\func{Re}W+$i$\func{Im}W$ as the $2N$--D vector function with real valued
components $\left( \func{Re}W,\func{Im}W\right) :=\left( W_{1},W_{2}\right)
:=W\in \mathbb{R}^{2N}$. Also, below for any complex number $a\in \mathbb{C}$
we denote $\overline{a}$ its complex conjugate.

Denote $w^{h}\left( z\right) $ the vector function $w^{h}\left( z\right)
=\left\{ w_{p,q}\left( z\right) \right\} _{p,q=0}^{Z_{h}},$ where $%
w_{p,q}^{h}\left( z\right) =w\left( x_{p},y_{q},z\right) .$ We introduce the
Hilbert spaces $H_{2N}^{2,h}=H_{2N}^{2,h}\left( \Omega _{h}\right) $ and $%
L_{2N}^{2,h}=L_{2N}^{2,h}\left( \Omega _{h}\right) $ of semi-discrete
complex valued functions as follows: 
\begin{align*}
	& H_{2N}^{2,h}=\left\{ w^{h}\left( z\right) :\left\Vert w^{h}\right\Vert
	_{H_{2N}^{2,h}}^{2}:=\sum_{p,q=1}^{Z_{h}-1}\sum_{m=0}^{2}h^{2}\int_{-b}^{b}%
	\left\vert \partial _{z}^{m}w_{p,q}^{h}\left( z\right) \right\vert
	^{2}dz<\infty \right\} , \\
	& L_{2N}^{2,h}=\left\{ w^{h}\left( z\right) :\left\Vert w^{h}\right\Vert
	_{L_{2N}^{2,h}}^{2}:=\sum_{p,q=1}^{Z_{h}-1}h^{2}\int_{-b}^{b}\left\vert
	w_{p,q}^{h}\left( z\right) \right\vert ^{2}dz<\infty \right\} .
\end{align*}%
We also define the subspace $H_{2N,0}^{2,h}\subset H_{2N}^{2,h}$ as 
\begin{equation*}
	H_{2N,0}^{2,h}=\left\{ w^{h}\left( z\right) \in H_{2N}^{2,h}:\left. \nabla
	^{h}w_{p,q}^{h}\left( z\right) \right\vert _{\partial \Omega _{h}\backslash
		\Gamma _{h}}\cdot \text{n}=0,\left. w_{p,q}^{h}\left( z\right) \right\vert
	_{\Gamma _{h}}=\left. \partial _{z}w_{p,q}^{h}\left( z\right) \right\vert
	_{\Gamma _{h}}=0\right\} .
\end{equation*}%
Let $h_{0}>0$ be a fixed positive constant. We assume below that 
\begin{equation}
h\geq h_{0}>0.  \label{190}
\end{equation}%
Thus, (\ref{190}) means that we do not allow the grid step size to tend to
zero.

\subsubsection{The semi-discrete form of the weighted cost functional and
	Carleman estimate}

Following (\ref{3}), denote%
\begin{equation*}
	V^{h}\left( z\right) =%
	\begin{pmatrix}
		v_{0}^{h}\left( z\right) & v_{1}^{h}\left( z\right) & \cdots & 
		v_{N-1}^{h}\left( z\right)%
	\end{pmatrix}%
	^{T}.
\end{equation*}%
We now consider the following weighted Tikhonov-like functional $J_{\lambda
}:H_{2N}^{2,h}\left( \Omega _{h}\right) \rightarrow \mathbb{R}_{+}$, 
\begin{equation}
J_{h,\lambda }\left( V^{h}\right)
=\sum_{p,q=1}^{Z_{h}-1}h^{2}\int_{-b}^{b}\left\vert L^{h}\left( V^{h}\left(
z\right) \right) \right\vert ^{2}\mu _{\lambda }\left( z\right) dz,
\label{eq:J}
\end{equation}%
where the CWF $\mu _{\lambda }\left( z\right) $ is defined in \eqref{eq:CWF}
and $L^{h}\left( V^{h}\left( z\right) \right) $ is the operator \eqref{eq:LL}
written in partial finite differences, 
\begin{equation}
L^{h}\left( V^{h}\left( z\right) \right) =\Delta ^{h}V^{h}\left( z\right)
+K\left( \nabla ^{h}V^{h}\left( z\right) \right) .  \label{98}
\end{equation}

Let $M>0$ be an arbitrary number. We define the set $B\left( M\right)
\subset H_{2N}^{2,h}\left( \Omega _{h}\right) $ as 
\begin{equation}
B\left( M\right) :=\left\{ V^{h}\in H_{2N}^{2,h}:\left\Vert V^{h}\right\Vert
_{H_{2N}^{2,h}}<M,\text{ }V_{\Gamma _{h}}^{h}=\psi _{0}^{h},\text{ }\partial
_{z}V_{\Gamma _{h}}^{h}=\psi _{1}^{h}\right\} .  \label{eq:BM}
\end{equation}%
The embedding theorem and (\ref{190}) imply that 
\begin{equation}
\overline{B\left( M\right) }\subset C_{2N}^{1,h}\left( \overline{\Omega _{h}}%
\right) \;\text{and }\left\Vert V^{h}\right\Vert _{C_{2N}^{1,h}\left( 
	\overline{\Omega _{h}}\right) }\leq C_{1}\quad \text{for all }V^{h}\in 
\overline{B\left( M\right) }.  \label{eq:BMM}
\end{equation}%
Here the number $C_{1}>0$ depends only on $M$ and $h_{0}.$ Since the lower
estimate $h_{0}$ of our grid step size $h$ is fixed by (\ref{190}), we
neglect below the $h-$dependence of constants used in proofs of our results.

\subsection*{Minimization problem}

This problem is formulated as: \emph{Minimize the cost functional }$%
J_{\lambda }\left( V^{h}\right) $\emph{\ on the set }$\overline{B\left(
	M\right) }$\emph{.}

Below we prove a one-dimensional Carleman estimate. This estimate resembles
Lemma 3.1 of \cite{Klibanov2017a}. However, the different Carleman Weight
Function we use here requires a different proof of the target estimate. In
formulations and proofs of Carleman estimates below only for real valued
functions $u$ are used, since $\left\vert u\right\vert ^{2}=\left( \func{Re}%
u\right) ^{2}+\left( \func{Im}u\right) ^{2}$ for complex valued ones.

\begin{lem}
	\label{lem:carle}For all real valued functions $u\in H^{2}\left( -b,b\right) 
	$ such that $u\left( -b\right) =u^{\prime }\left( -b\right) =0$ and for all $%
	\lambda \geq 1$ the following Carleman estimate holds: 
	\begin{align}
		\int_{-b}^{b}\left( u^{\prime \prime }\right) ^{2}\mu _{\lambda }\left(
		z\right) dz& \geq C\int_{-b}^{b}\left( u^{\prime \prime }\right) ^{2}\mu
		_{\lambda }\left( z\right) dz  \notag \\
		& +C\lambda \int_{-b}^{b}\left( u^{\prime }\right) ^{2}\mu _{\lambda }\left(
		z\right) dz+C\lambda ^{3}\int_{-b}^{b}u^{2}\mu _{\lambda }\left( z\right) dz.
		\label{eq:Car1}
	\end{align}%
	Here and below the constant $C>0$ depends only on numbers $r$ and $b$.
\end{lem}

\emph{Proof.} We prove this estimate for functions $u\in C^{2}\left[ -b,b%
\right] $ satisfying $u\left( -b\right) =u^{\prime }\left( -b\right) =0.$
Extension for the case $u\in H^{2}\left( -b,b\right) $ can be done using
density arguments. Recall that $r>b.$ Introduce the function $v=u\exp \left(
\lambda \left( z-\theta \right) ^{2}\right) $. Then $u=v\exp \left( -\lambda
\left( z-\theta \right) ^{2}\right) $. Hence, 
\begin{equation*}
	u_{zz}=\left( v_{zz}-4\lambda \left( z-\theta \right) v_{z}+4\lambda
	^{2}\left( z-\theta \right) ^{2}\left( 1+\mathcal{O}\left( 1/\lambda \right)
	\right) v\right) e^{-\lambda \left( z-\theta \right) ^{2}}.
\end{equation*}%
Here and below $\mathcal{O}\left( 1/\lambda \right) $ denotes different $%
C^{2}-$functions satisfying the estimate $\left\vert \mathcal{O}\left(
1/\lambda \right) \right\vert \leq C/\lambda ,\forall \lambda \geq 1$
together with their derivatives up to the second order. Therefore, 
\begin{align}
	\left( u_{zz}\right) ^{2}e^{2\lambda \left( z-\theta \right) ^{2}}& =\left(
	v_{zz}-4\lambda \left( z-\theta \right) v_{z}+4\lambda ^{2}\left( z-\theta
	\right) ^{2}\left( 1+\mathcal{O}\left( 1/\lambda \right) \right) v\right)
	^{2}  \notag \\
	& \geq -8\lambda \left( z-\theta \right) v_{z}\left( v_{zz}+4\lambda
	^{2}\left( z-\theta \right) ^{2}\left( 1+\mathcal{O}\left( 1/\lambda \right)
	\right) v\right)  \notag \\
	& =\left( -4\lambda \left( z-\theta \right) v_{z}^{2}\right) _{z}+4\lambda
	v_{z}^{2}+\left( -16\lambda ^{3}\left( z-\theta \right) ^{3}\left( 1+%
	\mathcal{O}\left( 1/\lambda \right) \right) v^{2}\right) _{z}  \notag \\
	& +48\lambda ^{3}\left( z-\theta \right) ^{2}\left( 1+\mathcal{O}\left(
	1/\lambda \right) \right) v^{2}.  \label{eq:22}
\end{align}%
Since $u\left( -b\right) =u^{\prime }\left( -b\right) =0$, then integrating
the estimate \eqref{eq:22} over $z\in \left( -b,b\right) $, we obtain 
\begin{equation}
\int_{-b}^{b}\left( u_{zz}\right) ^{2}e^{2\lambda \left( z-\theta \right)
	^{2}}dz\geq 4\lambda \int_{-b}^{b}v_{z}^{2}dz+47\lambda
^{3}\int_{-b}^{b}\left( z-\theta \right) ^{2}u^{2}e^{2\lambda \left(
	z-\theta \right) ^{2}}dz.  \label{eq:23-1}
\end{equation}%
We have used here the fact that 
\begin{align*}
	& \int_{-b}^{b}\left( -4\lambda \left( z-\theta \right) v_{z}^{2}\right)
	_{z}dz=\left. -4\lambda \left( z-\theta \right) \left( u_{z}e^{\lambda
		\left( z-\theta \right) ^{2}}+2\lambda \left( z-\theta \right) ue^{\lambda
		\left( z-\theta \right) ^{2}}\right) ^{2}\right\vert _{z=-b}^{z=b} \\
	& =-4\lambda \left( b-\theta \right) \left( u_{z}\left( b\right) +2\lambda
	\left( b-\theta \right) u\left( b\right) \right) ^{2}e^{2\lambda \left(
		b-\theta \right) ^{2}}\geq 0, \\
	& \int_{-b}^{b}\left( -16\lambda ^{3}\left( z-\theta \right) ^{3}\left( 1+%
	\mathcal{O}\left( 1/\lambda \right) \right) v^{2}\right) _{z}dz \\
	& =\left. -16\lambda ^{3}\left( z-\theta \right) ^{3}u^{2}e^{2\lambda \left(
		z-\theta \right) ^{2}}\left( 1+\mathcal{O}\left( 1/\lambda \right) \right)
	\right\vert _{z=-b}^{z=b}\geq 15\lambda ^{3}\left( r-\theta \right)
	^{3}u^{2}\left( b\right) e^{2\lambda \left( b-\theta \right) ^{2}}\geq 0.
\end{align*}

Now we take into account the first term on the right-hand side of %
\eqref{eq:23-1}. Using the Cauchy--Schwarz inequality, we estimate it from
the below as: 
\begin{equation*}
	4\lambda v_{z}^{2}=4\lambda \left[ u_{z}e^{\lambda \left( z-\theta \right)
		^{2}}+2\lambda \left( z-\theta \right) ue^{\lambda \left( z-\theta \right)
		^{2}}\right] ^{2}\geq 4\lambda \left( \frac{1}{2}u_{z}^{2}e^{2\lambda \left(
		z-\theta \right) ^{2}}-4\lambda ^{2}\left( z-\theta \right)
	^{2}u^{2}e^{2\lambda \left( z-\right) ^{2}}\right) .
\end{equation*}%
Hence, using \eqref{eq:23-1}, we obtain 
\begin{equation*}
	\int_{-b}^{b}\left( u_{zz}\right) ^{2}e^{2\lambda \left( z-\theta \right)
		^{2}}dz\geq 2\lambda \int_{-b}^{b}u_{z}^{2}e^{2\lambda \left( z-\theta
		\right) ^{2}}dz+31\lambda ^{3}\left( b-\theta \right)
	^{2}\int_{-b}^{b}u^{2}e^{2\lambda \left( z-\theta \right) ^{2}}dz.
\end{equation*}%
It easily follows from this estimate that 
\begin{equation*}
	\int_{-b}^{b}\left( u_{zz}\right) ^{2}e^{2\lambda \left( z-\theta \right)
		^{2}}dz\geq \frac{1}{2}\int_{-b}^{b}\left( u_{zz}\right) ^{2}e^{2\lambda
		\left( z-\theta \right) ^{2}}dz+\lambda \int_{-b}^{b}u_{z}^{2}e^{2\lambda
		\left( z-\theta \right) ^{2}}dz+15\lambda ^{3}\left( b-\theta \right)
	^{2}\int_{-b}^{b}u^{2}e^{2\lambda \left( z-\theta \right) ^{2}}dz.
\end{equation*}

Hence, we complete the proof of the lemma. $\quad \square$

We now derive a Carleman estimate for the Laplace operator in partial finite
differences. A similar Carleman estimate was proven in Theorem 7.1 of \cite%
{Klibanov2019b}. However, the CWF $\varphi _{\lambda }\left( z\right)
=e^{-2\lambda z}$ in \cite{Klibanov2019b} is different from the one we use.

\begin{thm}[Carleman estimate in partial finite differences]
	\label{thm:carle}There exist a sufficient large constant $\lambda
	_{0}=\lambda _{0}\left( \Omega ,\theta \right) \geq 1$ and a number $%
	C=C\left( \Omega ,r,b\right) >0$ depending only on listed parameters such
	that for all $\lambda \geq \lambda _{0}$ and for all vector functions $%
	u^{h}\in H_{2N,0}^{2,h}\left( \Omega _{h}\right) $ the following Carleman
	estimate holds: 
	\begin{align}
		& \sum_{p,q=1}^{Z_{h}-1}h^{2}\int_{-b}^{b}\left( \Delta
		^{h}u_{p,q}^{h}\left( z\right) \right) ^{2}\mu _{\lambda }\left( z\right) dz
		\label{eq:Carleman} \\
		& \geq C\sum_{p,q=1}^{Z_{h}-1}h^{2}\int_{-b}^{b}\left( \partial
		_{z}^{2}u_{p,q}^{h}\left( z\right) \right) ^{2}\mu _{\lambda }\left(
		z\right) dz+C\lambda \sum_{p,q=1}^{Z_{h}-1}h^{2}\int_{-b}^{b}\left( \partial
		_{z}u_{p,q}^{h}\left( z\right) \right) ^{2}\mu _{\lambda }\left( z\right) dz
		\notag \\
		& +C\lambda ^{3}\sum_{p,q=1}^{Z_{h}-1}h^{2}\int_{-b}^{b}\left( \left( \nabla
		^{h}u_{p,q}^{h}\left( z\right) \right) ^{2}+\left( u_{p,q}^{h}\left(
		z\right) \right) ^{2}\right) \mu _{\lambda }\left( z\right) dz.  \notag
	\end{align}
\end{thm}

\emph{Proof.} Again, it is sufficient to prove \eqref{eq:Carleman} for $%
u\left( x_{p},y_{q},z\right) \in C^{2}\left[ -b,b\right] $ for all $\left(
x_{p},y_{q}\right) .$ Using the definitions of operators in partial finite
differences (subsection \ref{subsec:Preliminaries-on-settings}) we obtain: 
\begin{align*}
	& \sum_{p,q=1}^{Z_{h}-1}h^{2}\int_{-b}^{b}\left\vert \Delta
	^{h}u_{p,q}^{h}\left( z\right) \right\vert ^{2}\mu _{\lambda }\left(
	z\right) dz=\sum_{p,q=1}^{Z_{h}-1}h^{2}\int_{-b}^{b}\left\vert \left(
	u_{zz}^{h}+u_{xx}^{h}+u_{yy}^{h}\right) \left( x_{p},y_{q},z\right)
	\right\vert ^{2}\mu _{\lambda }\left( z\right) dz \\
	& \geq \frac{1}{2}\sum_{p,q=1}^{Z_{h}-1}h^{2}\int_{-b}^{b}\left\vert
	\partial _{z}^{2}u_{p,q}^{h}\left( z\right) \right\vert ^{2}\mu _{\lambda
	}\left( z\right) dz-h^{2}\sum_{p,q=1}^{Z_{h}-1}\int_{-b}^{b}\left\vert
	\left( u_{xx}^{h}+u_{yy}^{h}\right) \left( x_{p},y_{q},z\right) \right\vert
	^{2}\mu _{\lambda }\left( z\right) dz \\
	& \geq \frac{1}{2}\sum_{p,q=1}^{Z_{h}-1}h^{2}\int_{-b}^{b}\left\vert
	\partial _{z}^{2}u_{p,q}^{h}\left( z\right) \right\vert ^{2}\mu _{\lambda
	}\left( z\right) dz-C\sum_{p,q=1}^{Z_{h}-1}\int_{-b}^{b}\left\vert
	u_{p,q}^{h}\left( z\right) \right\vert ^{2}\mu _{\lambda }\left( z\right) dz.
\end{align*}%
Thus, the rest of the proof follows from Lemma \ref{lem:carle}. $\quad
\square $

\subsubsection{Global strict convexity of the functional $J_{h,\protect%
		\lambda }\left( V^{h}\right) $ on the set $\overline{B\left( M\right) }$}

Below, $\left( \cdot ,\cdot \right) $ is the scalar product in the space $%
H_{2N}^{2,h}\left( \Omega _{h}\right) $.

\begin{thm}[Global strict convexity: the central theorem of this paper]
	\label{thm:convex}Let $\lambda _{0}>1$ be the number of Theorem \ref%
	{thm:carle}. For any $\lambda >0$ the functional $J_{h,\lambda }\left(
	V^{h}\right) $ defined in \eqref{eq:J} has its Frechét derivative $%
	J_{h,\lambda }^{\prime }\left( V^{h}\right) \in H_{2N,0}^{2,h}$ at any point 
	$V^{h}\in \overline{B\left( M\right) }$. In addition, there exist numbers $%
	\lambda _{1}=\lambda _{1}\left( \Omega ,\theta ,N,M\right) \geq \lambda
	_{0}>1$ and $C_{2}=C_{2}\left( \Omega _{h},\theta ,N,M\right) >0$ depending
	only on listed parameters such that for all $\lambda \geq \lambda _{1}$ the
	functional $J_{h,\lambda }\left( V^{h}\right) $ is strictly convex on $%
	\overline{B\left( M\right) }$. More precisely, the following estimate holds: 
	\begin{equation}
	J_{h,\lambda }\left( V^{h}+r^{h}\right) -J_{h,\lambda }\left( V^{h}\right)
	-J_{h,\lambda }^{\prime }\left( V^{h}\right) \left( r^{h}\right) \geq
	C_{2}e^{2\lambda \left( b-\theta \right) ^{2}}\left\Vert r^{h}\right\Vert
	_{H_{2N}^{2,h}}^{2}\quad \text{for all }V^{h},V^{h}+r^{h}\in \overline{%
		B\left( M\right) }.  \label{eq:convex}
	\end{equation}
\end{thm}

\emph{Proof.} Denote $\mathbf{x}^{h}=\left( x_{p},y_{q},z\right) $. For
brevity, we do not indicate here the dependence of $\mathbf{x}^{h}$ on
indices $p$ and $q$. Below $C_{2}=C_{2}\left( \Omega _{h},\theta ,N,M\right)
>0$ denotes different positive constants depending only on listed
parameters. Denote 
\begin{equation}
B_{0}\left( M\right) =\left\{ W^{h}\in H_{2N,0}^{2,h}:\left\Vert
W^{h}\right\Vert _{H_{2N}^{2,h}\left( \Omega _{h}\right) }<M\right\} \subset
H_{2N,0}^{2,h}.  \label{99}
\end{equation}%
Let $V_{\left( 1\right) }^{h},V_{\left( 2\right) }^{h}\in \overline{B\left(
	M\right) }$ be two arbitrary points. Denote $r^{h}=\left(
r_{1}^{h},r_{2}^{h}\right) =V_{\left( 2\right) }^{h}-V_{\left( 1\right)
}^{h} $. Then 
\begin{equation}
r^{h}\in \overline{B_{0}\left( 2M\right) }.  \label{100}
\end{equation}%
Obviously, $\left\vert L\left( V_{\left( 2\right) }^{h}\right) \right\vert
^{2}=\left\vert L\left( V_{\left( 1\right) }^{h}+r^{h}\right) \right\vert
^{2}$, where $L$ is the operator defined in \eqref{eq:LL}. Recall that by %
\eqref{eq:17}, the nonlinear term $K\left( \nabla V\right) $ in $L$ is
quadratic with respect to the components of $\nabla V$. Hence, 
\begin{align*}
	L\left( V_{\left( 1\right) }^{h}+r^{h}\right) & =\Delta ^{h}V_{\left(
		1\right) }^{h}+\Delta ^{h}r^{h}+K\left( \nabla ^{h}V_{\left( 1\right)
	}^{h}+\nabla ^{h}r^{h}\right) \\
	& =\Delta ^{h}V_{\left( 1\right) }^{h}+\Delta ^{h}r^{h}+K\left( \nabla
	^{h}V_{\left( 1\right) }^{h}\right) +K_{1}\left( \mathbf{x}^{h}\right)
	\nabla ^{h}r^{h}+K_{2}\left( \mathbf{x}^{h},\nabla ^{h}r^{h}\right) .
\end{align*}%
Thus,%
\begin{equation}
L\left( V_{\left( 1\right) }^{h}+r^{h}\right) =L\left( V_{\left( 1\right)
}^{h}\right) +\Delta ^{h}r^{h}+K_{1}\left( \mathbf{x}^{h}\right) \nabla
^{h}r^{h}+K_{2}\left( \mathbf{x}^{h},\nabla ^{h}r^{h}\right) .
\label{eq:28-1}
\end{equation}%
Here, the vector functions $K_{1},K_{2}$ are continuous with respect to $%
\mathbf{x}^{h}$ in $\overline{\Omega }$. Moreover, $K_{1}\left( \mathbf{x}%
^{h}\right) $ is independent of $r^{h}$ and $K_{2}\left( \mathbf{x}%
^{h},\nabla ^{h}r^{h}\right) $ is quadratic with respect to the components
of $\nabla ^{h}r^{h}$. It follows from \eqref{eq:BM} and the analog of %
\eqref{eq:BMM} for $\overline{B_{0}\left( 2M\right) }$ that 
\begin{equation}
\left\vert K_{2}\left( \mathbf{x}^{h},\nabla ^{h}r^{h}\right) \right\vert
\leq C_{2}\left\vert \nabla ^{h}r^{h}\right\vert ^{2}\quad \text{for all }%
\mathbf{x}\in \overline{\Omega }.  \label{282}
\end{equation}%
Obviously, for any pair $z_{1},z_{2}\in \mathbb{C}$ the following is true: $%
\left\vert z_{1}+z_{2}\right\vert ^{2}=\left\vert z_{1}\right\vert
^{2}+\left\vert z_{2}\right\vert ^{2}+2\func{Re}\left( \overline{z_{1}}%
z_{2}\right) $. Using this, we square the absolute value of both sides of %
\eqref{eq:28-1} and obtain 
\begin{align*}
	\left\vert L^{h}\left( V_{\left( 1\right) }^{h}+r^{h}\right) \right\vert
	^{2}& =\left\vert L^{h}\left( V_{\left( 1\right) }^{h}\right) \right\vert
	^{2}+2\func{Re}\left\{ \overline{L^{h}\left( V_{\left( 1\right) }^{h}\right) 
	}\left[ \Delta ^{h}r^{h}+K_{1}\left( \mathbf{x}^{h}\right) \nabla
	^{h}r^{h}+K_{2}\left( \mathbf{x}^{h},\nabla ^{h}r^{h}\right) \right] \right\}
	\\
	& +\left\vert \Delta ^{h}r^{h}+K_{1}\left( \mathbf{x}^{h}\right) \nabla
	^{h}r^{h}+K_{2}\left( \mathbf{x}^{h},\nabla ^{h}r^{h}\right) \right\vert
	^{2},
\end{align*}%
which leads to 
\begin{align}
	& \left\vert L^{h}\left( V_{\left( 1\right) }^{h}+r^{h}\right) \right\vert
	^{2}-\left\vert L^{h}\left( V_{\left( 1\right) }^{h}\right) \right\vert ^{2}
	\label{eq:29} \\
	& =2\func{Re}\left\{ \overline{L^{h}\left( V_{\left( 1\right) }^{h}\right) }%
	\left[ \Delta ^{h}r^{h}+K_{1}\left( \mathbf{x}^{h}\right) \nabla ^{h}r^{h}%
	\right] \right\} +2\func{Re}\left\{ \overline{L^{h}\left( V_{\left( 1\right)
		}^{h}\right) }K_{2}\left( \mathbf{x}^{h},\nabla ^{h}r^{h}\right) \right\} 
	\notag \\
	& +\left\vert \Delta ^{h}r^{h}\right\vert ^{2}+2\func{Re}\left\{ \overline{%
		\Delta ^{h}r^{h}}\left[ K_{1}\left( \mathbf{x}^{h}\right) \nabla
	^{h}r^{h}+K_{2}\left( \mathbf{x}^{h},\nabla ^{h}r^{h}\right) \right]
	\right\} +\left\vert K_{1}\left( \mathbf{x}^{h}\right) \nabla
	^{h}r^{h}+K_{2}\left( \mathbf{x}^{h},\nabla ^{h}r^{h}\right) \right\vert
	^{2}.  \notag
\end{align}

The first term in the right hand side of equation \eqref{eq:29} is linear
with respect to $r$. Hence, using \eqref{eq:J}, we obtain 
\begin{align}
	& J_{h,\lambda }\left( V_{\left( 1\right) }^{h}+r^{h}\right) -J_{h,\lambda
	}\left( V_{\left( 1\right) }^{h}\right) =\text{Lin}\left( r^{h}\right) 
	\notag \\
	& +\sum_{p,q=1}^{Z_{h}-1}h^{2}\int_{-b}^{b}\left( \left\vert \Delta
	^{h}r^{h}\right\vert ^{2}+2\func{Re}\left\{ \overline{\Delta ^{h}r}\left[
	\left( K_{1}\left( \mathbf{x}^{h}\right) \nabla ^{h}r^{h}\right)
	+K_{2}\left( \mathbf{x}^{h},\nabla ^{h}r^{h}\right) \right] \right\} \right)
	\mu _{\lambda }\left( z\right) dz  \notag \\
	& +\sum_{p,q=1}^{Z_{h}-1}h^{2}\int_{-b}^{b}\left[ 2\func{Re}\left\{ 
	\overline{L\left( V_{\left( 1\right) }^{h}\right) }K_{2}\left( \mathbf{x}%
	^{h},\nabla ^{h}r^{h}\right) \right\} +\left\vert K_{1}\left( \mathbf{x}%
	^{h}\right) \nabla ^{h}r^{h}+K_{2}\left( \mathbf{x}^{h},\nabla
	^{h}r^{h}\right) \right\vert ^{2}\right] \mu _{\lambda }\left( z\right) dz.
	\label{eq:J-J}
\end{align}%
In \eqref{eq:J-J}, the linear functional $\text{Lin}\left( r^{h}\right)
:H_{2N,0}^{2,h}\rightarrow \mathbb{R}$ is given by 
\begin{equation*}
	\text{Lin}\left( r\right) =2\sum_{p,q=1}^{Z_{h}-1}h^{2}\int_{-b}^{b}\func{Re}%
	\left\{ \overline{L\left( V_{\left( 1\right) }^{h}\right) }\left[ \Delta
	^{h}r+K_{1}\left( \mathbf{x}^{h}\right) \nabla ^{h}r\right] \right\} \mu
	_{\lambda }\left( z\right) dz.
\end{equation*}%
We now estimate this functional from the above. Since $V_{\left( 1\right)
}^{h}\in \overline{B\left( M\right) }$ then the structure of $K$ and %
\eqref{eq:BMM} imply that $\left\vert K\left( \nabla V_{\left( 1\right)
}^{h}\right) \right\vert \leq C_{2}$. Hence, the H\H{o}lder inequality and (%
\ref{max}) imply that 
\begin{align*}
	\left\vert \text{Lin}\left( r\right) \right\vert & \leq
	2\sum_{p,q=1}^{Z_{h}-1}h^{2}\left( \int_{-b}^{b}\left\vert \Delta V_{\left(
		1\right) }^{h}+K\left( \nabla V_{\left( 1\right) }^{h}\right) \right\vert
	^{2}dz\right) ^{1/2}\left( \int_{-b}^{b}\left\vert \Delta
	^{h}r^{h}+K_{1}\left( \mathbf{x}^{h}\right) \nabla ^{h}r^{h}\right\vert
	^{2}dz\right) ^{1/2} \\
	& \leq 4C_{2}e^{2\lambda \left( b+\theta \right) ^{2}}\left\Vert
	r^{h}\right\Vert _{H_{2N}^{2,h}}.
\end{align*}%
Thus, the functional $\text{Lin}\left( r^{h}\right)
:H_{2N,0}^{2,h}\rightarrow \mathbb{R}$ is linear and bounded. Hence, by the
Riesz theorem, there exists a unique point $\mathbf{P}^{h}\in H_{2N,0}^{2,h}$
independent of $r^{h}$ such that 
\begin{equation*}
	\text{Lin}\left( r^{h}\right) =\left( \mathbf{P}^{h},r^{h}\right) \quad 
	\text{for all }r^{h}\in H_{2N,0}^{2,h}.
\end{equation*}

Next, using the Cauchy--Schwarz inequality and (\ref{282}), we find that 
\begin{equation}
2\func{Re}\left\{ \overline{\Delta ^{h}r^{h}}\left[ \left( K_{1}\left( 
\mathbf{x}^{h}\right) \nabla ^{h}r^{h}\right) +K_{2}\left( \mathbf{x}%
^{h},\nabla ^{h}r^{h}\right) \right] \right\} \leq \frac{1}{2}\left\vert
\Delta ^{h}r^{h}\right\vert ^{2}+C_{2}\left\vert \nabla ^{h}r^{h}\right\vert
^{2}.  \label{eq:31}
\end{equation}%
In addition, we have

\begin{equation}
\left\vert 2\func{Re}\left\{ \overline{L^{h}\left( V_{\left( 1\right)
	}^{h}\right) }K_{2}\left( \mathbf{x}^{h},\nabla ^{h}r^{h}\right) \right\}
\right\vert +\left\vert K_{1}\left( \mathbf{x}^{h}\right) \nabla
^{h}r^{h}+K_{2}\left( \mathbf{x}^{h},\nabla ^{h}r^{h}\right) \right\vert
^{2}\leq C_{2}\left\vert \nabla ^{h}r^{h}\right\vert ^{2}.  \label{eq:32}
\end{equation}%
Combining \eqref{eq:J-J}, \eqref{eq:31} and \eqref{eq:32}, we easily obtain 
\begin{equation*}
	\left\vert J_{h,\lambda }\left( V_{\left( 1\right) }^{h}+r\right)
	-J_{h,\lambda }\left( V_{\left( 1\right) }^{h}\right) -\left( \mathbf{P}%
	^{h},r\right) \right\vert \leq C_{2}\left\Vert r^{h}\right\Vert
	_{H_{2N}^{2,h}}^{2}.
\end{equation*}%
Henceforth, $\mathbf{P}^{h}\in H_{2N,0}^{2,h}$ is actually the Frechét
derivative of the cost functional $J_{h,\lambda }$ at the point $V_{\left(
	1\right) }^{h}\in \overline{B\left( M\right) }$, i.e. $\mathbf{P}%
^{h}=J_{h,\lambda }^{\prime }\left( V_{\left( 1\right) }^{h}\right) \in
H_{2N,0}^{2,h}$.

We now focus on the proof of the target estimate \eqref{eq:convex}. To do
so, we estimate from below the second and third terms on the right hand side
of \eqref{eq:J-J}. In fact, using \eqref{eq:31} we get 
\begin{align}
	& \sum_{p,q=1}^{Z_{h}-1}h^{2}\int_{-b}^{b}\left( \left\vert \Delta
	^{h}r^{h}\right\vert ^{2}+2\func{Re}\left\{ \overline{\Delta ^{h}r^{h}}\left[
	\left( K_{1}\left( \mathbf{x}^{h}\right) \nabla ^{h}r^{h}\right)
	+K_{2}\left( \mathbf{x}^{h},\nabla ^{h}r^{h}\right) \right] \right\} \right)
	\mu _{\lambda }\left( z\right) dz  \notag \\
	& \geq \frac{1}{2}\sum_{p,q=1}^{Z_{h}-1}h^{2}\int_{-b}^{b}\left\vert \Delta
	^{h}r^{h}\right\vert ^{2}\mu _{\lambda }\left( z\right)
	dz-C_{2}\sum_{p,q=1}^{Z_{h}-1}h^{2}\int_{-b}^{b}\left\vert \nabla
	^{h}r^{h}\right\vert ^{2}\mu _{\lambda }\left( z\right) dz.  \label{34}
\end{align}%
Similarly, using (\ref{282}), we estimate the third term in the right hand
side of \eqref{eq:J-J} as 
\begin{align}
	& \sum_{p,q=1}^{Z_{h}-1}h^{2}\int_{-b}^{b}\left[ 2\func{Re}\left\{ \overline{%
		L^{h}\left( V_{\left( 1\right) }^{h}\right) }K_{2}\left( \mathbf{x}%
	^{h},\nabla ^{h}r^{h}\right) \right\} +\left\vert K_{1}\left( \mathbf{x}%
	^{h}\right) \nabla ^{h}r^{h}+K_{2}\left( \mathbf{x}^{h},\nabla
	^{h}r^{h}\right) \right\vert ^{2}\right] \mu _{\lambda }\left( z\right) dz 
	\notag \\
	& \geq -C_{2}\sum_{p,q=1}^{Z_{h}-1}h^{2}\int_{-b}^{b}\left\vert \nabla
	^{h}r^{h}\right\vert ^{2}\mu _{\lambda }\left( z\right) dz.  \label{eq:35}
\end{align}%
Thus, combining \eqref{eq:J-J} and \eqref{eq:35}, we obtain 
\begin{align}
	& J_{h,\lambda }\left( V_{\left( 1\right) }^{h}+r^{h}\right) -J_{h,\lambda
	}\left( V_{\left( 1\right) }^{h}\right) -\left( J_{h,\lambda }^{\prime
	}\left( V_{\left( 1\right) }^{h}\right) ,r^{h}\right)  \notag \\
	& \geq \left[ \frac{1}{2}\sum_{p,q=1}^{Z_{h}-1}h^{2}\int_{-b}^{b}\left\vert
	\Delta ^{h}r_{p,q}^{h}\right\vert ^{2}\mu _{\lambda }\left( z\right)
	dz-C_{2}\sum_{p,q=1}^{Z_{h}-1}h^{2}\int_{-b}^{b}\left\vert \nabla
	^{h}r_{p,q}^{h}\right\vert ^{2}\mu _{\lambda }\left( z\right) dz\right] .
	\label{eq:36}
\end{align}%
Since the function $r\in H_{2N,0}^{2,h}$, we now can apply to \eqref{eq:36}
the Carleman estimate \eqref{eq:Carleman}. Prior to that, we note that we
can find a sufficiently large number $\tilde{\lambda}_{0}=\tilde{\lambda}%
_{0}\left( \Omega ,\theta ,N,M\right) \geq \lambda _{0}>1$ such that for all 
$\lambda \geq \tilde{\lambda}_{0}$ 
\begin{align*}
	& \frac{1}{2}\sum_{p,q=1}^{Z_{h}-1}h^{2}\int_{-b}^{b}\left\vert \Delta
	^{h}r_{p,q}^{h}\right\vert ^{2}\mu _{\lambda }\left( z\right)
	dz-C_{2}\sum_{p,q=1}^{Z_{h}-1}h^{2}\int_{-b}^{b}\left\vert \nabla
	^{h}r_{p,q}^{h}\right\vert ^{2}\mu _{\lambda }\left( z\right) dz \\
	& \geq C\sum_{p,q=1}^{Z_{h}-1}h^{2}\int_{-b}^{b}\left\vert \partial
	_{z}^{2}r_{p,q}^{h}\left( z\right) \right\vert ^{2}\mu _{\lambda }\left(
	z\right) dz+C\lambda \sum_{p,q=1}^{Z_{h}-1}h^{2}\int_{-b}^{b}\left\vert
	\partial _{z}r_{p,q}^{h}\left( z\right) \right\vert ^{2}\mu _{\lambda
	}\left( z\right) dz \\
	& +C\lambda ^{3}\sum_{p,q=1}^{Z_{h}-1}h^{2}\int_{-b}^{b}\left( \left\vert
	\nabla ^{h}r_{p,q}^{h}\left( z\right) \right\vert ^{2}+\left\vert
	r_{p,q}^{h}\left( z\right) \right\vert ^{2}\right) \mu _{\lambda }\left(
	z\right) dz-C_{2}\sum_{p,q=1}^{Z_{h}-1}h^{2}\int_{-b}^{b}\left\vert \nabla
	^{h}r_{p,q}^{h}\right\vert ^{2}\mu _{\lambda }\left( z\right) dz.
\end{align*}%
Hence, choosing $\lambda _{1}=\lambda _{1}\left( \Omega _{h},\theta
,N,M\right) \geq \tilde{\lambda}_{0}>1$ such that $C\lambda _{1}>2C_{2}$, we
obtain 
\begin{align}
	& \frac{1}{2}\sum_{p,q=1}^{Z_{h}-1}h^{2}\int_{-b}^{b}\left\vert \Delta
	^{h}r_{p,q}^{h}\right\vert ^{2}\mu _{\lambda _{1}}\left( z\right)
	dz-C_{2}\sum_{p,q=1}^{Z_{h}-1}h^{2}\int_{-b}^{b}\left\vert \nabla
	^{h}r_{p,q}^{h}\right\vert ^{2}\mu _{\lambda _{1}}\left( z\right) dz
	\label{eq:37} \\
	& \geq C_{2}\sum_{p,q=1}^{Z_{h}-1}h^{2}\int_{-b}^{b}\left\vert \partial
	_{z}^{2}r_{p,q}^{h}\left( z\right) \right\vert ^{2}\mu _{\lambda _{1}}\left(
	z\right) dz+C_{2}\lambda \sum_{p,q=1}^{Z_{h}-1}h^{2}\int_{-b}^{b}\left(
	\left\vert \partial _{z}r_{p,q}^{h}\left( z\right) \right\vert
	^{2}+\left\vert r^{h}\left( z\right) \right\vert ^{2}\right) \mu _{\lambda
		_{1}}\left( z\right) dz.  \notag
\end{align}%
Thus, it follows from \eqref{eq:36},\ \eqref{eq:37} and (\ref{max1}) that
for all $\lambda \geq \lambda _{1}$ 
\begin{equation}
J_{h,\lambda }\left( V_{\left( 1\right) }^{h}+r^{h}\right) -J_{h,\lambda
}\left( V_{\left( 1\right) }^{h}\right) -\left( J_{h,\lambda }^{\prime
}\left( V_{\left( 1\right) }^{h}\right) ,r^{h}\right) \geq C_{2}e^{2\lambda
	\left( b-\theta \right) ^{2}}\left\Vert r^{h}\right\Vert _{H_{2N}^{2,h}}^{2}.
\label{4}
\end{equation}%
Estimate (\ref{4}) completes the proof of this theorem. $\quad \text{ }%
\square $

We now formulate a theorem about the Lipschitz continuity of the Frechét
derivative $J_{h,\lambda }^{\prime }\left( V^{h}\right) $ on $\overline{%
	B\left( M\right) }$. We omit the proof of this result because it is similar
to the proof of Theorem 3.1 in \cite{Bakushinsky2017}.

\begin{thm}
	\label{thm:lip}The Frechét derivative $J_{h,\lambda }^{\prime }\left(
	V^{h}\right) $ constructed in the proof of Theorem \ref{thm:convex}
	satisfies the Lipschitz continuity condition on the set $\overline{B\left(
		M\right) }$. More precisely, there exists a number $\widehat{C}=\widehat{C}%
	\left( \Omega ,\theta ,N,M,\lambda \right) >0$ depending only on listed
	parameters such that for any pair $V_{\left( 1\right) }^{h},V_{\left(
		2\right) }^{h}\in \overline{B\left( M\right) }$ the following estimate
	holds: 
	\begin{equation*}
		\left\Vert J_{h,\lambda }^{\prime }\left( V_{\left( 2\right) }^{h}\right)
		-J_{h,\lambda }^{\prime }\left( V_{\left( 1\right) }^{h}\right) \right\Vert
		_{H_{2N}^{2,h}}\leq \tilde{C}\left\Vert V_{\left( 2\right) }^{h}-V_{\left(
			1\right) }^{h}\right\Vert _{H_{2N}^{2,h}}.
	\end{equation*}
\end{thm}

As to the existence and uniqueness of the minimizer, they are established in
the following theorem. This theorem follows from a combination of Theorems %
\ref{thm:convex} and \ref{thm:lip} with Lemma 2.1 and Theorem 2.1 of \cite%
{Bakushinsky2017}. Therefore, we omit its proof.

\begin{thm}
	\label{thm:6} Let $\lambda _{1}>1$ be the number chosen in Theorem \ref%
	{thm:convex}. Then there exists a unique minimizer $V_{\min ,\lambda
	}^{h}\in \overline{B\left( M\right) }$ of the functional $J_{h,\lambda
	}^{\prime }\left( V^{h}\right) $ on the set $\overline{B\left( M\right) }$.
	Furthermore, the following inequality holds: 
	\begin{equation}
	\left( J_{h,\lambda }^{\prime }\left( V_{\min ,\lambda }^{h}\right) ,V_{\min
		,\lambda }^{h}-Q\right) \leq 0\quad \text{for all } Q\in \overline{B\left(
		M\right) }.  \label{3.0}
	\end{equation}
\end{thm}

\subsection{Convergence rate of regularized solutions}

Using (\ref{98}), we obtain the following analog of problem (\ref{eq:pdeV})--%
\eqref{eq:bdrV2} in partial finite differences is:%
\begin{align}
	& L^{h}\left( V^{h}\left( z\right) \right) =\Delta ^{h}V^{h}\left( \mathbf{x}%
	^{h}\right) +K\left( \nabla ^{h}V^{h}\left( \mathbf{x}^{h}\right) \right)
	=0\quad \text{for }\mathbf{x}^{h}\in \Omega _{h},  \label{3.1} \\
	& \nabla ^{h}V^{h}\left( \mathbf{x}^{h}\right) \cdot \text{n}=0\quad \text{%
		for }\mathbf{x}^{h}\in \partial \Omega _{h}\backslash \Gamma _{h},
	\label{3.2} \\
	& V\left( \mathbf{x}^{h}\right) =\psi _{0}^{h}\left( \mathbf{x}^{h}\right)
	,V_{z}\left( \mathbf{x}^{h}\right) =\psi _{1}\left( \mathbf{x}^{h}\right)
	\quad \text{for }\mathbf{x}^{h}\in \Gamma _{h}.  \label{3.3}
\end{align}%
Following the Tikhonov regularization concept (cf. e.g. \cite{Tikhonov1995}%
), we assume that there exists an exact solution $V_{\ast }^{h}\in $ $%
H_{2N}^{2,h}$ of problem (\ref{3.1})--(\ref{3.3}) with the noiseless data $%
\psi _{0\ast }^{h}\left( \mathbf{x}^{h}\right) $ and $\psi _{1\ast
}^{h}\left( \mathbf{x}^{h}\right) $. To this end, the subscript
\textquotedblleft $\ast $\textquotedblright\ is only used for the exact
solution. In applications the data $\psi _{0}^{h}\left( \mathbf{x}%
^{h}\right) $ and $\psi _{1}^{h}\left( \mathbf{x}^{h}\right) $ are
noise-contaminated.

We, therefore, denote by $\delta \in \left( 0,1\right) $ the level of noise
in those data. We assume that 
\begin{equation}
\left\Vert V_{\ast }^{h}\right\Vert _{H_{2N}^{2,h}}<M-\delta .  \label{3.40}
\end{equation}%
In accordance with the regularization theory (cf. e.g. \cite{Tikhonov1995}),
given the value of $\delta ,$ the minimizer $V_{\min ,\lambda }^{h}$ of the
functional $J_{h,\lambda }\left( V^{h}\right) ,$ which was found in Theorem %
\ref{thm:6}, is called the \emph{regularized solution} of problem (\ref{3.1}%
)--(\ref{3.3}). We want to estimate the $\delta -$dependence of the norm $%
\left\Vert V_{\ast }^{h}-V_{\min ,\lambda }^{h}\right\Vert _{H_{2N}^{2,h}}$,
i.e. we want to estimate the convergence rate of regularized solutions,
which also means the accuracy estimate of the minimizer $V_{\min ,\lambda
}^{h}.$ To do this, we assume that there exist two vector functions $\Psi
_{\ast }^{h},\Psi ^{h}\in H_{2N}^{2,h}$ such that 
\begin{align}
	& \nabla ^{h}\Psi _{\ast }^{h}\left( \mathbf{x}^{h}\right) \cdot \text{n}%
	=0,\quad \nabla ^{h}\Psi ^{h}\left( \mathbf{x}^{h}\right) \cdot \text{n}%
	=0\quad \text{for }\mathbf{x}^{h}\in \partial \Omega _{h}\backslash \Gamma
	_{h},  \label{3.4} \\
	& \Psi _{\ast }^{h}\left( \mathbf{x}^{h}\right) =\psi _{0\ast }^{h}\left( 
	\mathbf{x}^{h}\right) ,\quad \Psi _{z\ast }^{h}\left( \mathbf{x}^{h}\right)
	=\psi _{1\ast }\left( \mathbf{x}^{h}\right) \quad \text{for }\mathbf{x}%
	^{h}\in \Gamma _{h},  \label{3.5} \\
	& \Psi ^{h}\left( \mathbf{x}^{h}\right) =\psi _{0}^{h}\left( \mathbf{x}%
	^{h}\right) ,\quad \Psi _{z}^{h}\left( \mathbf{x}^{h}\right) =\psi
	_{1}\left( \mathbf{x}^{h}\right) \quad \text{for }\mathbf{x}^{h}\in \Gamma
	_{h},  \label{3.6} \\
	& \left\Vert \Psi _{\ast }^{h}\right\Vert _{H_{2N}^{2,h}}<M,\quad \left\Vert
	\Psi ^{h}\right\Vert _{H_{2N}^{2,h}}<M,  \label{3.7} \\
	& \left\Vert \Psi ^{h}-\Psi _{\ast }^{h}\right\Vert _{H_{2N}^{2,h}}<\delta .
	\label{3.8}
\end{align}

\begin{thm}[convergence rate of regularized solutions]
	\label{thm:7} Assume that conditions (\ref{3.40})--(\ref{3.8}) are valid.
	Let $\lambda _{1}=\lambda _{1}\left( \Omega ,\theta ,N,M\right) >1$ be the
	number of Theorem \ref{thm:lip}. Let $V_{\min ,\lambda }^{h}\in \overline{%
		B\left( M\right) }$ be the minimizer of functional (\ref{eq:J}) which is
	found in Theorem \ref{thm:6}. Then the following accuracy estimate holds for
	all $\lambda \geq \lambda _{1}$ 
	\begin{equation}
	\left\Vert V_{\min ,\lambda }^{h}-V_{\ast }^{h}\right\Vert
	_{H_{2N}^{2,h}}\leq C_{2}\delta e^{4\lambda b\theta }.  \label{3.80}
	\end{equation}
\end{thm}

\emph{Proof.} Denote 
\begin{equation}
W_{\min ,\lambda }^{h}=V_{\min ,\lambda }^{h}-\Psi ^{h},\quad W_{\ast
}^{h}=V_{\ast }^{h}-\Psi _{\ast }^{h}.  \label{3.9}
\end{equation}%
Hence, by (\ref{99}) and (\ref{3.7}) 
\begin{equation}
W_{\min ,\lambda }^{h},W_{\ast }^{h}\in B_{0}\left( 2M\right) .  \label{101}
\end{equation}%
It follows from (\ref{3.40})--(\ref{3.9}) and the triangle inequality that 
\begin{eqnarray*}
	\left\Vert W_{\ast }^{h}+\Psi ^{h}\right\Vert _{H_{2N}^{2,h}} &=&\left\Vert
	W_{\ast }^{h}+\Psi _{\ast }^{h}+\left( \Psi ^{h}-\Psi _{\ast }^{h}\right)
	\right\Vert _{H_{2N}^{2,h}}\leq \left\Vert W_{\ast }^{h}+\Psi _{\ast
	}^{h}\right\Vert _{H_{2N}^{2,h}}+\delta \\
	&<&M-\delta +\delta =M.
\end{eqnarray*}%
Denote 
\begin{equation}
\widetilde{V}_{\ast }^{h}=W_{\ast }^{h}+\Psi ^{h}.  \label{3.130}
\end{equation}%
Hence, (\ref{eq:BM}), (\ref{3.5}), (\ref{3.6}), (\ref{101}) and (\ref{3.130}%
) imply that%
\begin{equation}
W_{\ast }^{h}+\Psi ^{h}=\widetilde{V}_{\ast }^{h}\in B\left( M\right) .
\label{3.12}
\end{equation}

By (\ref{3.1}), it holds that $L^{h}\left( V_{\ast }^{h}\right) =0.$ Hence, $%
L^{h}\left( W_{\ast }^{h}+\Psi _{\ast }^{h}\right) =L^{h}\left( V_{\ast
}^{h}\right) =0.$ Hence, by (\ref{eq:J}) we have 
\begin{equation}
J_{h,\lambda }\left( V_{\ast }^{h}\right) =0.  \label{3.13}
\end{equation}%
Next, by (\ref{eq:convex}) and (\ref{3.12}) we estimate that 
\begin{equation}\label{3.14}
	J_{h,\lambda }\left( \widetilde{V}_{\ast }^{h}\right) -J_{h,\lambda }\left(
	V_{\min ,\lambda }^{h}\right) -J_{h,\lambda }^{\prime }\left( V_{\min
		,\lambda }^{h}\right) \left( \widetilde{V}_{\ast }^{h}-V_{\min ,\lambda
	}^{h}\right) \geq C_{2}e^{2\lambda \left( b-\theta \right) ^{2}}\left\Vert 
	\widetilde{V}_{\ast }^{h}-V_{\min ,\lambda }^{h}\right\Vert
	_{H_{2N}^{2,h}}^{2}.
\end{equation}%
It follows from (\ref{3.0}) that $-J_{h,\lambda }^{\prime }\left( V_{\min
	,\lambda }^{h}\right) \left( \widetilde{V}_{\ast }^{h}-V_{\min ,\lambda
}^{h}\right) \leq 0.$ This implies that%
\begin{equation*}
	J_{h,\lambda }\left( \widetilde{V}_{\ast }^{h}\right) -J_{h,\lambda }\left(
	V_{\min ,\lambda }^{h}\right) -J_{h,\lambda }^{\prime }\left( V_{\min
		,\lambda }^{h}\right) \left( \widetilde{V}_{\ast }^{h}-V_{\min ,\lambda
	}^{h}\right) \leq J_{h,\lambda }\left( \widetilde{V}_{\ast }^{h}\right) .
\end{equation*}%
The latter and (\ref{3.14}) lead to%
\begin{equation}
\left\Vert \widetilde{V}_{\ast }^{h}-V_{\min ,\lambda }^{h}\right\Vert
_{H_{2N}^{2,h}}^{2}\leq C_{2}e^{-2\lambda \left( b-\theta \right)
	^{2}}J_{h,\lambda }\left( \widetilde{V}_{\ast }^{h}\right) .  \label{3.15}
\end{equation}

We now estimate $J_{h,\lambda }\left( \widetilde{V}_{\ast }^{h}\right) $
from the above. By (\ref{eq:J}), (\ref{3.1}), (\ref{3.8}), (\ref{3.13}) and (%
\ref{3.130}) it yields 
\begin{eqnarray*}
	J_{h,\lambda }\left( \widetilde{V}_{\ast }^{h}\right)
	&=&\sum_{p,q=1}^{Z_{h}-1}h^{2}\int_{-b}^{b}\left\vert L^{h}\left( \widetilde{%
		V}_{\ast }^{h}\left( z\right) \right) \right\vert ^{2}\mu _{\lambda }\left(
	z\right) dz \\
	&=&\sum_{p,q=1}^{Z_{h}-1}h^{2}\int_{-b}^{b}\left\vert L^{h}\left( W_{\ast
	}^{h}\left( z\right) +\Psi _{\ast }^{h}\left( z\right) \right) +\left( \Psi
	^{h}\left( z\right) -\Psi _{\ast }^{h}\left( z\right) \right) \right\vert
	^{2}\mu _{\lambda }\left( z\right) dz \\
	&=&\sum_{p,q=1}^{Z_{h}-1}h^{2}\int_{-b}^{b}\left\vert L^{h}\left( V_{\ast
	}^{h}\left( z\right) \right) \right\vert ^{2}\mu _{\lambda }\left( z\right)
	dz+\sum_{p,q=1}^{Z_{h}-1}h^{2}\int_{-b}^{b}S^{h}\left( z\right) \mu
	_{\lambda }\left( z\right) dz \\
	&=&\sum_{p,q=1}^{Z_{h}-1}h^{2}\int_{-b}^{b}S^{h}\left( z\right) \mu
	_{\lambda }\left( z\right) dz.
\end{eqnarray*}%
In view of the fact that 
\begin{equation*}
	\sum_{p,q=1}^{Z_{h}-1}h^{2}\int_{-b}^{b}\left\vert S^{h}\left( z\right)
	\right\vert \mu _{\lambda }\left( z\right) dz\leq C_{2}\delta ^{2}\max_{ 
		\left[ -b,b\right] }\mu _{\lambda }\left( z\right) ,
\end{equation*}%
we then use (\ref{max}) to obtain%
\begin{equation}
J_{h,\lambda }\left( \widetilde{V}_{\ast }^{h}\right) \leq
\sum_{p,q=1}^{Z_{h}-1}h^{2}\int_{-b}^{b}\left\vert S^{h}\left( z\right)
\right\vert \mu _{\lambda }\left( z\right) dz\leq C_{2}\delta ^{2}\exp
\left( 2\lambda \left( b+\theta \right) ^{2}\right) .  \label{3.16}
\end{equation}%
Combining estimates (\ref{3.15}) and (\ref{3.16}), we obtain the target
estimate (\ref{3.80}) of this theorem. $\quad \square $

\subsection{Global convergence of the gradient projection method}

Just as in the previous section, we still assume the existence of vector
functions $V_{\ast }^{h},\Psi _{\ast }^{h}$ and $\Psi ^{h}$ satisfying
conditions formulated in that section. Similarly with (\ref{3.9}), for each $%
V^{h}\in B\left( M\right) $ consider the vector function $W^{h}=V^{h}-\Psi
^{h}.$ Then (\ref{3.7}) and the triangle inequality imply that, similarly
with (\ref{101}), 
\begin{align}
	& W^{h}\in B_{0}\left( 2M\right) \subset H_{2N,0}^{2,h}\quad \text{for all }%
	V^{h}\in B\left( M\right) ,  \label{3.17} \\
	& W^{h}+\Psi ^{h}\in B\left( 3M\right) \quad \text{for all }W^{h}\in
	B_{0}\left( 2M\right) .  \label{3.18}
\end{align}%
Consider the functional $I_{h,\lambda }:B_{0}\left( 2M\right) \rightarrow 
\mathbb{R}$ defined as%
\begin{equation}
I_{h,\lambda }\left( W^{h}\right) =J_{h,\lambda }\left( W^{h}+\Psi
^{h}\right) \quad \text{for all }W^{h}\in B_{0}\left( 2M\right) .
\label{3.19}
\end{equation}%
We omit the proof of Theorem \ref{thm:8} since it follows immediately from
Theorems \ref{thm:lip}--\ref{thm:7} and (\ref{3.17})--(\ref{3.19}).

\begin{thm}
	\label{thm:8} For any $\lambda >0$ the functional $I_{h,\lambda }\left(
	W^{h}\right) $ has its Frechét derivative $I_{h,\lambda }^{\prime }\left(
	W^{h}\right) \in H_{2N,0}^{2,h}$ at any point $W^{h}\in \overline{%
		B_{0}\left( 2M\right) }$ and this derivative is Lipschitz continuous on $%
	\overline{B_{0}\left( 2M\right) }.$ Let $\lambda _{1}=\lambda _{1}\left(
	\Omega ,\theta ,N,M\right) >1$ and $C_{2}=C_{2}\left( \Omega _{h},\theta
	,N,M\right) >0$ be the numbers of Theorem \ref{thm:lip}. Denote $\widetilde{%
		\lambda }=\lambda _{1}\left( \Omega ,\theta ,N,3M\right) >1$ and $\widetilde{%
		C}_{2}=C_{2}\left( \Omega _{h},\theta ,N,3M\right) >0$. Then for any $%
	\lambda \geq \widetilde{\lambda }$ the functional $I_{h,\lambda }\left(
	W^{h}\right) $ is strictly convex on the ball $B_{0}\left( 2M\right) \subset
	H_{2N,0}^{2,h},$ i.e. the following analog of estimate (\ref{eq:convex}) is
	valid for all $W^{h},W^{h}+r^{h}\in \overline{B_{0}\left( 2M\right) }:$ 
	\begin{equation*}
		I_{h,\lambda }\left( W^{h}+r^{h}\right) -I_{h,\lambda }\left( W^{h}\right)
		-I_{h,\lambda }^{\prime }\left( W^{h}\right) \left( r^{h}\right) \geq 
		\widetilde{C}_{2}e^{2\lambda \left( b-\theta \right) ^{2}}\left\Vert
		r^{h}\right\Vert _{H_{2N}^{2,h}}^{2}.
	\end{equation*}%
	Furthermore, there exists a unique minimizer $W_{\min ,\lambda }^{h}$ of the
	functional $I_{h,\lambda }\left( W^{h}\right) $ on the set $\overline{%
		B_{0}\left( 2M\right) }$ and the following inequality holds: 
	\begin{equation*}
		\left( I_{h,\lambda }^{\prime }\left( W_{\min ,\lambda }^{h}\right) ,W_{\min
			,\lambda }^{h}-Q\right) \leq 0\quad \text{for all }Q\in \overline{%
			B_{0}\left( 2M\right) }.
	\end{equation*}%
	On top of that, let $Y_{\min ,\lambda }^{h}=W_{\min ,\lambda }^{h}+\Psi
	^{h}. $ Then the direct analog of (\ref{3.80}) holds where $V_{\min ,\lambda
	}^{h}$ is replaced with $Y_{\min ,\lambda }^{h}$.
\end{thm}

We now construct the gradient projection method of the minimization of the
functional $I_{h,\lambda }\left( W^{h}\right) $ on the set $\overline{%
	B_{0}\left( 2M\right) }.$ Let $P:H_{2N,0}^{2,h}\rightarrow \overline{%
	B_{0}\left( 2M\right) }$ be the orthogonal projection operator. Let $%
W_{0}^{h}\in B_{0}\left( 2M\right) $ be an arbitrary point of this ball. Let 
$\gamma \in \left( 0,1\right) $ be a number, which we will choose in Theorem
\ref{thm:9}. The sequence of the gradient projection method is:%
\begin{equation}
W_{n,\lambda ,\gamma }^{h}=\mathbf{P}\left( W_{n-1,\lambda ,\gamma }^{h}-\gamma
I_{h,\lambda }^{\prime }\left( W_{n-1,\lambda ,\gamma }^{h}\right) \right) ,%
\text{ }n=1,2,...  \label{3.20}
\end{equation}%
Note that since by Theorem \ref{thm:8}, we have $I_{h,\lambda }^{\prime
}\left( W_{n-1,\lambda ,\gamma }^{h}\right) \in H_{2N,0}^{2,h}$ and also
since $W_{n-1,\lambda ,\gamma }^{h}\in \overline{B_{0}\left( 2M\right) }%
\subset H_{2N,0}^{2,h}$, then all three terms in (\ref{3.20}) belong to $%
H_{2N,0}^{2,h},$ i.e. zero boundary conditions (\ref{3.2}), (\ref{3.3}) are
satisfied for these functions.

\begin{thm}[the global convergence of the gradient projection method]
	\label{thm:9} Assume that conditions of Theorem \ref{thm:8} hold and let $%
	\lambda \geq \widetilde{\lambda }$. Then there exists a number $\gamma
	_{0}=\gamma _{0}\left( \Omega _{h},\theta ,N,3M\right) \in \left( 0,1\right) 
	$ depending only on listed parameters such that for every $\gamma \in \left(
	0,\gamma _{0}\right) $ there exists a number $\xi =\xi \left( \gamma \right)
	\in \left( 0,\gamma _{0}\right) $ depending on $\gamma $ such that for these
	values of $\gamma $ the sequence (\ref{3.20}) converges to $W_{\min ,\lambda
	}^{h}$ and the following convergence rate holds: 
	\begin{equation}
	\left\Vert W_{n,\lambda ,\gamma }^{h}-W_{\min ,\lambda }^{h}\right\Vert
	_{H_{2N}^{2,h}}\leq \xi ^{n}\left\Vert W_{\min ,\lambda
	}^{h}-W_{0}^{h}\right\Vert _{H_{2N}^{2,h}}.  \label{3.21}
	\end{equation}%
	In addition, 
	\begin{equation}
	\left\Vert W_{\ast }^{h}-W_{n,\lambda ,\gamma }^{h}\right\Vert
	_{H_{2N}^{2,h}}\leq \widetilde{C}_{2}\delta e^{4\lambda b\theta }+\xi
	^{n}\left\Vert W_{\min ,\lambda }^{h}-W_{0}^{h}\right\Vert _{H_{2N}^{2,h}}.
	\label{3.22}
	\end{equation}%
	Now, let $c^{h}\left( \mathbf{x}^{h}\right) $ be the function $c\left( 
	\mathbf{x}\right) $ obtained after the substitution of components of the
	vector function $V^{h}$ in equation (\ref{eq:v}) for a certain position of
	the point source $\mathbf{x}_{\alpha }$ (Remark \ref{rem:X}). More
	precisely, let $c_{n,\lambda ,\gamma }^{h}\left( \mathbf{x}^{h}\right) $ be
	obtained after the substitution of the components of the vector function $%
	V_{n,\lambda ,\gamma }^{h}=W_{n,\lambda ,\gamma }^{h}+\Psi ^{h}$ and let $%
	c_{\ast }^{h}\left( \mathbf{x}^{h}\right) $ be obtained after the
	substitution of the components of the vector function $V_{\ast }^{h}=W_{\ast
	}^{h}+\Psi _{\ast }^{h}.$ Then the following convergence estimate is valid 
	\begin{equation}
	\left\Vert c_{\ast }^{h}-c_{n,\lambda ,\gamma }^{h}\right\Vert
	_{L_{2N}^{2,h}}\leq \widetilde{C}_{2}\delta e^{4\lambda b\theta }+\xi
	^{n}\left\Vert W_{\min ,\lambda }^{h}-W_{0}^{h}\right\Vert _{H_{2N}^{2,h}}.
	\label{3.23}
	\end{equation}
\end{thm}

\emph{Proof}. Estimate (\ref{3.21}) follows immediately from a combination
of Theorem \ref{thm:8} with Theorem 2.1 of \cite{Bakushinsky2017}. Estimate (%
\ref{3.22}) is immediately implied by (\ref{3.21}) and an obvious
modification of the proof of estimate (\ref{3.80}). More precisely, in this
modification the minimizer $V_{\min ,\lambda }^{h}$ of functional (\ref{eq:J}%
) should be replaced with $Y_{\min ,\lambda }^{h}=W_{\min ,\lambda
}^{h}+\Psi ^{h}:$ recall that $W_{\min ,\lambda }^{h}$\emph{\ }is the
minimizer of the functional $I_{h,\lambda }\left( W^{h}\right) $ on the set%
\emph{\ }$\overline{B_{0}\left( 2M\right) }.$ Finally, (\ref{3.23}) follows
from (\ref{eq:v}) and (\ref{3.22}). $\quad \square $

\begin{rem}
	Since the radius of the ball $\overline{B_{0}\left(
		2M\right) }$ is $2M$, where $M$ is an arbitrary number and since the
	starting point $W_{0}^{h}$ of iterations of the sequence (\ref{3.20}) is an
	arbitrary point of $B_{0}\left( 2M\right) ,$ then Theorem \ref{thm:9} actually claims
	the \emph{global convergence }of the sequence (\ref{3.20}) to the exact
	solution, as long as the level $\delta $ of the noise in the data tends to
	zero; see the second paragraph of the Introduction for the definition of the
	global convergence.
\end{rem} 

\section{Experimental Results\label{sec:Experimental-results}}

\subsection{Experimental setup\label{subsec:Experimental-setup}}

\begin{figure*}[t]
\begin{centering}
\subfloat[A photograph of our experimental setup$\qquad\qquad\qquad$\label{fig:Experimental-setup}]{\begin{raggedright}
\includegraphics[scale=0.35]{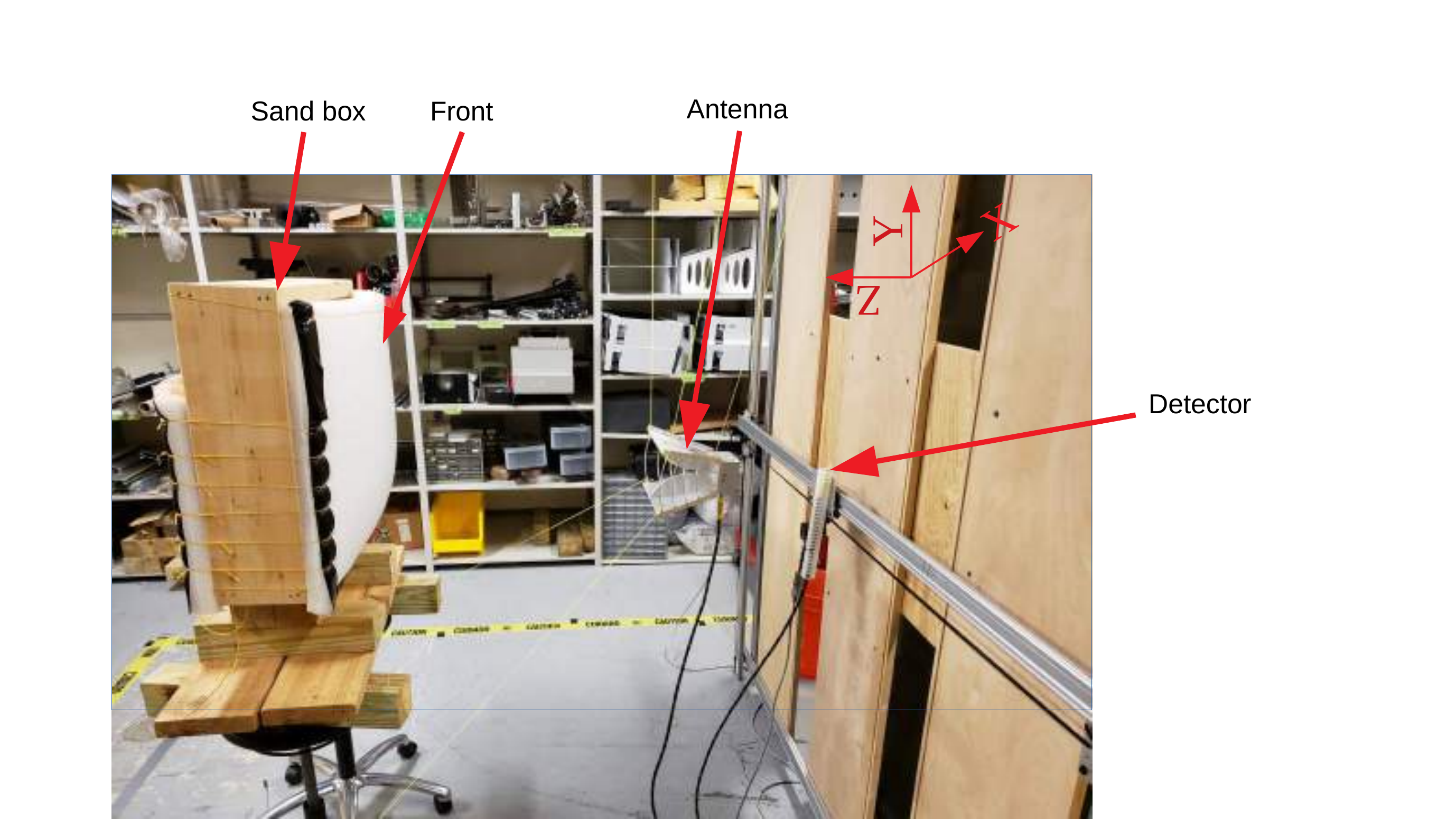}
\par\end{raggedright}
}\subfloat[A schematic diagram of sources/detectors locations in
our experimental setup$\qquad\qquad\qquad$\label{fig:Mesh-refinement}]{\begin{centering}
\includegraphics[scale=0.3]{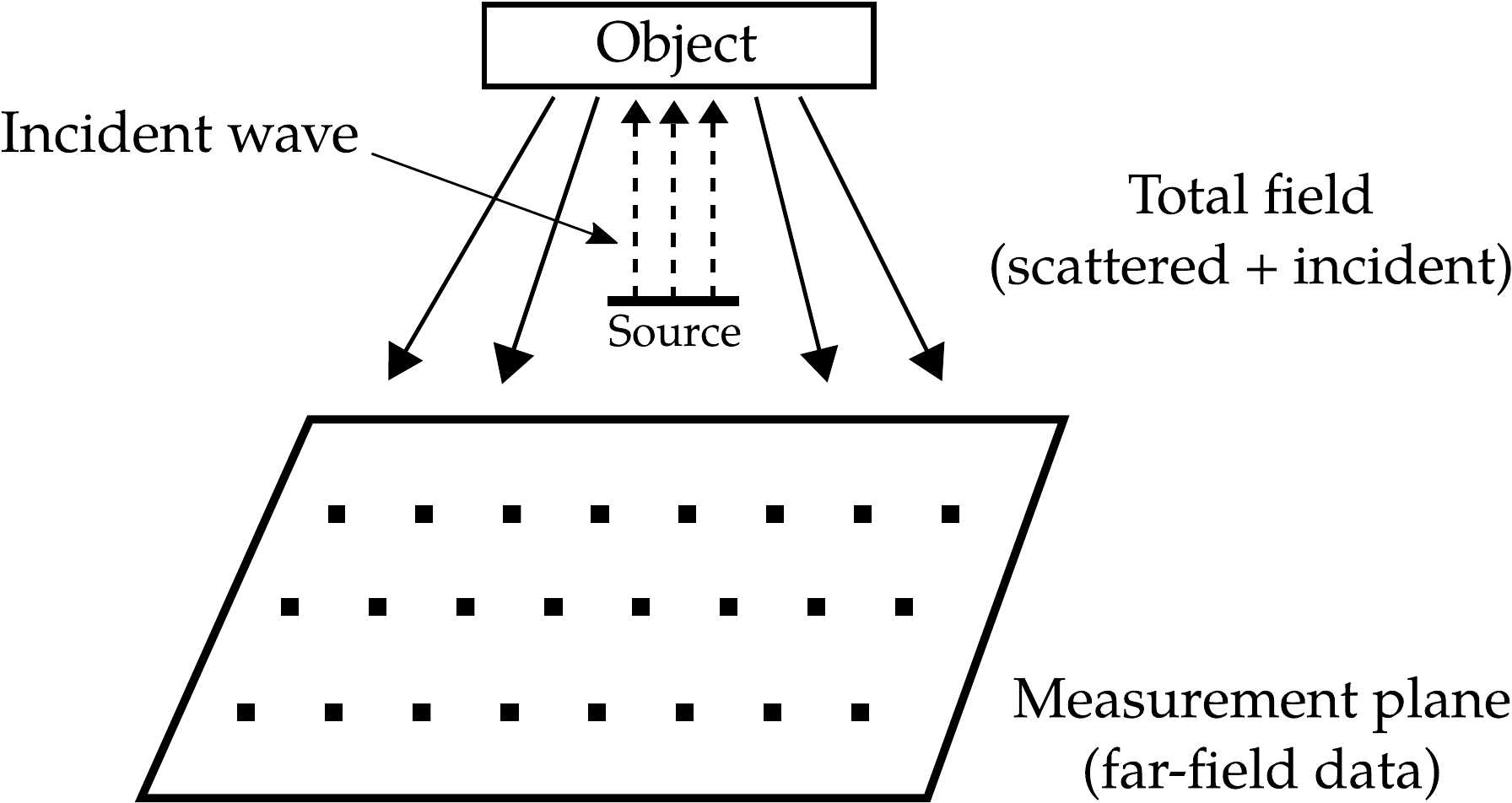}
\par\end{centering}
}
\par\end{centering}
\caption{Our experimental setup (left) and a schematic
	diagram of our measurements (right).}
\label{fig:1}
\end{figure*}

We now explain expound our experimental setup and data acquisition at the
microwave facility of University of North Carolina at Charlotte (UNCC).
Keeping in mind our target application mentioned in the first paragraph of
Introduction to imaging of explosive-like devices, We have collected
experimental for objects buried in a sandbox. More precisely, we have placed
the targets of interest inside of a wooden framed box filled with the
moisture free sand. Besides, we cover the front and back sides of the
sandbox by Styrofoam whose dielectric constant is close to 1, i.e. to the
dielectric constant of air. Hence, Styrofoam should not affect neither the
incident nor the scattered electric waves. Here, the front surface is
physically defined as the foam layer closer to the transmitter fixed at a
given position. On the other hand, the burial depths of objects do not
exceed 10 cm, which really mimics a scanning and detecting action for
shallow mine-like targets. Typically the sizes of antipersonnel land mines
and improvised explosive devices are between 5 and do 15 centimeters (cm),
see, e.g. \cite{Nguyen2018}. The transmitter is a standard horn antenna,
whose length is about 20 cm, and the detector is essentially a point probe.
To get a better insight into the description we have detailed, the reader
can take a look at Figure \ref{fig:1}.

It is worth mentioning that there are several challenges that we confront in
this configuration, which actually reflect the difficulties met in the
realistic detection of land mines. We now name some central challenges:

\begin{itemize}
\item \textbf{Distractions.} Cf. Figure \ref{fig:Experimental-setup}, we
deliberately keep many other devices and items (made of different materials)
on the desks outside the yellow caution bands. In other words, we do not use
any isolations of our device from the outside World. This is reasonable
since no isolation conditions can be created on a battlefield. Obviously,
such unwanted obstacles and furniture can affect the quality of the raw
backscatter. The presence of the Wi-Fi signal is also unavoidable in the
room where we conduct the experiments. Moreover, it is technically very hard
to place the antenna behind the measurement site. Therefore, the backscatter
wave hits the antenna first and only then comes to detectors, which is
another complicating factor.

\item \textbf{Random noise factor.} When facing real experiments, one cannot
rarely estimate the noise level as well as its frequency dependent dynamics
since they depends on hundreds of factors such as measurement process,
unknown true data, distracting signals, etc.
\end{itemize}

\subsection{Buried targets to be imaged\label{subsec:Experimental-examples}}

We present here five (5) examples of computational reconstructions of buried
objects mimicking typical metallic and non-metallic land mines. The tested
objects we use in the experiments are basic in-store items that one can
easily purchase. The burial depth of any target is not of an interest here
since all depths are just a few centimeters. The most valuable information
for the engineering part is in estimating the values of dielectric constants
of targets as well as their shapes.

Our five examples are:

\begin{itemize}
\item Example 1: An aluminum cylinder (see Figure \ref{fig:tube-1}). As
metallic mines usually caught in military services, this object can be
shaped as the NO-MZ 2B, a Vietnamese anti-personnel fragmentation mine; cf.
e.g. \cite{Banks1998}. It is known that metallic objects can be
characterized by large values of dielectric constants \cite{Kuzh}. Hence, we
suppose that the true values of dielectric constants of metallic objects are
large and are not fixed.

\item Example 2: A glass bottle filled with the clear water (see Figure \ref%
{fig:bottle-1}). This object is more complicated than the one of Example 2
due to the presence of the cap on the top of the bottle. Example 2 is a good
fit of the usual Glassmine 43 (cf. \cite{ChiefofOrdnance1945}), a
non-metallic anti-personnel land mine largely with a glass body that the
Germans used to make detection harder in the World War II era. The true
value of the dielectric constant in this case was measured to be 23.8 \cite%
{Thanh2015}.

\item Example 3: An U-shaped piece of a dry wood (see Figure \ref{fig:Walnut}%
). This example is our next attempt to deal with a non-metallic object. Note
that the shape is non convex now. In the spirit of Example 2, this
wood-based object is well-suited (in terms of the material) to the case of
Schu-mine 42, an anti-personnel blast mine that the Germans developed during
the World War II. The augmented complexity of the geometry of the object is
just our purpose of this work since we wish to see how the reconstruction
works with different front shapes. In this circumstance, the maximal
achievable value of the dielectric constant which we see in \cite%
{Table} should be 6.

\item Examples 4 and 5: Metallic letters \textquotedblleft A'' and
\textquotedblleft O'' (see Figures \ref{fig:A} and \ref{fig:O}). Shapes are
non convex. These two tests are different from the above examples because
they were \textbf{blind} tests. This means that we did not know any other
information except of the measured data and the fact that these objects were
buried close to the sand surface. Since they are metallic, the true contrast
should be large as in Example 1.
\end{itemize}

\subsection{The necessity of data propagation}

In the experimental setup, our observed and measured data are the source
dependent backscattering data of the electric field. Although our
experimental device measures the backscattering data with varied frequencies
for each location of the point source, we use only a single frequency for
each experiment when solving our CIP. Basically, these are varied far-field
data; see Figures \ref{fig:1}. However, these data are deficient, i.e. it  is unlikely
that these data can be reasonably inverted; see Figures \ref{fig:Raw1}--\ref{fig:Raw5}. In fact,
the same observation was made in previous publications of our research group
on experimental data \cite{Klibanov2019b,Nguyen2018,Thanh2015}. Hence, to
make our data feasible for inversion, we apply the well known data
propagation procedure which approximates the near field data. These
approximate data form actual inputs of our minimization process. A rigorous
justification of the data propagation procedure can be found in \cite%
{Nguyen2018}.

It is our experience that the good quality near-field data are not
always obtained well enough from any far field data after the propagation.
This requires a substantial workload in choosing proper data among a large
amount of frequency dependent data sets. In other words, we have no choice
but to select an acceptable frequency for each particular target we work
with. So, we select its own frequency for each considered target. Then we
use this frequency for all positions of the source we work with; see section
\ref{sec:4.4} about particular choices of frequencies. A particular choice of an
admissible and acceptable set of data has been illustrated in \cite%
{Nguyen2018}, and below this strategy will be confirmed again.

\subsection{Data propagation revisited\label%
{subsec:Data-propagation-revisited}}

We know in advance that the half space $\left\{ z<-b\right\} \subset \mathbb{%
R}^{3}$ is homogeneous, i.e. $c\left( \mathbf{x}\right) =1$ in this half
space. Therefore, the function $u_{s}$ is a backscatter wave in $\left\{
z<-b\right\} $ and it satisfies the following conditions: 
\begin{equation}
\begin{cases}
\Delta u_{s}+k^{2}u_{s}=0 & \text{for }\mathbf{x}\in \left\{ z<-b\right\} ,
\\ 
\partial _{r}u_{s}-\text{i}ku_{s}=\mathcal{O}\left( r^{-1}\right)  & \text{%
for }r=\left\vert \mathbf{x}-\mathbf{x}_{\alpha }\right\vert ,\text{i}=\sqrt{%
-1}.%
\end{cases}
\label{eq:scat}
\end{equation}

As was mentioned in section 3.1, we actually measure the far field data,
i.e. the function $u_{s}\left( x,y,-D,\mathbf{x}_{\alpha }\right) $, where
the number $D>b$. Having the function $u_{s}\left( x,y,-D,\mathbf{x}_{\alpha
}\right) ,$ we want approximate the function $u_{s}\left( x,y,-b,\mathbf{x}%
_{\alpha }\right) ,$ i.e. we want to approximate the wave field in the near
field zone. The data propagation procedure does exactly this. Denote 
\begin{equation}
u_{s}\left( x,y,-b,\mathbf{x}_{\alpha }\right) =\mathbf{U}\left( x,y,\mathbf{%
x}_{\alpha }\right) \;\text{and}\;u_{s}\left( x,y,-D,\mathbf{x}_{\alpha
}\right) =\mathbf{V}\left( x,y,\mathbf{x}_{\alpha }\right) .  \label{eq:data}
\end{equation}

In this work, we rely on the data propagation procedure to unveil this
difficulty as it has been successfully exploited in \cite{Nguyen2018}.
First, we apply the Fourier transform of the scattered field with respect to 
$x,y$, assuming that the corresponding integral converges: 
\begin{equation}\label{399}
\hat{u}_{s}\left( \rho _{1},\rho _{2},z,\mathbf{x}_{\alpha }\right) =\frac{1%
}{2\pi }\int_{\mathbb{R}^{2}}u_{s}\left( x,y,z,\mathbf{x}_{\alpha }\right)
e^{-\text{i}\left( x\rho _{1}+y\rho _{2}\right) }dxdy\quad \text{for }\rho
_{1},\rho _{2}\in \mathbb{R}.
\end{equation}%
Next, we apply this Fourier transform to the PDE in \eqref{eq:scat} and
arrive at a second order ODE with respect to $z$: 
\begin{equation}
\partial _{zz}^{2}\hat{u}_{s}+\left( k^{2}-\rho _{1}^{2}-\rho
_{2}^{2}\right) \hat{u}_{s}=0\quad \text{for }z<-b.  \label{400}
\end{equation}%
By \eqref{eq:data}, we also have 
\begin{equation*}
\hat{u}_{s}\left( \rho _{1},\rho _{2},-b,\mathbf{x}_{\alpha }\right) =\hat{%
\mathbf{U}}\left( \rho _{1},\rho _{2},\mathbf{x}_{\alpha }\right) \;\text{and%
}\;\hat{u}_{s}\left( \rho _{1},\rho _{2},-D,\mathbf{x}_{\alpha }\right) =%
\hat{\mathbf{V}}\left( \rho _{1},\rho _{2},\mathbf{x}_{\alpha }\right) .
\end{equation*}

It follows from (\ref{400}) that 
\begin{equation*}
\hat{u}_{s}\left( \rho _{1},\rho _{2},z,\mathbf{x}_{\alpha }\right) =\left\{ 
\begin{array}{c}
\hat{\mathbf{U}}\left( \rho _{1},\rho _{2},\mathbf{x}_{\alpha }\right) e^{%
\sqrt{\rho _{1}^{2}+\rho _{2}^{2}-k^{2}}\left( z+b\right) }\text{ if }\rho
_{1}^{2}+\rho _{2}^{2}>k^{2}, \\ 
C_{1}e^{-\text{i}\sqrt{k^{2}-\rho _{1}^{2}-\rho _{2}^{2}}\left( z+b\right)
}+C_{2}e^{\text{i}\sqrt{k^{2}-\rho _{1}^{2}-\rho _{2}^{2}}\left( z+b\right) }%
\text{ if }\rho _{1}^{2}+\rho _{2}^{2}<k^{2},%
\end{array}%
\right. 
\end{equation*}%
where $z<-b.$ It is not immediately clear which of two terms in the second
line of the last formula should be taken. However, it was proven in Theorem
4.1 of \cite{Nguyen2018} that only the first term which should be taken and
one should set $C_{2}:=0.$ Thus, for $z<-b$ 
\begin{equation}
\hat{u}_{s}\left( \rho _{1},\rho _{2},z,\mathbf{x}_{\alpha }\right) =%
\begin{cases}
\hat{\mathbf{U}}\left( \rho _{1},\rho _{2},\mathbf{x}_{\alpha }\right) e^{%
\sqrt{\rho _{1}^{2}+\rho _{2}^{2}-k^{2}}\left( z+b\right) } & \text{if }\rho
_{1}^{2}+\rho _{2}^{2}>k^{2}, \\ 
\hat{\mathbf{U}}\left( \rho _{1},\rho _{2},\mathbf{x}_{\alpha }\right) e^{-%
\text{i}\sqrt{k^{2}-\rho _{1}^{2}-\rho _{2}^{2}}\left( z+b\right) } & \text{%
otherwise}.%
\end{cases}
\label{eq:prop}
\end{equation}

Observe that if the Fourier frequency satisfies $\rho _{1}^{2}+\rho
_{2}^{2}>k^{2}$, the function $\hat{u}_{s}\left( \rho _{1},\rho _{2},z,%
\mathbf{x}_{\alpha }\right) $ decays exponentially with respect to $%
z\rightarrow -\infty $. Therefore, if the measurement surface is far away
from the domain of interest, i.e. $D$ is large, then we can neglect the term
in the first line of (\ref{eq:prop}). In other words, we can neglect high
spatial frequencies in (\ref{eq:prop}). Thus, we take $z=-D$ in %
\eqref{eq:prop} to get 
\begin{equation*}
\hat{\mathbf{U}}\left( \rho _{1},\rho _{2},\mathbf{x}_{\alpha }\right) =\hat{%
\mathbf{V}}\left( \rho _{1},\rho _{2},\mathbf{x}_{\alpha }\right) e^{\text{i}%
\sqrt{k^{2}-\rho _{1}^{2}-\rho _{2}^{2}}\left( -D+b\right) }\quad \text{for }%
\rho _{1}^{2}+\rho _{2}^{2}<k^{2}.
\end{equation*}
Using the inverse Fourier transform, we obtain 
\begin{align}
& \mathbf{U}\left( x,y,\mathbf{x}_{\alpha }\right) =u_{s}\left( x,y,-b,%
\mathbf{x}_{\alpha }\right)   \label{eq:propdata} \\
& =\frac{1}{2\pi }\dint\limits_{\rho _{1}^{2}+\rho _{2}^{2}<k^{2}}\hat{%
\mathbf{V}}\left( \rho _{1},\rho _{2},\mathbf{x}_{\alpha }\right) e^{\text{i}%
\sqrt{k^{2}-\rho _{1}^{2}-\rho _{2}^{2}}D}e^{\text{i}\left( x\rho _{1}+y\rho
_{2}\right) }d\rho _{1}d\rho _{2}  \notag \\
& =\frac{1}{\left( 2\pi \right) ^{2}}\dint\limits_{\rho _{1}^{2}+\rho
_{2}^{2}<k^{2}}\left[ \dint\limits_{\mathbb{R}^{2}}u_{s}\left( \tilde{x},%
\tilde{y},-D,\mathbf{x}_{\alpha }\right) e^{-\text{i}\left( \tilde{x}\rho
_{1}+\tilde{y}\rho _{2}\right) }d\tilde{x}d\tilde{y}\right] e^{\text{i}\sqrt{%
k^{2}-\rho _{1}^{2}-\rho _{2}^{2}}\left( -D+b\right) }e^{\text{i}\left(
x\rho _{1}+y\rho _{2}\right) }d\rho _{1}d\rho _{2}.  \notag
\end{align}
The last formula of \eqref{eq:propdata} is the actual data propagation
procedure we will use in this work.

\subsection{Computational setup}\label{sec:4.4}

We introduce dimensionless variables as $\mathbf{x}^{\prime }=\mathbf{x}/(10 \text{ cm})
$ and keep the same notations as before, for brevity. This means that the
dimensions we use in computations are 10 times less than the real ones in
centimeters. We illustrate the choice of the coordinate system on Figures
2c, 3c, 4c, 5c and 6c: the $x-$ and $y-$axis are horizontal and vertical
sides respectively and $z-$axis is orthogonal to the measurement plane.

The far-field data are measured on a rectangular surface of dimensions 100
cm $\times $ 100 cm, i.e. $10\times 10$ in dimensionless regime. Cf. Figure %
\ref{fig:Mesh-refinement} as to our mesh grid of the measurement plane, each
step is 2 cm (0.2) over 100 cm (10) length row. The total number of steps in
a row is 50, and the total number of steps in a column is also 50. The
distance between the measurement plane and the sandbox with the foam layer,
whose thickness is 5 cm, is about 110.5 cm (11.05). The length in the $z$
direction of the sandbox without the foam is approximately 44 cm, but due to
the bending foam layer, we reduce 10\% of this length. Henceforth, our
choice of the domain $\Omega $ should be 
\begin{equation*}
\Omega =\left\{ \mathbf{x}\in \mathbb{R}^{3}:\left\vert x\right\vert
,\left\vert y\right\vert <5,\left\vert z\right\vert <2\right\} ,
\end{equation*}%
which implies that $R=5$ and $b=2$. The near-field or propagated measurement
site is then assigned as 
\begin{equation*}
\Gamma :=\left\{ \mathbf{x}\in \mathbb{R}^{3}:\left\vert x\right\vert
,\left\vert y\right\vert <5,z=-2\right\} .
\end{equation*}%
Also, we take $D=14$ for the far-field measurement site as we estimate the
distance between this site and the zero point. Meanwhile, for all objects,
for the line of sources $L_{\text{src}}$ defined in section \ref%
{sec:Statement-of-the} we have $d=9$, $a_{1}=0.1$ and $a_{2}=0.6$. Besides, we take $\theta = 4$.

It remains to obtain the wavenumber $k$ corresponding to the dimensionless spatial variables we are working with. It is well-known that
the relation between the wavelength ($\tilde{\lambda}$) and the wavenumber
is expressed by $k=2\pi /\tilde{\lambda}$. Basically, the wavelength can be
computed via the formulation $\tilde{\lambda}=\tilde{v}/\tilde{f}$, where $%
\tilde{v}=299792458$ (m/s) is the speed of light in vacuum and $\tilde{f}$
is the frequency in Hertz (Hz or $\text{s}^{-1}$). Hence, in the computational setting we compute (in $\text{cm}^{-1}$)
\begin{equation*}
k=\frac{2\pi }{2997924580}\tilde{f}.
\end{equation*}%
The choice of $k$ relies on the performance of the
data after preprocessing. More precisely, our criterion is heuristically based upon the best visualization of the propagated data that we obtain using the data propagation. For each example below we then use its own frequency, which we specify in Table \ref{table:1}.
Note that for each location of the detector we measure the backscatter data
for 300 frequency points uniformly distributed between 1 GHz and 10 GHz.

\begin{table}
\begin{center}
\begin{tabular}{|c|c|c|c|c|c|}
\hline Example & 1 & 2 & 3 & 4 & 5 \tabularnewline
\hline 
$k$ & 8.51 & 6.62 & 11.43 & 9.55 & 8.79 \tabularnewline
\hline 
Frequency (GHz) & 4.06 & 3.16 & 5.45 & 4.55 & 4.19\tabularnewline
\hline 

\end{tabular}
\end{center}
\caption{Wave numbers and frequencies for Examples 1--5.\label{table:1}}
\end{table}

Now, we summarize the crucial steps of the data preprocessing to obtain fine
data for our inversion method from the raw ones.

\begin{itemize}
\item \textbf{Step 1.} For every frequency and for every location of the
source, we subtract the reference data from the far-field measured data. A
similar procedure was implemented in \cite{Klibanov2019b,Nguyen2018}. The
reference data are the background ones measured when the sandbox is without
a target. This subtraction helps to extract the pure signals from buried
objects from the whole signal. Therefore, we reduce the noise this way.

\item \textbf{Step 2.} We apply the data propagation procedure as in
subsection \ref{subsec:Data-propagation-revisited} to obtain the near-field
data. This procedure provides a significantly better estimation for $x,y$
coordinates of buried objects, and reduces the size of the computational
domain in the $z$--direction.

\item \textbf{Step 3.} We truncate the so obtained near-field data to get
rid of random oscillations. The oscillations appear randomly during the data
propagation and may cause unnecessary issues during our inversion procedure.
This data truncation was developed in \cite{Nguyen2018} and now we
improve it using the following two steps, given a function $g\left( x,y,\alpha
\right)$ to be truncated:

\begin{itemize}
\item For each point source, we replace the function $g\left( x,y,\alpha
\right) $ with a function $\tilde{g}\left( x,y,\alpha \right) $ defined as:
\begin{equation*}
\tilde{g}\left( x,y,\alpha \right) =%
\begin{cases}
g\left( x,y,\alpha \right)  & \text{if }\left\vert g\left( x,y,\alpha
\right) \right\vert \geq \kappa _{1}\max_{\left\vert x\right\vert
,\left\vert y\right\vert \leq R}\left\vert g\left( x,y,\alpha \right)
\right\vert , \\ 
0 & \text{otherwise.}%
\end{cases}%
\end{equation*}

Here, we call $\kappa _{1}>0$ the truncation number. Even though this number
should be dependent of the source position $\alpha $ and should be different
from every single choice of the frequency point, we apply the same
truncation number to all the examples below. By the trial and error procedure, we have chosen $\kappa _{1}=0.4$, which means that we only preserve those
propagated near-field data whose values are  least 40 percents of the global
maximum value.

\item The next step would be smoothing the function $\tilde{g}$ using the
Gaussian filter. However, we notice that when doing so, the maximum value of 
$\tilde{g}$ will be smaller than that of $g$. In order to preserve this
important \textquotedblleft peak\textquotedblright\ of $g$ after truncation,
we add back some percents of $\tilde{g}$ in the following manner: 
\begin{equation}
\tilde{g}_{\text{new}}\left( x,y,\alpha \right) =\kappa _{2}\tilde{g}_{\text{%
old}}\left( x,y,\alpha \right) .  \label{eq:28}
\end{equation}%
Here, we call $\kappa _{2}>0$ the retrieval number. This number is computed
by $\kappa _{2}=\max \left( \left\vert \tilde{g}\right\vert \right) /\tilde{m%
}$, where $\tilde{m}$ is the maximal absolute value of the smoothed $\tilde{g%
}_{\text{old}}$.
\end{itemize}
\end{itemize}

\subsubsection*{Fully discrete setting}

We now present our numerical approach of the approximation of the right
hand side of formula \eqref{eq:propdata} in order to use it for our
experimental data. First, we adapt the conventional Riemannian sum
approximation to compute the Fourier transform of the function $\mathbf{V}$.
Using (\ref{399}) and the samples $\left\{ u_{s}\left( \tilde{x}_{i},\tilde{y}_{j},-D,%
\mathbf{x}_{\alpha }\right) \right\} _{i,j=0}^{\tilde{N}-1}$ over a 2D
finite domain, where we are experimentally measuring the far-field data, we
find that 
\begin{equation*}
\hat{\mathbf{V}}\left( \rho _{1},\rho _{2},\mathbf{x}_{\alpha }\right)
\approx \omega ^{2}\sum_{i,j=0}^{\tilde{N}-1}u_{s}\left( \tilde{x}_{i},%
\tilde{y}_{j},-D,\mathbf{x}_{\alpha }\right) \exp \left( -\text{i}\left( 
\tilde{x}_{i}\rho _{1}+\tilde{y}_{j}\rho _{2}\right) \right) .
\end{equation*}%
Here, a uniform sampling rate, i.e. $\tilde{x}_{i}=i\Delta \tilde{x}_{i},%
\tilde{y}_{j}=j\Delta \tilde{y}_{j}$, is used with $\Delta \tilde{x}%
_{i}=\Delta \tilde{y}_{j}=\omega $ for a number $\omega \in \left( 0,1\right) $.

Next, we define the following truncated Fourier domain in 2D: 
\begin{equation*}
\Theta _{k}:=\left\{ \left( \rho _{1},\rho _{2}\right) \in \mathbb{R}%
^{2}:\rho _{1}^{2}+\rho _{2}^{2}<k^{2}\right\} .
\end{equation*}%
We sample this truncated Fourier domain at uniformly discrete points $\rho
_{1m_{1}}=m_{1}\omega _{\rho },\rho _{2m_{2}}=m_{2}\omega _{\rho }$ for a number 
$\omega _{\rho }\in \left( 0,1\right) $ and $0\leq m_{1},m_{2}\leq \tilde{M}%
-1$, provided that these points are in the set $\Theta _{k}$. Thus, we
conclude that%
\begin{align}
&\mathbf{U}\left( x_{p},y_{q},\alpha _{l}\right)\label{29} \\ &   
\approx \frac{1}{\left( 2\pi \right) ^{2}}\omega _{\rho
}^{2}\sum_{m_{1},m_{2}=0}^{\tilde{M}-1}\hat{\mathbf{V}}\left( \rho
_{1m_{1}},\rho _{2m_{2}},\alpha _{l}\right) \exp \left( \text{i}\sqrt{%
k^{2}-\rho _{1m_{1}}^{2}-\rho _{2m_{2}}^{2}}\left( -D+b\right) \right) \exp
\left( \text{i}\left( x_{p}\rho _{1m_{1}}+y_{q}\rho _{2m_{2}}\right) \right)
.\notag 
\end{align}

In our experimental data, we have $N_{p}=N_{q}=51,$ where $N_{p}$ and $N_{q}$
are the number of discrete points in $x$ and $y$ directions respectively.
Therefore, we take $\tilde{N}=\tilde{M}=51$, which gives $\omega =\omega
_{\rho }=1/50$. Thus, (\ref{29}) gives us the approximate Dirichlet boundary
condition $V\left( \mathbf{x}^{h}\right) =\psi _{0}^{h}\left( \mathbf{x}%
^{h}\right) $ at $\left\{ z=-b\right\} $ in (\ref{3.3}). Since we also need
the function $V_{z}\left( \mathbf{x}^{h}\right) =\psi _{1}^{h}\left( \mathbf{%
x}^{h}\right) $ at $\left\{ z=-b\right\} $ in (\ref{3.3}), then to obtain
it, we formally replace in (\ref{29}) $b$ with $z$, differentiate the right
hand side of the obtained equality with respect to $z$, then set again $z:=-b
$ and calculate the resulting sum. The result is $\psi _{1}^{h}\left( 
\mathbf{x}^{h}\right) $ at $\left\{ z=-b\right\} $ in (\ref{3.3}).

Hence, the Cauchy boundary data in (\ref{3.3}) are in the fully discrete
form now. Then we write the functional $J_{h,\lambda }\left( V^{h}\right) $
defined in (\ref{eq:J}) in the fully discrete form, similarly to the
semi-discrete form in (\ref{eq:J}). In this fully discrete setting we take
into account the grid points in $x,y,z$ directions, $\left\{ \left(
x_{p},y_{q},z_{s}\right) \right\} _{p,q,s=0}^{Z_{h}}$. For brevity, we do
not bring in here this fully discrete form of $J_{h,\lambda }\left( V^{h}\right) .
$ 

After the global minimum $V_{p,q,s}$ of the functional $J_{h,\lambda }\left(
V^{h}\right) $ (in its discrete form) is obtained, we compute an
approximation of the unknown dielectric constant $c_{p,q,s}$ using the
following formula: 
\begin{equation*}
c_{p,q,s}=\text{mean}_{\alpha _{l}}\left\vert -\frac{\Delta
^{h}v_{p,q,s,\alpha _{l}}+\left( \nabla ^{h}v_{p,q,s,\alpha _{l}}\right)
^{2}+2\nabla ^{h}v_{p,q,s,\alpha _{l}}\cdot \tilde{\mathbf{x}}_{p,q,s,\alpha
_{l}}}{k^{2}}\right\vert +1,
\end{equation*}%
which is resulted from \eqref{eq:v}; see Remark \ref{rem:X}. Here, $v_{p,q,s,\alpha
_{l}}=v\left( x_{p},y_{q},z_{s},\alpha _{j}\right) .$ Recall that $\tilde{%
\mathbf{x}}_{p,q,s,\alpha _{l}}$ denote vectors $\tilde{\mathbf{x}}_{\alpha }
$ at $\left( x_{p},y_{q},z_{s}\right) $ for every $\alpha _{l};$ see subsection
\ref{subsec:3.1}. Since the number of point sources is small, we apply the
Gauss--Legendre quadrature method to compute the measured data in the
Fourier mode. This was mentioned already in our previous work with simulated
data; cf. \cite{Khoa2019}.

\begin{figure}[H]
\begin{centering}
		\subfloat[Raw data at $\alpha=0.4$\label{fig:Raw1}]{\begin{centering}
				\includegraphics[scale=0.5]{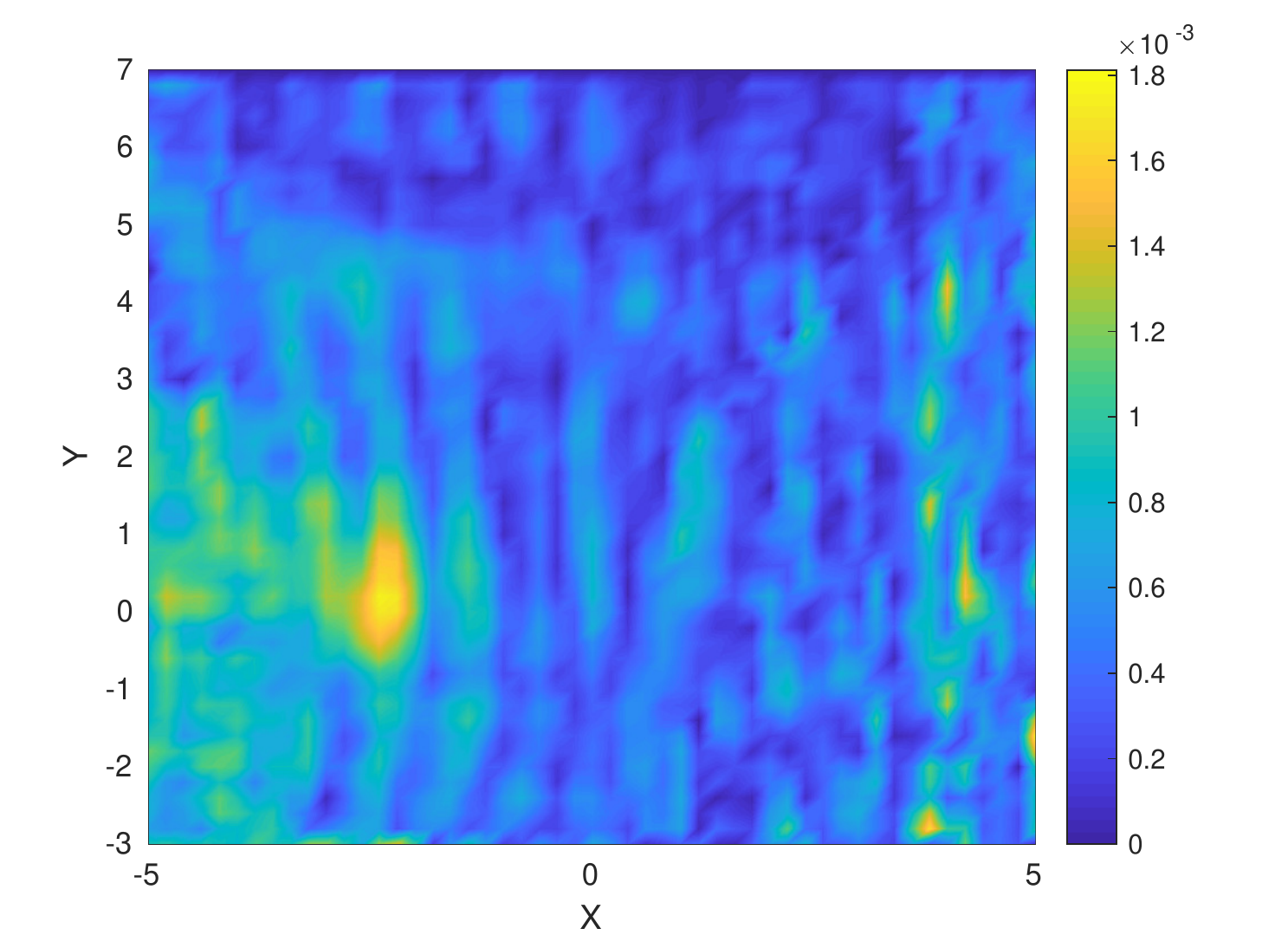}
				\par\end{centering}
		}\subfloat[Propagated data at $\alpha=0.4$\label{fig:Prop1}]{\begin{centering}
				\includegraphics[scale=0.5]{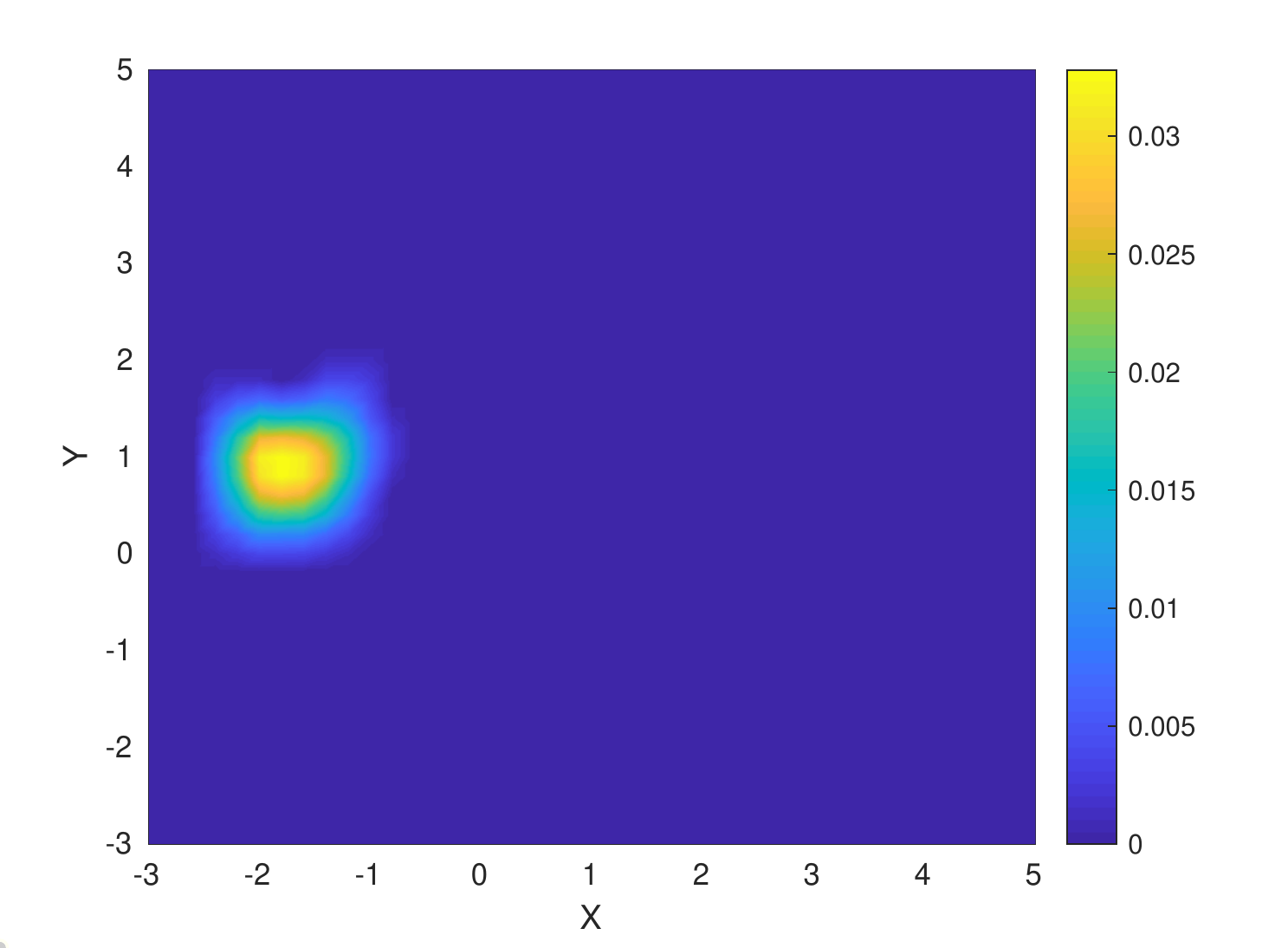}
				\par\end{centering}
		}
		\par\end{centering}
\begin{centering}
		\subfloat[Aluminum cylinder\label{fig:tube-1}]{\begin{centering}
				\includegraphics[scale=0.27]{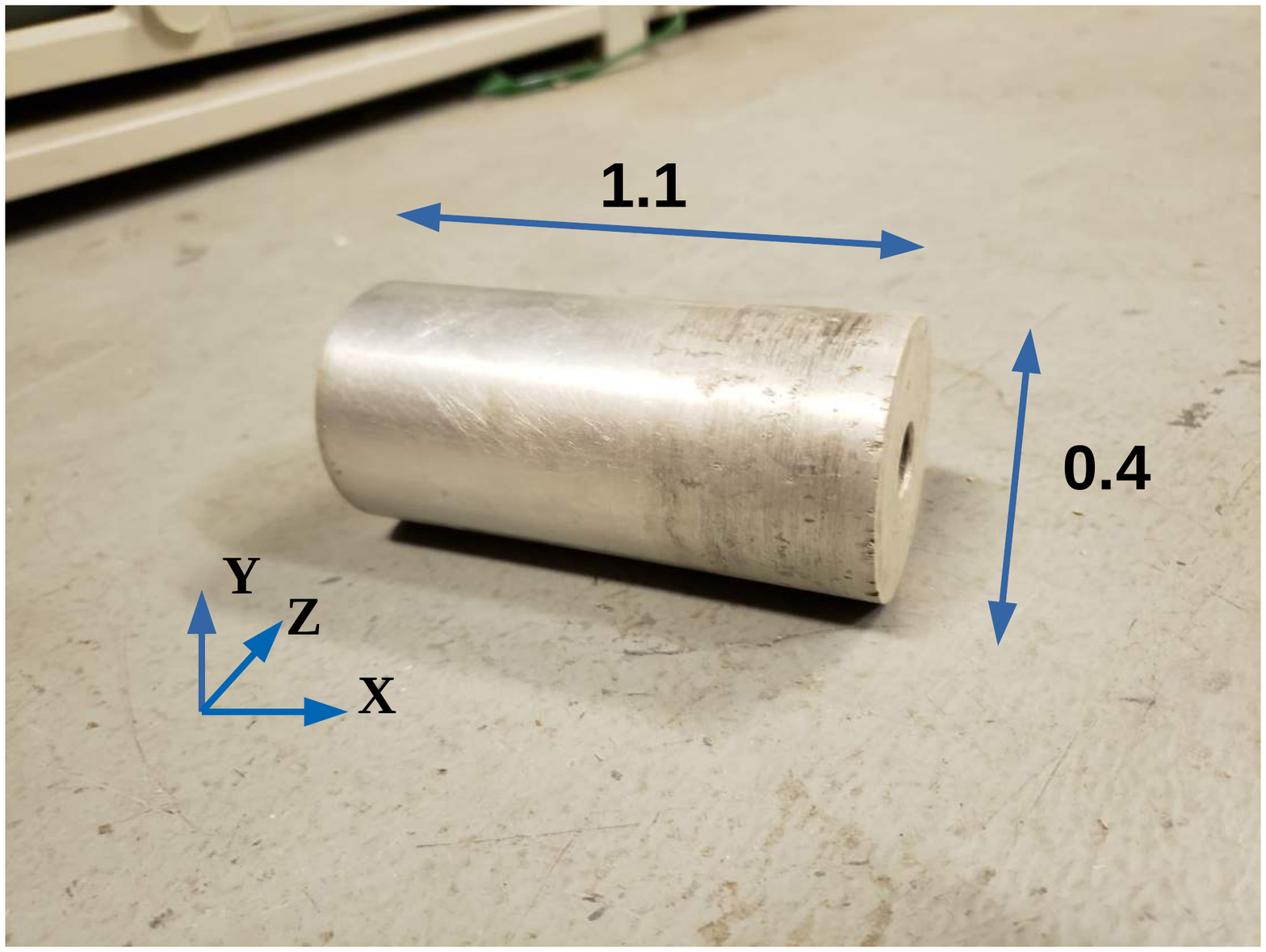}
				\par\end{centering}
		}\subfloat[Computed inclusion\label{fig:Comp1}]{\begin{centering}
				\includegraphics[scale=0.6]{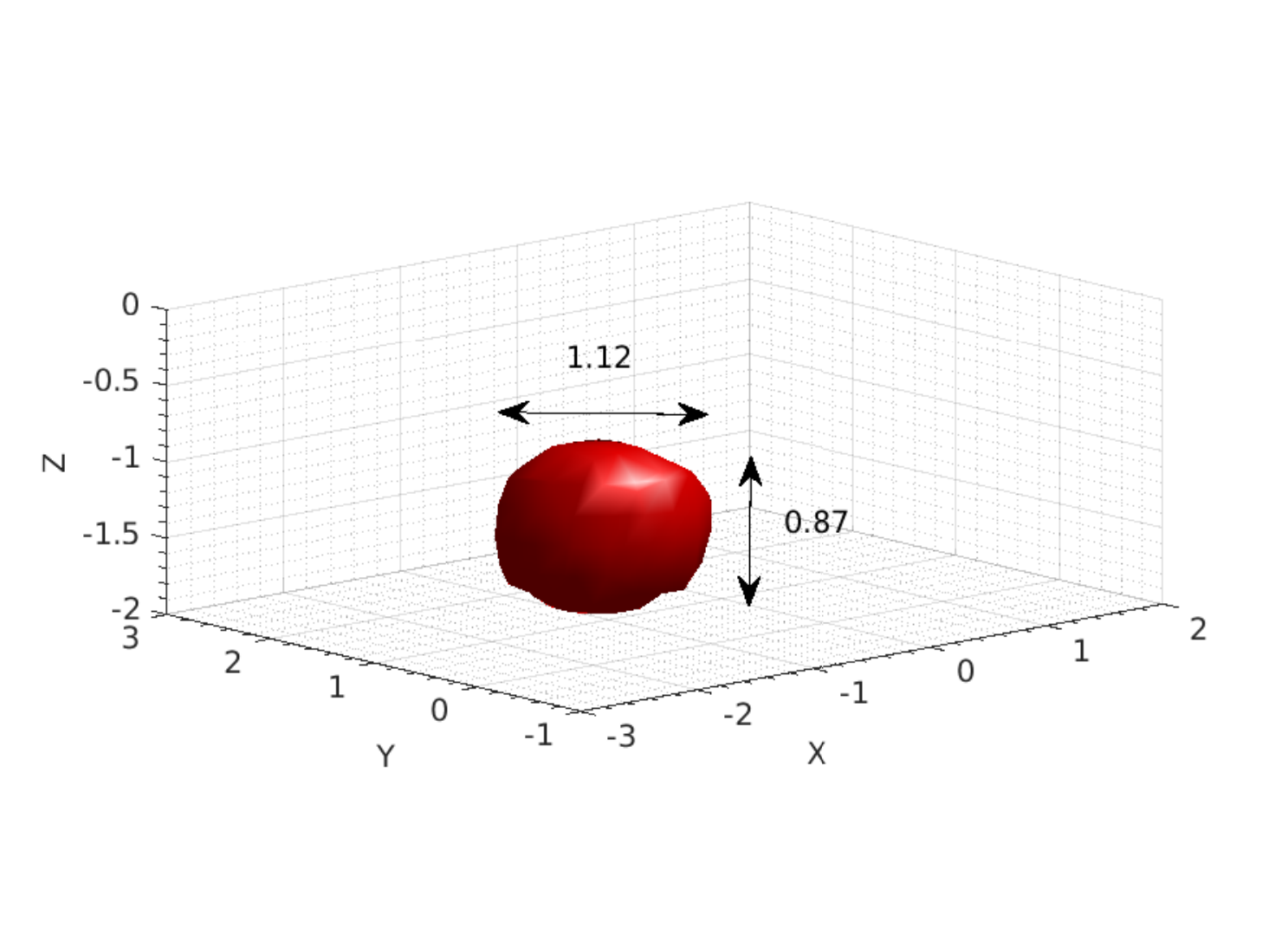}
				\par\end{centering}
		}
		\par\end{centering}
\caption{Reconstruction results of Test 1 (aluminum tube). (a) Illustration
of the absolute value of the raw far-field data; (b) Illustration of the absolute value of the near-field data after the data 
propagation procedure; (c) Photo of the experimental object; (d) The computed image of (c). All images are in the  dimensionless variables.}
\end{figure}

Since this work focuses on the detection and identification of antipersonnel
land mines and IEDs, we know that the sizes of these targets are between 5
and 15 cm, cf. e.g. \cite{Nguyen2018}. Therefore, we search for targets in a
sub-domain of $\Omega $ with only 20 cm in depth in the $z-$direction.
Denote this sub-domain by $\Omega _{1}=\left\{ -b\leq z\leq -b+2\right\} $.
We consider the following vector $V_{0}^{h}=V_{0}\left(
x_{p},y_{q},z_{s}\right) $ as the starting point of iterations in the
minimization of the functional $J_{\lambda,h }\left( V^{h}\right) $: 
\begin{equation}
V_{0}^{h}=%
\begin{pmatrix}
v_{00}^{h} & v_{01}^{h} & \cdots  & v_{0\left( N-1\right) }^{h}%
\end{pmatrix}%
^{T},\quad v_{0n}^{h}=\left( \psi _{0n}^{h}+\psi _{1n}^{h}\left( z+b\right)
\right) \chi \left( z\right) .  \label{30}
\end{equation}%

\begin{figure}[H]
	\begin{centering}
		\subfloat[Raw data at $\alpha=0.5$\label{fig:Raw2}]{\begin{centering}
				\includegraphics[scale=0.5]{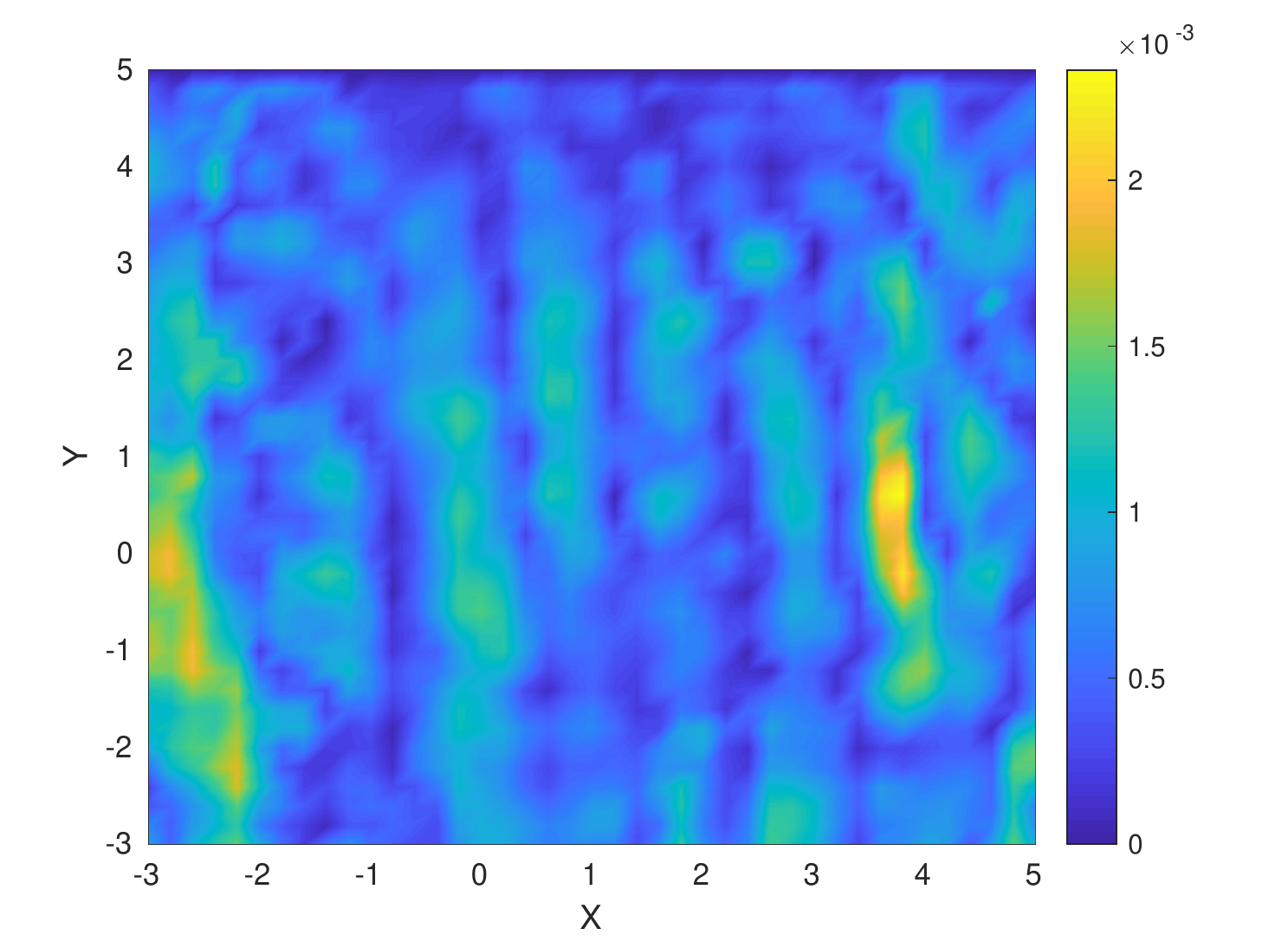}
				\par\end{centering}
		}\subfloat[Propagated data at $\alpha=0.5$\label{fig:Prop2}]{\begin{centering}
				\includegraphics[scale=0.5]{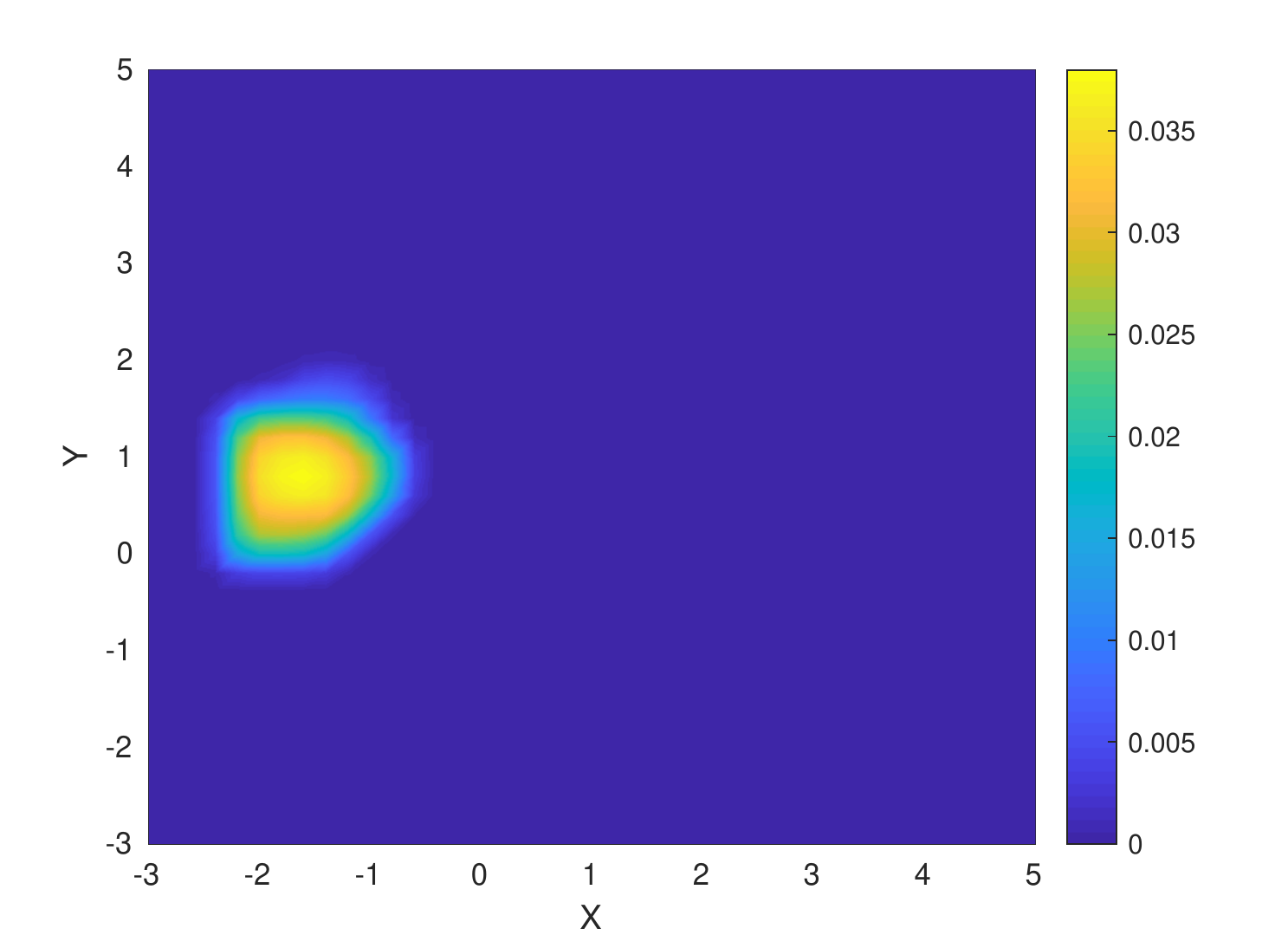}
				\par\end{centering}
		}
		\par\end{centering}
	\begin{centering}
		\subfloat[Glass bottle\label{fig:bottle-1}]{\noindent \begin{centering}
				\includegraphics[scale=0.27]{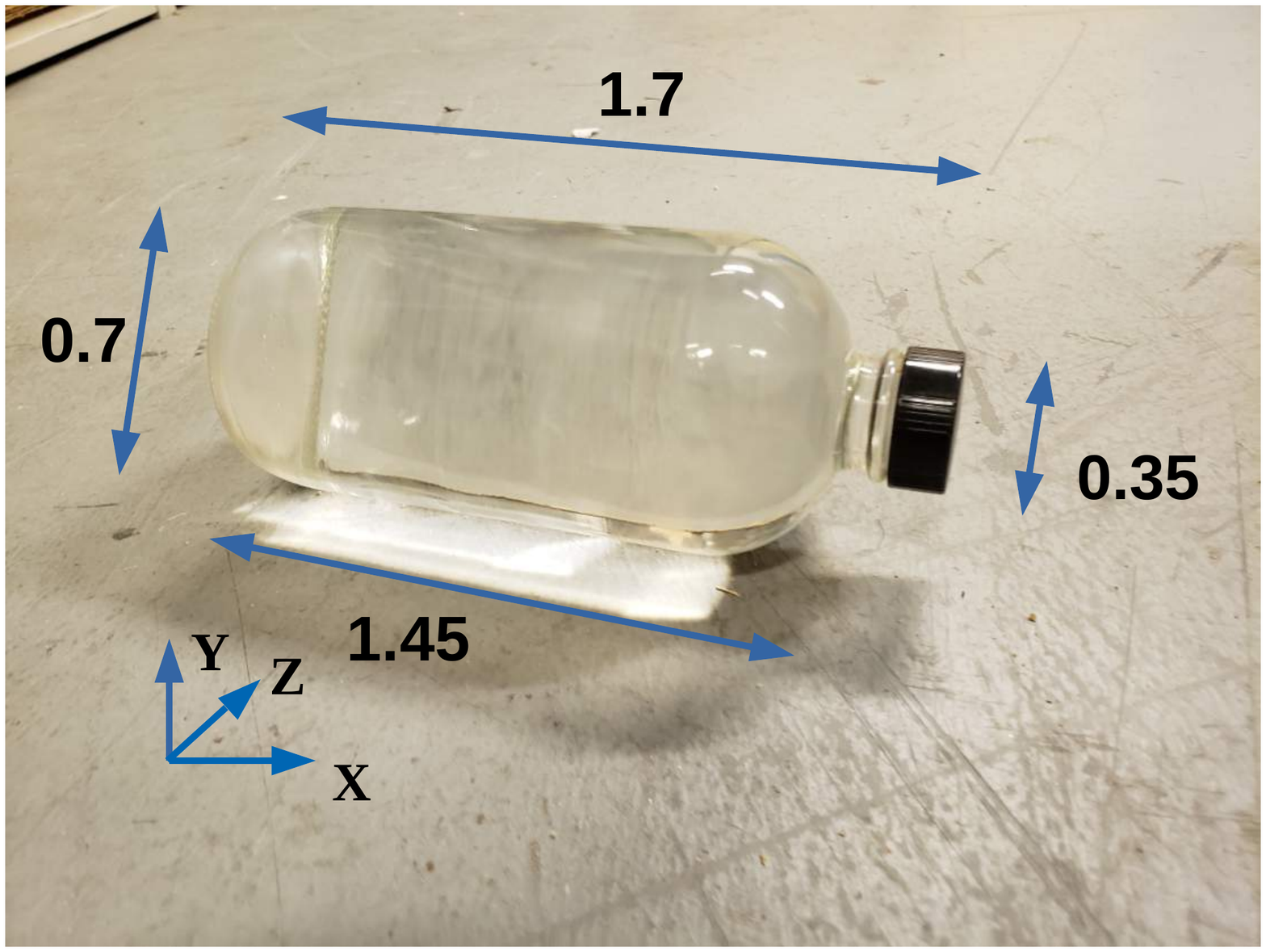}
				\par\end{centering}
		}\subfloat[Computed inclusion. Note that the cap of the bottle is clearly seen.\label{fig:Comp2}]{\begin{centering}
				\includegraphics[scale=0.6]{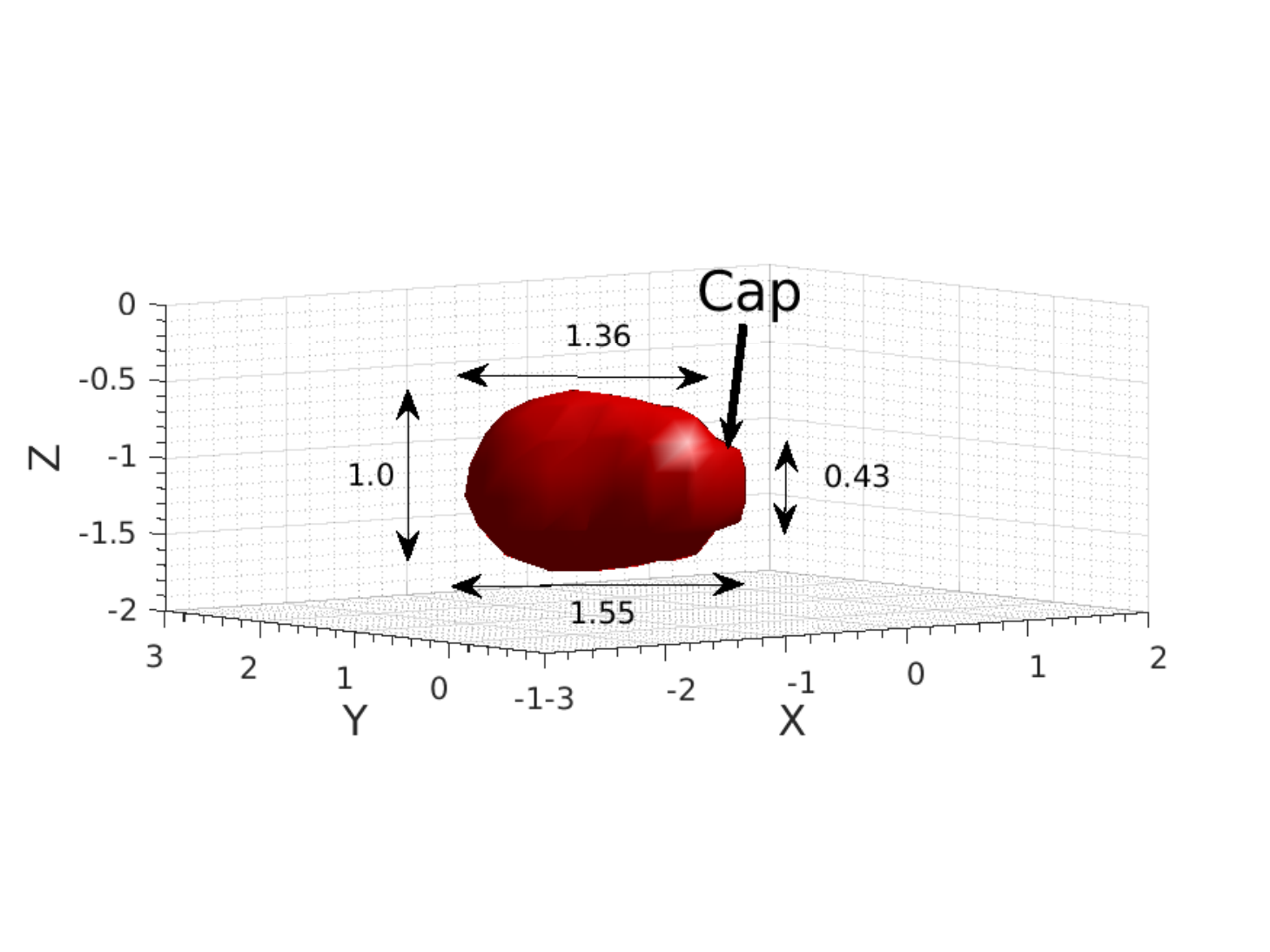}
				\par\end{centering}
		}
		\par\end{centering}
	\caption{Reconstruction results of Test 2 (a glass bottle filled with clear water). (a) Illustration
		of the absolute value of the raw far-field data; (b) Illustration of the absolute value of the near-field data after the data 
		propagation procedure; (c) Photo of the experimental object; (d) The computed image of (c). All images are in the dimensionless
		variables. An interesting point here is that we can even see the cap of the bottle in (d), which is challenging to image.}
\end{figure}

Recall that $\psi _{0n}^{h}$ and $\psi _{1n}^{h}$ are the Fourier
coefficients of the propagated data in (\ref{3.3}). Here, $\chi :\left[ -b,b%
\right] \rightarrow \mathbb{R}$ is the smooth function given by 
\begin{equation*}
\chi \left( z\right) =%
\begin{cases}
\exp \left( \frac{2\left( z+b\right) ^{2}}{\left( z+b\right) ^{2}-b^{2}}%
\right)  & \text{if }z<0, \\ 
0 & \text{otherwise}.%
\end{cases}%
\end{equation*}%
This function attains the maximum value 1 at $z=-b$ where the propagated
data are given. Then, it is easy to see that $v_{0n}^{h}\mid _{z=-b}=\psi
_{0n}^{h}$, $\partial _{z}v_{0n}^{h}\mid _{z=-b}=\psi _{1n}^{h}$. On the
other hand, $\chi $ tends to 0 as $z\rightarrow 0^{+},$ which, in
particular, means that $v_{0n}^{h}\mid _{z=b}=\partial _{z}v_{0n}^{h}\mid
_{z=b}=0.$ Thus, this starting point (\ref{30}) of iterations satisfies the
boundary conditions (\ref{3.3}). 

\begin{figure}[H]
	\begin{centering}
		\subfloat[Raw data at $\alpha=0.5$\label{fig:Raw3}]{\begin{centering}
				\includegraphics[scale=0.5]{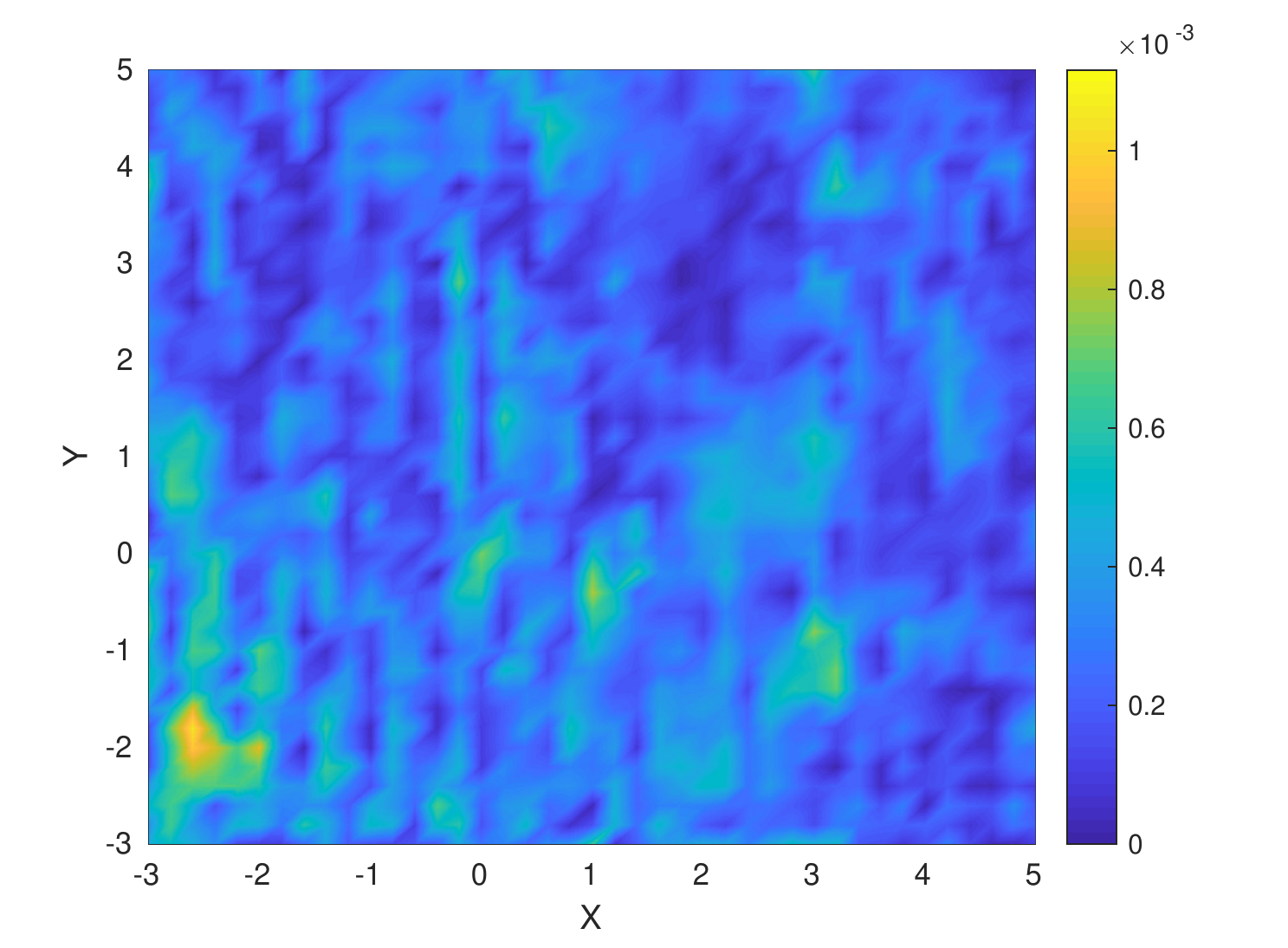}
				\par\end{centering}
		}\subfloat[Propagated data at $\alpha=0.5$\label{fig:Prop3}]{\begin{centering}
				\includegraphics[scale=0.5]{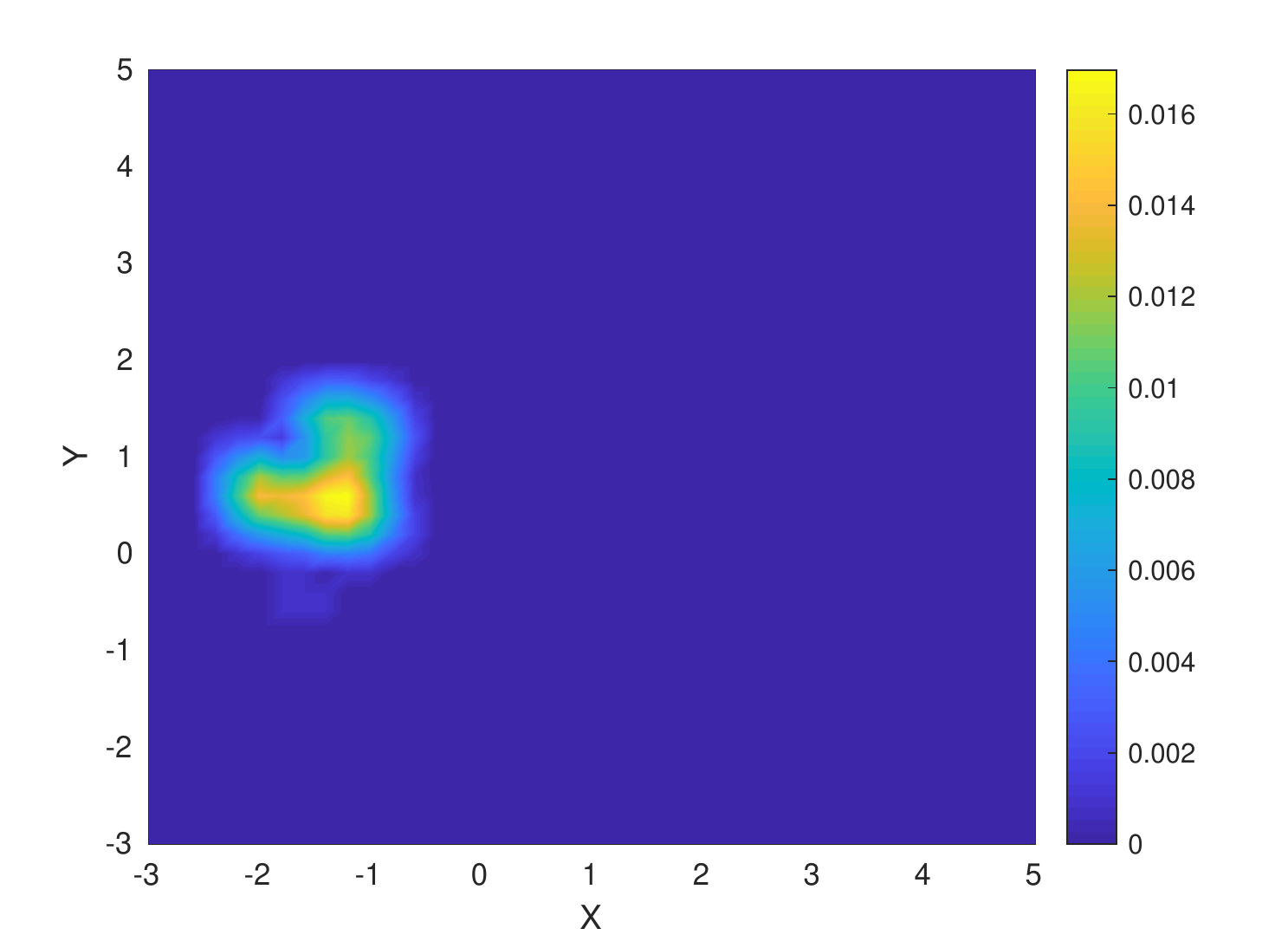}
				\par\end{centering}
		}
		\par\end{centering}
	\begin{centering}
		\subfloat[U-shaped piece of dry wood\label{fig:Walnut}]{\noindent \begin{centering}
				\includegraphics[scale=0.35]{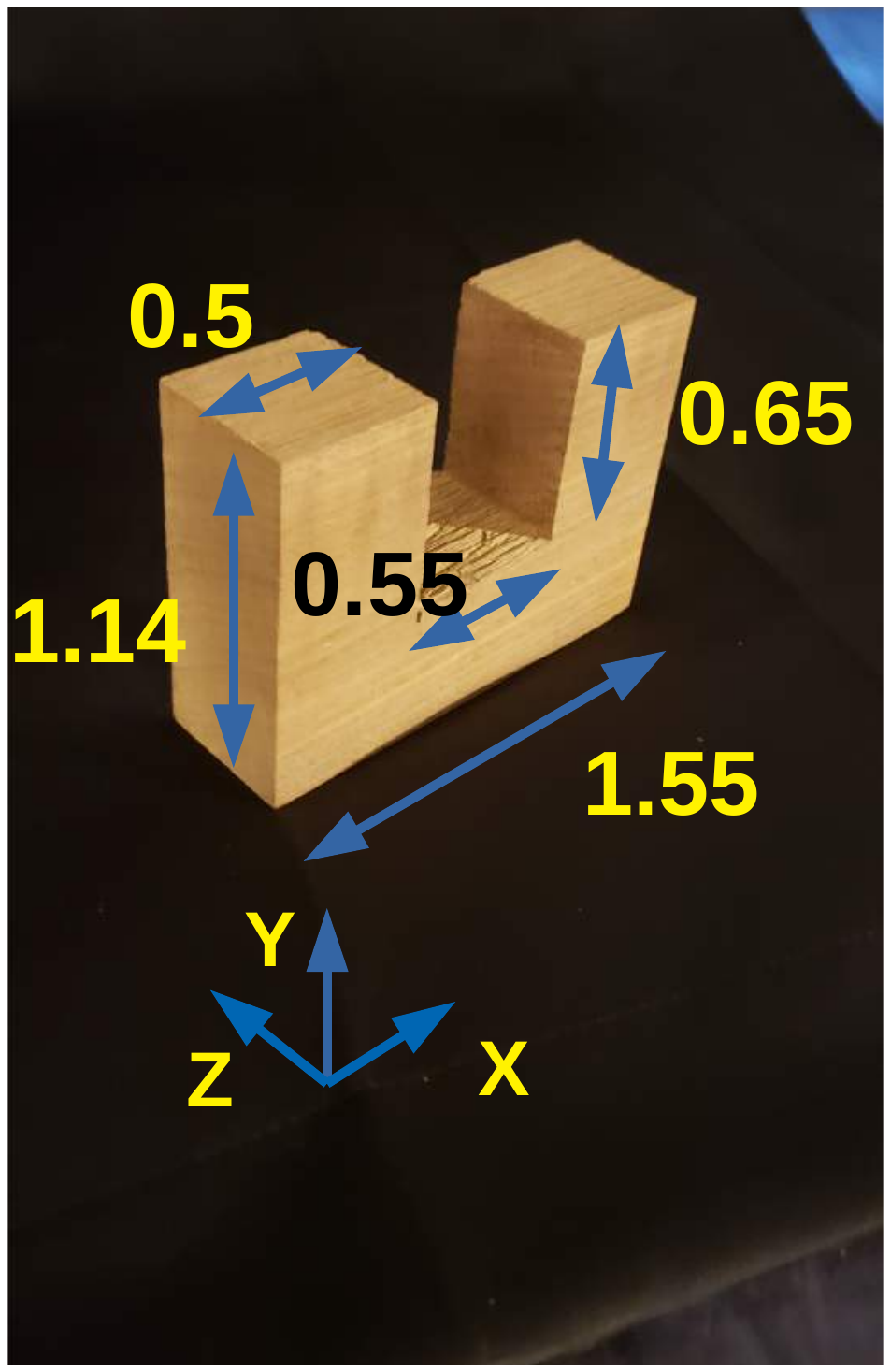}
				\par\end{centering}
		}\subfloat[Computed inclusion\label{fig:Comp3}]{\begin{centering}
				\includegraphics[scale=0.6]{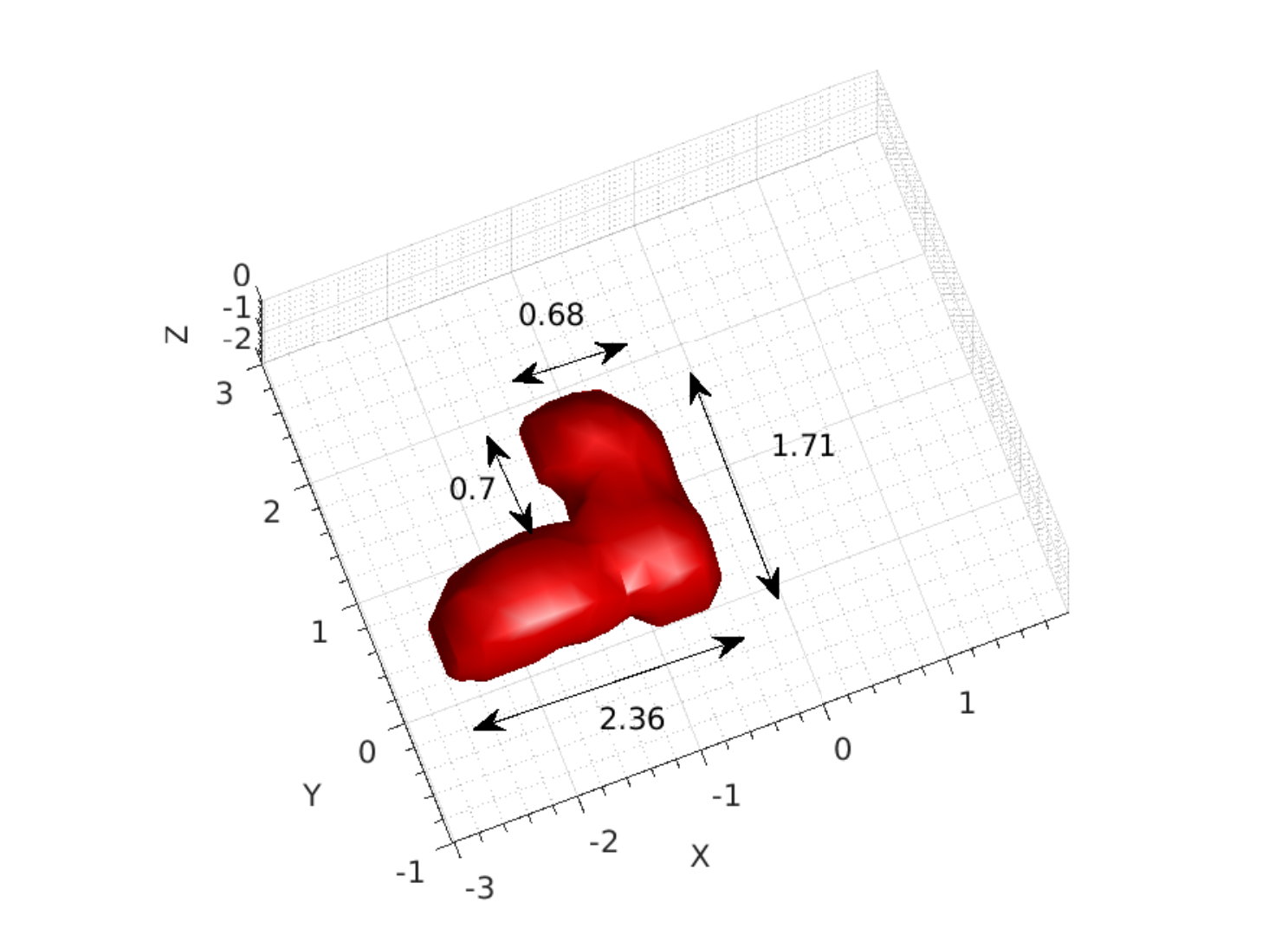}
				\par\end{centering}
		}
		\par\end{centering}
	\caption{Reconstruction results of Test 3 (U-shaped piece of dry wood). Note that the shape is non-convex, which is difficult to image. (a) Illustration
		of the absolute value of the raw far-field data; (b) Illustration of the absolute value of the near-field data after the data 
		propagation procedure; (c) Photo of the experimental object; (d) The
		computed image of (c). Note that the void is clearly seen which is difficult to image. Our axes on (d) are oriented differently from ones on (c) due to some technical problem of the imaging software. These axes are comparable.  All images are in the dimensionless
		variables.\label{fig:4}}
\end{figure}

Although Theorem \ref{thm:9} claims the global convergence of the gradient projection
method, we have successfully used the gradient descent method for the
minimization of the target functional $J_{h,\lambda }\left( V^{h}\right) $
of (\ref{eq:J}). Clearly, the gradient descent method is easier to implement
than the gradient projection method. Our success in working with the
gradient descent method is similar with the success in all previous
publications discussing the numerical studies of the convexification \cite%
{Khoa2019,Klibanov2017a,convIPnew,Klibanov2019,Klibanov2019b,Klibanov2019a,Klibhyp}%
. As to the value of the parameter $\lambda $ in $J_{h,\lambda }\left(
V^{h}\right) ,$ even though the analysis requires large values of $\lambda ,$
our numerical experience tells us that we can choose a moderate value of $%
\lambda :$
\begin{equation*}
\lambda =1.1.
\end{equation*}%
Again similar values of $\lambda \in \left[ 1,3\right] $ were chosen in
the above cited publications on the convexification.

As to the step size $\gamma $ of the gradient descent method, we start from $%
\gamma _{1}=10^{-1}$. For each iteration step $m\geq 1$, the following step
size $\gamma _{m}$ is reduced by the factor of 2 if the value of the
functional on the step $m$ exceeds its value of the previous step.
Otherwise, $\gamma _{m+1}=\gamma _{m}$. The minimization process is stopped
when either $\gamma _{m}<10^{-10}$ or $\left\vert J_{h,\lambda }\left(
V_{m}^{h}\right) -J_{h,\lambda }\left( V_{m-1}^{h}\right) \right\vert
<10^{-10}$. As to the gradient $J_{h,\lambda }^{\prime }$ of the discrete
functional $J_{h,\lambda }$, we apply the technique of Kronecker deltas (cf.
e.g. \cite{Kuzhuget2010}) to derive its explicit formula, which
significantly reduces the computational time. For brevity, we do not provide
this formula here.

After the minimization procedure is stopped, we obtain numerically the
coefficient of $c_{p,q,s}$, denoted by $\tilde{c}$. Our reconstructed
solution, denoted by $c_{\text{comp}}$, is concluded after we smooth $\tilde{%
c}$ by the standard filtering via the \texttt{smooth3} built-in function in
MATLAB. In fact, we find $c_{\text{comp}}$ by using $c_{\text{comp}}=\hat{%
\varrho}\text{smooth}\left( \left\vert \tilde{c}\right\vert \right) $, for
some $\hat{\varrho}>0$ depending on every single example. This step is
definitely similar to the smoothing procedure discussed in \eqref{eq:28} and
we do not repeat how to find $\hat{\varrho}$ here. We use this step to get 
better images.

\subsection{Reconstruction results}

Values of $\max \left( c_{\text{true}}\right) $ and $\max \left( c_{\text{%
comp}}\right) $ for all five tests are tabulated in Table \ref{table:2}.
Values of $\max \left( c_{\text{true}}\right) $ for all tests are mentioned
in subsection \ref{subsec:Experimental-examples}, which were used, are
published ones \cite{Kuzh,Table,Thanh2015}. More precisely, as to the
metallic targets of Example 1 (aluminum cylinder), Example 4 (metallic
letter ``A'') and Example 5 (metallic letter
``O''), it was numerically established that one can treat
metals as materials with large values of the dielectric constant in the
interval $c\in \left[ 10,30\right] ,$ see the formula (7.2) of \cite{Kuzh}.
As to the Example 2, the dielectric constant of the clear water for our
frequency range was directly measured in \cite{Thanh2015}, and it was 23.8:
see the first line of Table 2 of \cite{Thanh2015}. As to the Example 3 (an
U-shaped piece of a dry wood), the table of dielectric constants \cite{Table}
tells one that in the dielectric constant of a dry wood is $c\in \left[ 2,6%
\right] .$

\begin{table}
\begin{center}
\begin{tabular}{|c|c|c|c|c|c|}
\hline Example number & 1 & 2 & 3 & 4 & 5 \tabularnewline
\hline 
Object  & Metal Cylinder & Water & Wood & Metal letter \textquotedblleft A'' & 
Metal letter \textquotedblleft O'' \tabularnewline
\hline  
$\max \left( c_{\text{comp}}\right) $ & 18.72 & 23.29 & 6.56 & 15.01 & 16.25
\tabularnewline
\hline 
$c_{\text{true}}$ & [10,30] & 23.8 & $\left[ 2,6%
\right] $ & [10,30]  & [10,30] %
\tabularnewline
\hline 
Reference & \cite{Kuzh} & \cite{Thanh2015} & \cite{Table} & \cite{Kuzh} & \cite{Kuzh}%
\tabularnewline
\hline 
\end{tabular}
\end{center}
\caption{True $c_{\text{true}}$ and computed $\max \left( c_{\text{%
			comp}}\right) $ dielectric constants of Examples 1--5 of experimental data.
	True values were taken from: (a) Examples 1, 4, 5: formula (7.2) of \cite{Kuzh}%
	, (b)  Example 2 (clear water) \cite{Thanh2015}, (c) Example 3 \cite{Table}. \label{table:2}}
\end{table}

\begin{figure}[H]
	\begin{centering}
		\subfloat[Raw data at $\alpha=0.2$\label{fig:Raw4}]{\begin{centering}
				\includegraphics[scale=0.5]{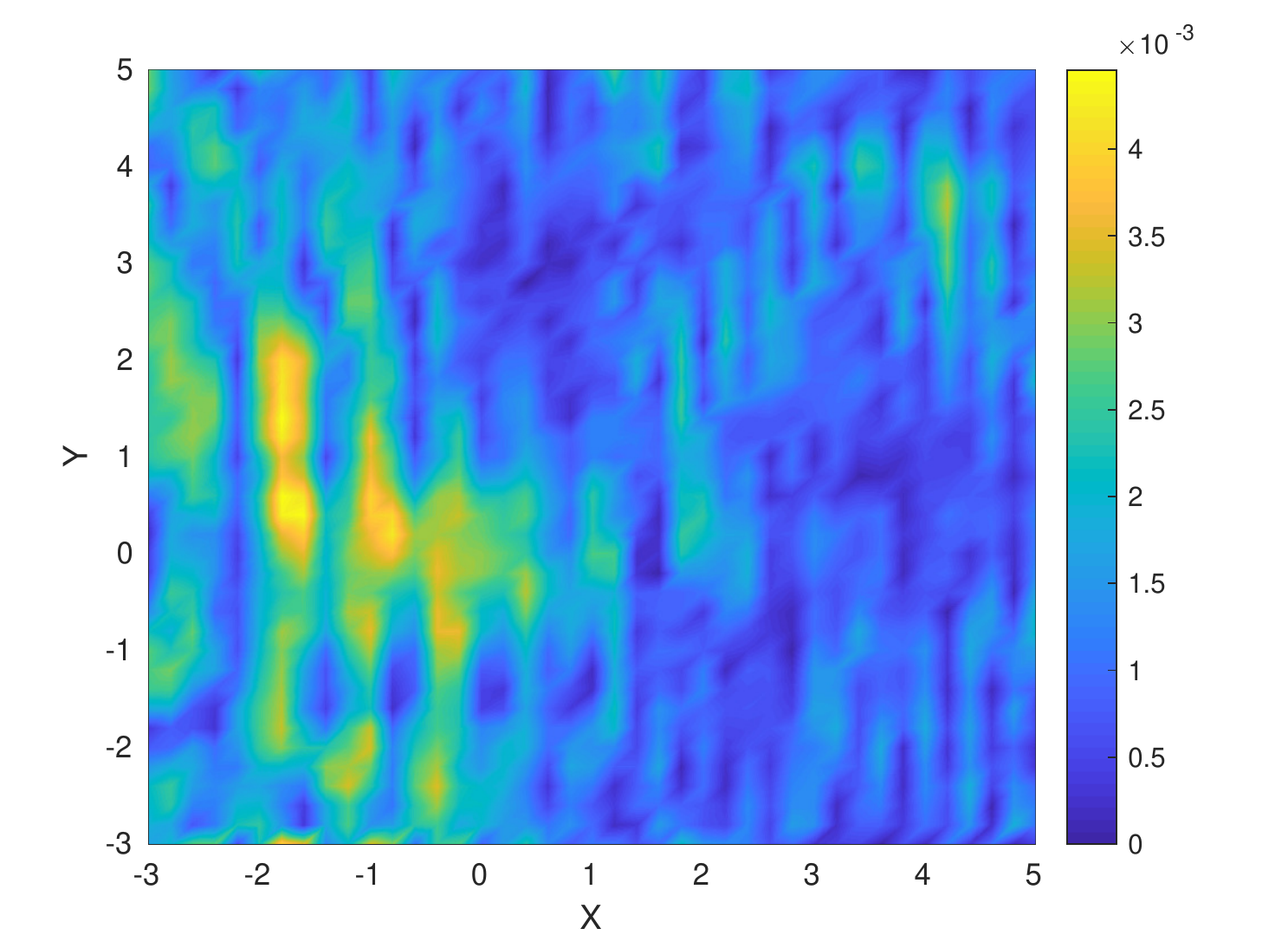}
				\par\end{centering}
		}\subfloat[Propagated data at $\alpha=0.2$\label{fig:Prop4}]{\begin{centering}
				\includegraphics[scale=0.5]{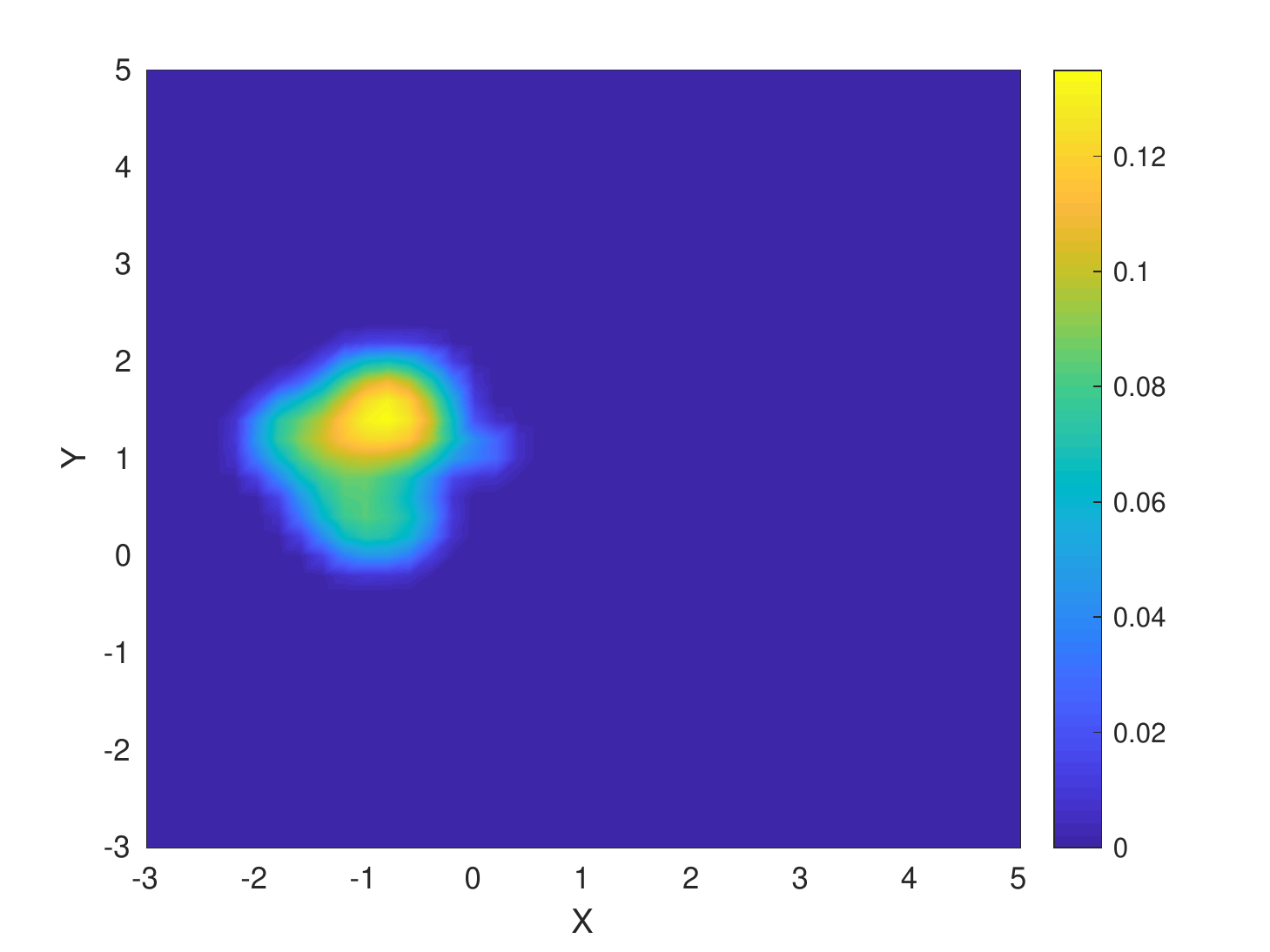}
				\par\end{centering}
		}
		\par\end{centering}
	\begin{centering}
		\subfloat[Metallic letter ``A''. Blind test. The shape is non-convex, which is difficult to image.\label{fig:A}]{\noindent \begin{centering}
				\includegraphics[scale=0.3]{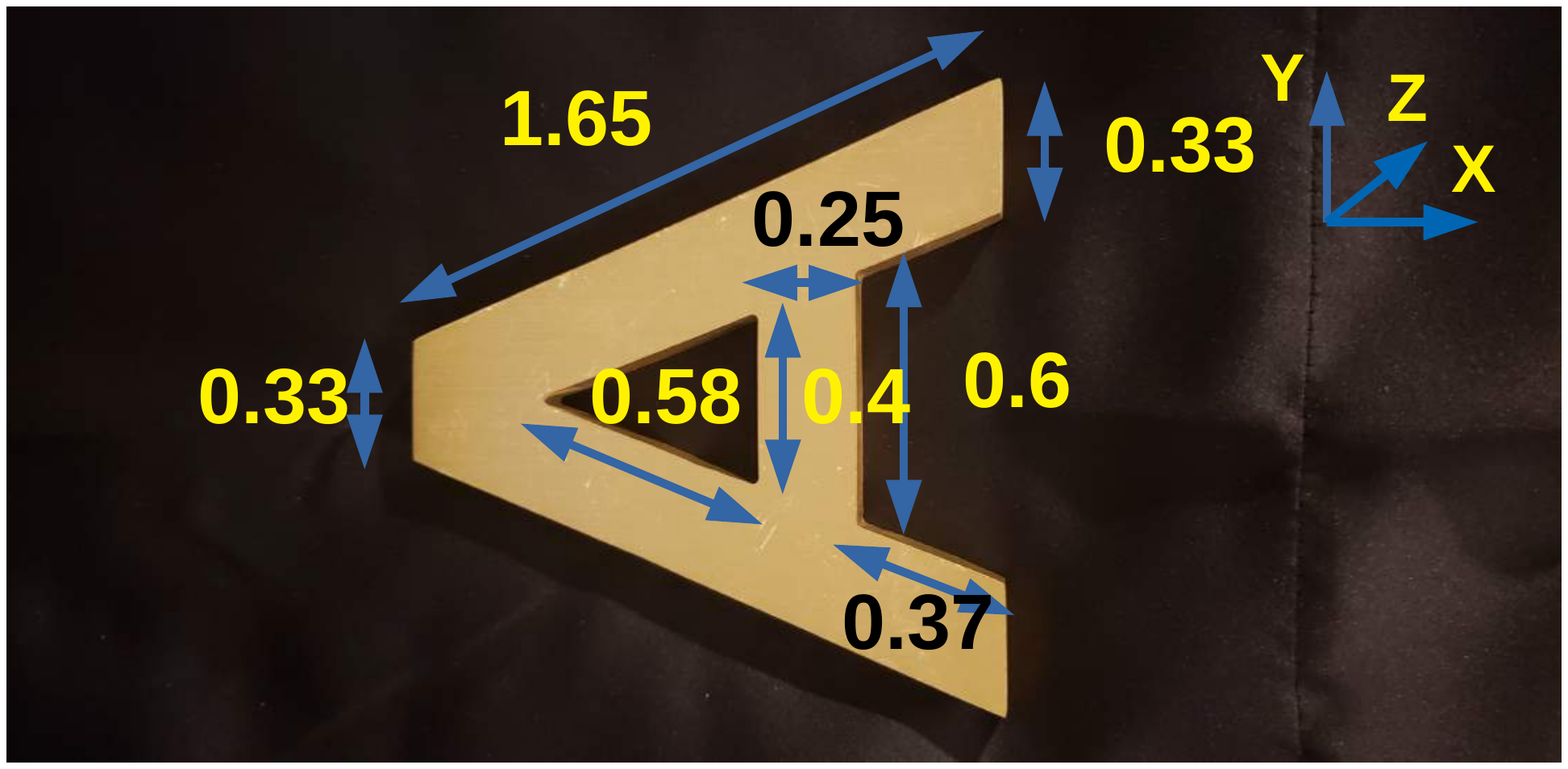}
				\par\end{centering}
		}\subfloat[Computed inclusion\label{fig:Comp4}]{\begin{centering}
				\includegraphics[scale=0.6]{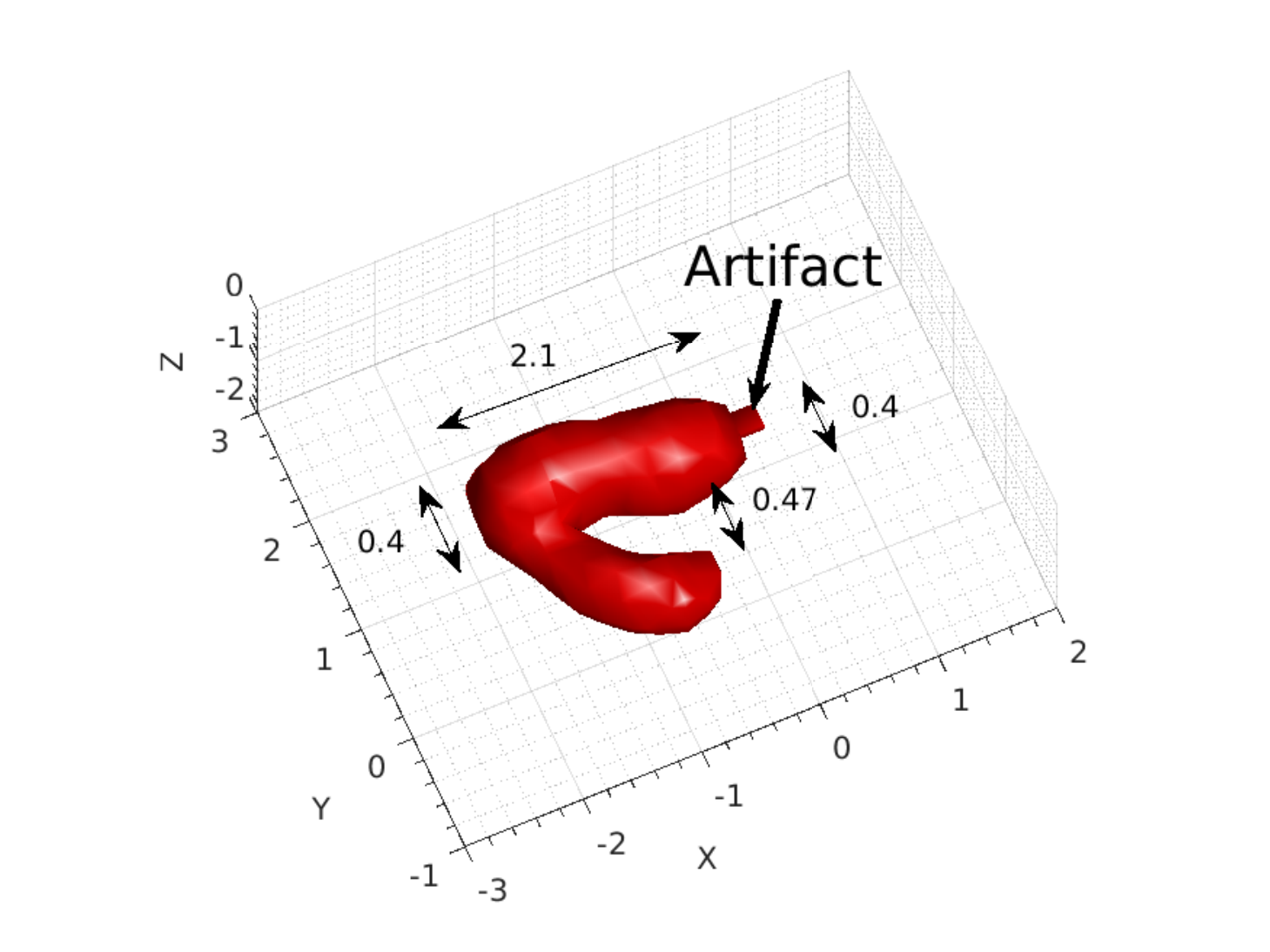}
				\par\end{centering}
		}
		\par\end{centering}
	\caption{Reconstruction results of Test 4 (metallic letter ``A''). This is a blind test. (a) Illustration
		of the absolute value of the raw far-field data; (b) Illustration of the absolute value of the near-field data after the data 
		propagation procedure; (d) The computed image of (c). Note that the void is clearly seen, which is challenging to image. Also, sizes of the imaged target are close to the true ones. The strip of ``A'' is not seen since its width is 2.5 cm, which is less than the wavelength of 10.4 cm we have used with $k=9.55$. All images are in the dimensionless
		variables. \label{fig:5}}
\end{figure}

\begin{figure}[H]
	\begin{centering}
		\subfloat[Raw data at $\alpha=0.6$\label{fig:Raw5}]{\begin{centering}
				\includegraphics[scale=0.5]{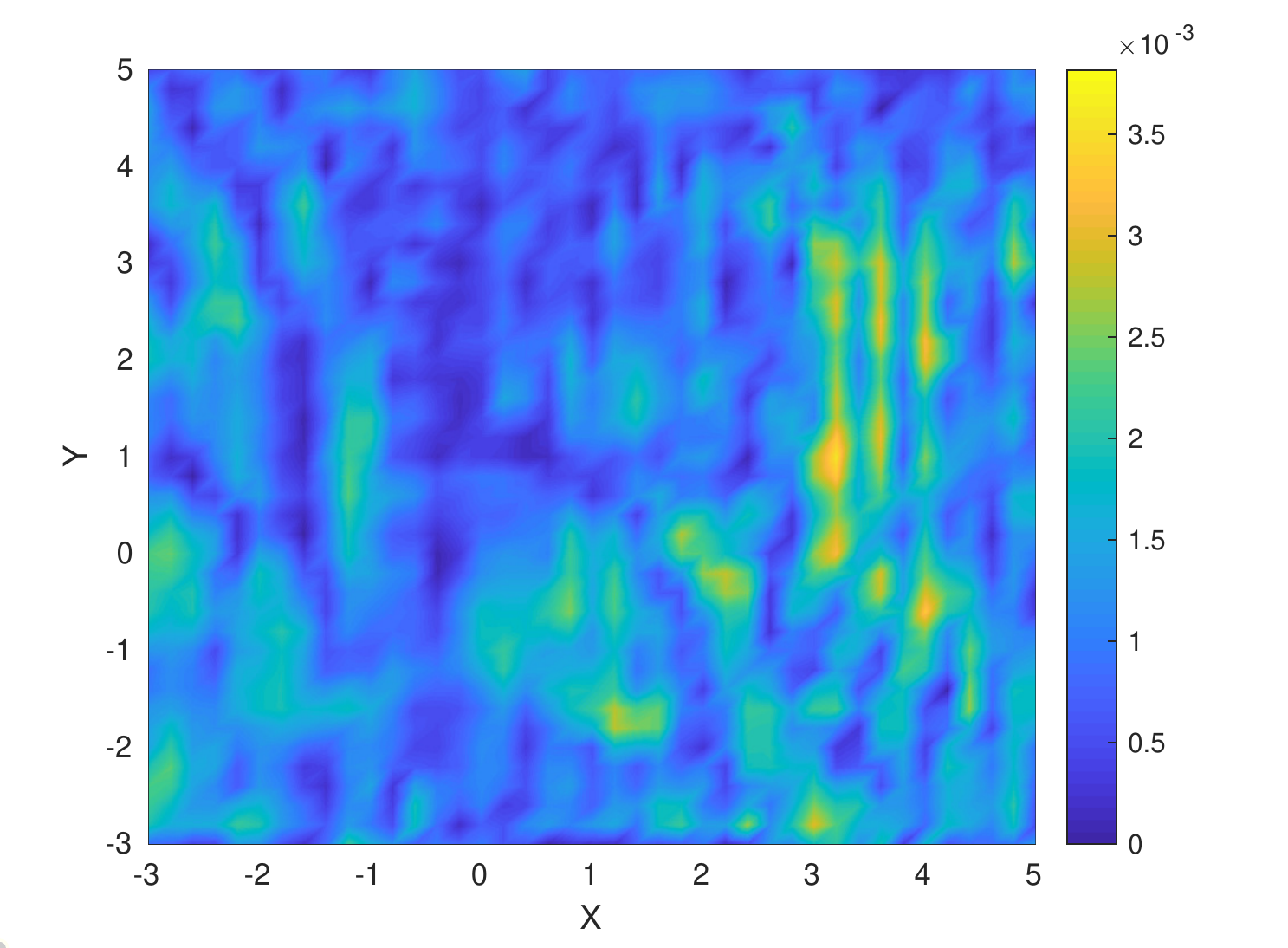}
				\par\end{centering}
		}\subfloat[Propagated data at $\alpha=0.6$\label{fig:Prop5}]{\begin{centering}
				\includegraphics[scale=0.5]{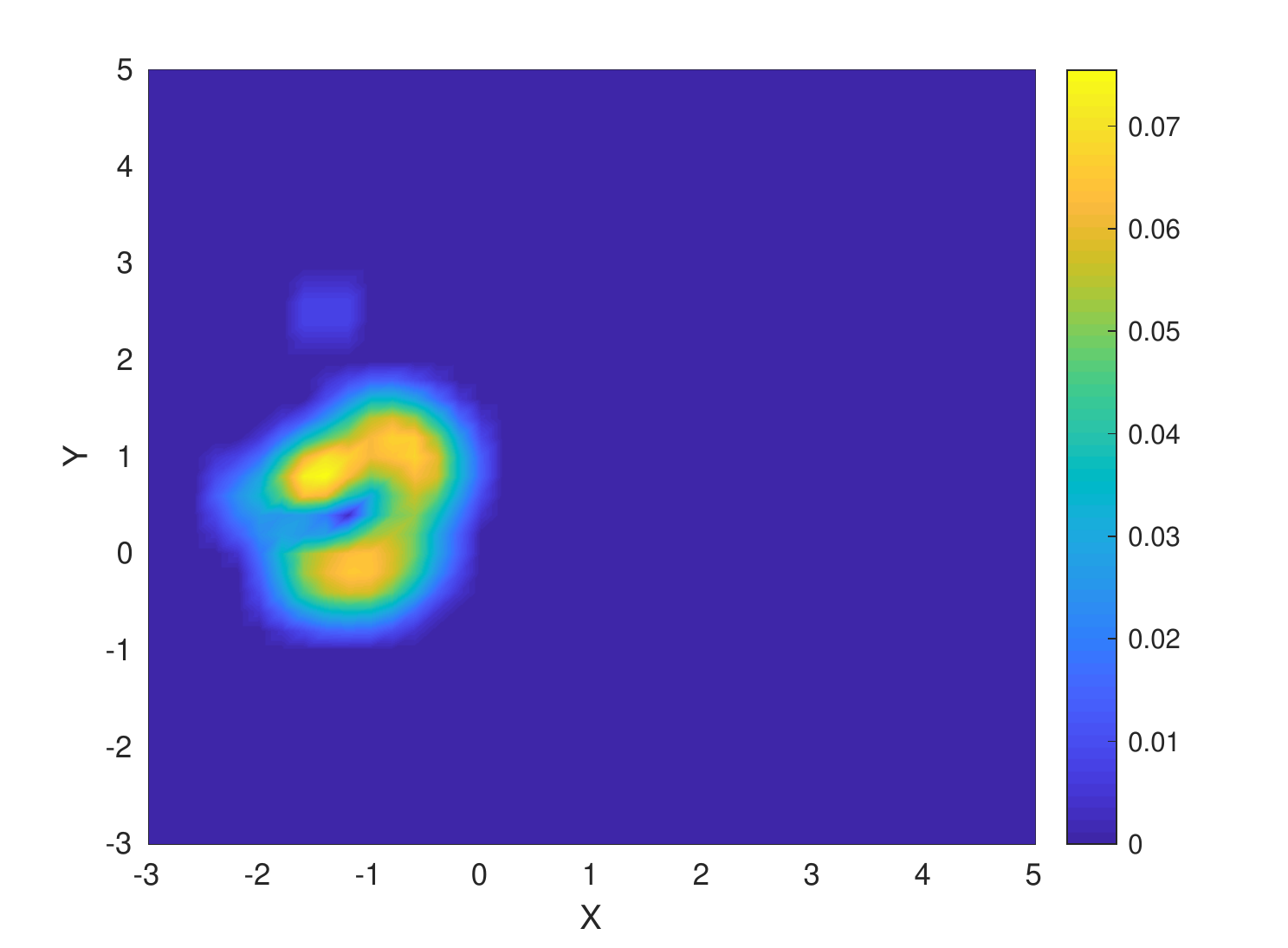}
				\par\end{centering}
		}
		\par\end{centering}
	\begin{centering}
		\subfloat[Metallic letter ``O''. Blind test. The shape is non-convex, which is challenging to image.\label{fig:O}]{\noindent \begin{centering}
				\includegraphics[scale=0.3]{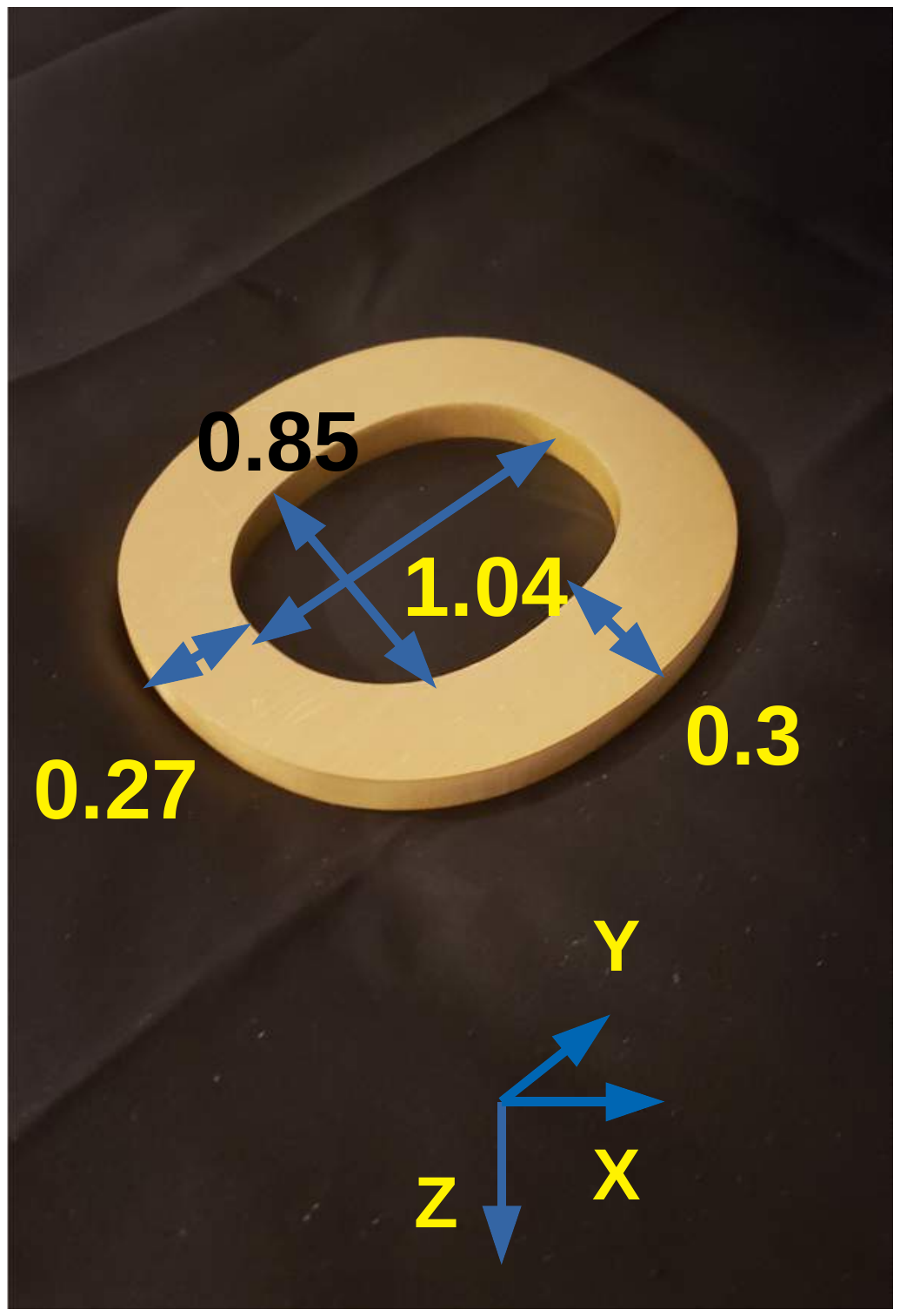}
				\par\end{centering}
		}\subfloat[Computed inclusion\label{fig:Comp5}]{\begin{centering}
				\includegraphics[scale=0.6]{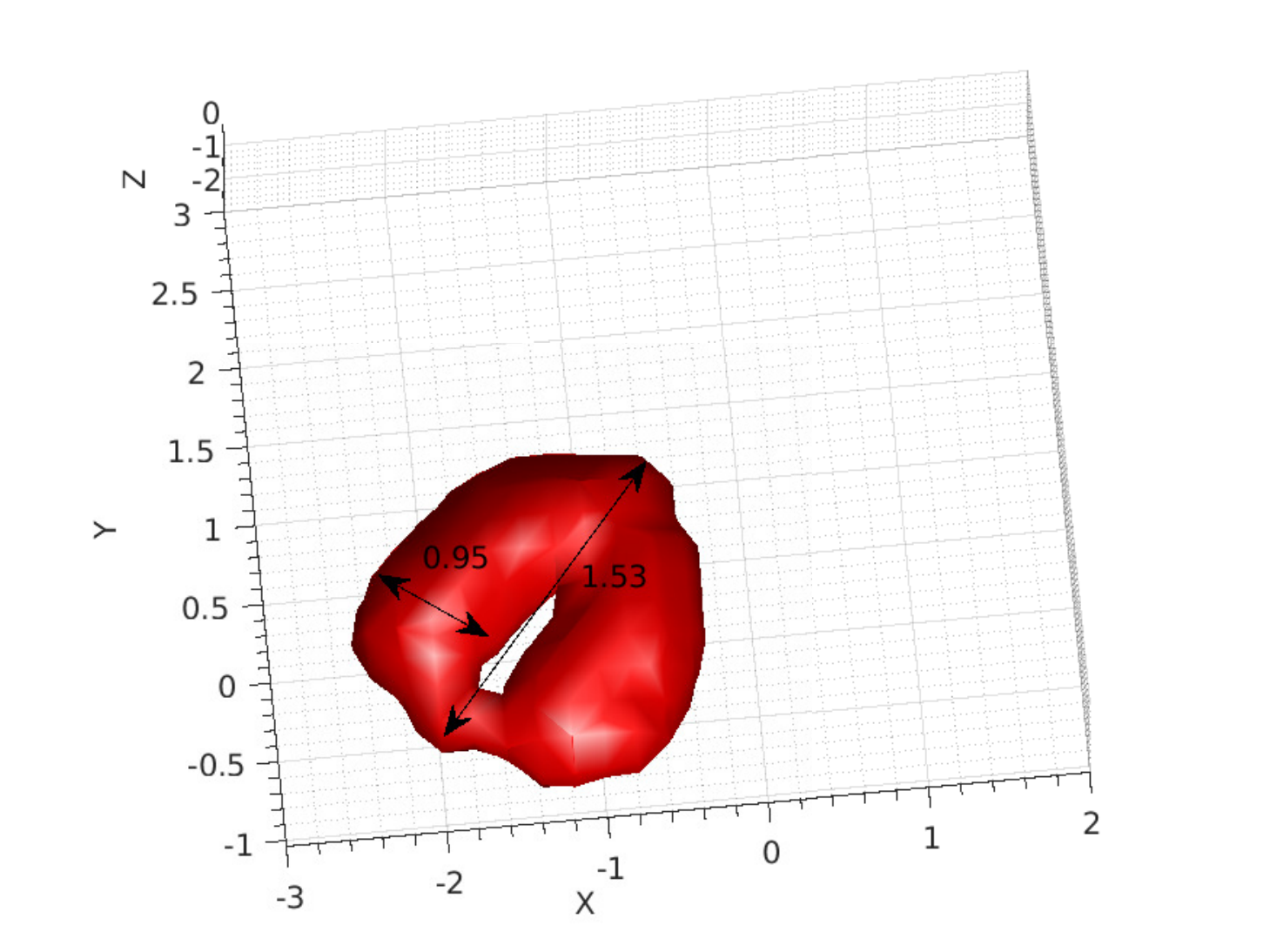}
				\par\end{centering}
		}
		\par\end{centering}
	\caption{Reconstruction results of Test 5 (metallic letter ``O''). This is a blind test. (a) Illustration
		of the absolute value of the raw far-field data; (b) Illustration of the absolute value of the near-field data after the data 
		propagation procedure; (c) Photo of the experimental object; (d) The
		computed image of (c). Note that the void is clearly seen, which is not easy to image. All images are in the  dimensionless
		variables. \label{fig:6}}
\end{figure}

Figures \ref{fig:Raw1}, \ref{fig:Raw2}, \ref{fig:Raw3}, \ref{fig:Raw4} and %
\ref{fig:Raw5} show how \textquotedblleft bad\textquotedblright\ the
far-field data look like. It is clear from these figures that something
should be done to the data to have a proper inversion. On the other hand,
one can see the good shapes of the corresponding images after the data
propagation procedure; see Figures \ref{fig:Prop1}, \ref{fig:Prop2}, \ref%
{fig:Prop3}, \ref{fig:Prop4} and \ref{fig:Prop5}. For every test, we
deliberately show the 2D illustrations (raw and propagated) of the data at a
specific point source, where the images of the propagated data and the
computed inclusion are congruent with each other. 

3D images of computed inclusions are depicted by using the \texttt{isosurface%
} function in MATLAB with the associated \texttt{isovalue} being 10\% of the
maximal value; see Figures \ref{fig:Comp1}, \ref{fig:Comp2}, \ref{fig:Comp3}%
, \ref{fig:Comp4} and \ref{fig:Comp5}. The most challenging targets to image
were: (1) The U-shaped piece of dry wood, see Figure \ref{fig:4}, (2) The metallic
letter ``A'', see Figure \ref{fig:5}, and (3) the metallic letter
``O'', see Figure \ref{fig:6}. This is because these targets have the
most complicated geometries. Nevertheless, we are still able to see their
characteristic shapes in the images of computed inclusion. 

Most notably, one can see voids in imaged letters ``A'' and
``O''. The latter is usually difficult to achieve. The
``strip'' of the letter ``A'' is not imaged
since its width was 2.5 cm, which is less than the used wavelength of 10.4 cm
with $k=9.55$. Another interesting observation is that we can even see the cap
on the bottle of water on Figure \ref{fig:Comp2}. 

We also find that the lengths of parts of true and computed inclusions are
quite compatible with each others. Note that even though the computed
inclusions here are slightly larger (just a few centimeters) than the true
ones, it is still useful in detection and identification of land mines and further in the
mine-clearing operations. In fact, having information of smaller sizes is
rather dangerous. Hence, we conclude that the dimensions of the computed
inclusions are acceptable.

Finally, we can accurately obtain approximations of the
dielectric constants. Aside from the dielectric constant of metallic
targets, we notice from Table \ref{table:2} that the relative errors obtained
for the bottle with water and for the wooden target are 2.14\% and 9.33\%,
respectively.

\section{Summary}

We have developed \ a new version of the globally convergent convexification
method for the case of a 3D CIP for the Helmholtz equation. In our case, the
point source is moving along an interval of a straight line and the
frequency is fixed. For each position of the point source we measure one
component of the backscattering electric wave field at a part of a plane.
Thus, our data depend on three variables, which means that they are non
overdetermined ones. We use the partial finite differences and construct a
weighted cost functional with the Carleman Weight Function in it. The use of
partial finite differences enables us to avoid the use of the penalty
regularization term. The latter is the major analytical novelty here. We
prove that our functional is strictly convex on a finite set of an arbitrary
size. This theorem leads to the theorem about the \emph{global convergence}
to the exact solution of the gradient projection method of the minimization
of this functional, as long as the level of noise in the data tends to zero.
The global convergence property is the most important feature of the
convexification method.

We have tested our method numerically on backscattering experimentally
collected data. Our testing reveals that we can accurately image both
dielectric constants and shapes of targets of interest. Including even
rather complicated geometries. This is an advantage compared with the
previously considered version of the convexification in which the point
source was fixed and the frequency was varied. Indeed, while in the latter
case the dielectric constants were computed accurately, shapes were not accurately
imaged; see, e.g., \cite{Klibanov2019b} for the case of experimental data. 



\end{document}